\documentclass[11pt]{amsart}
\usepackage{amsmath,amssymb,amsfonts}
\usepackage{enumerate}
\usepackage{amscd, amssymb, latexsym, amsmath, amscd}
\usepackage[all]{xy}
\usepackage{pb-diagram}
\usepackage{multicol,lipsum}

\newtheorem{thm}[equation]{Theorem}
\newtheorem{prop}[equation]{Proposition}
\newtheorem{cor}[equation]{Corollary}

\newtheorem{lemma}[equation]{Lemma}

\numberwithin{equation}{section}

\newcommand{\Q}{\mathbb Q}
\newcommand{\Z}{\mathbb Z}
\newcommand{\R}{\mathbb R}
\newcommand{\C}{\mathbb C}
\newcommand{\B}{\mathbb B}

\newcommand{\G}{\mathfrak G}

\newcommand{\A}{\mathbb A}

\newcommand{\N}{\mathbb N}

\def\Hom{{\rm Hom}}

\def\G{{\rm G}}
\def\SL{{\rm SL}}

\def\GSp{{\rm GSp}}

\def\Sp{{\rm Sp}}

\def\U{{\rm U}}

\def\GL{{\rm GL}}
\def\PGL{{\rm PGL}}

\def\SO{{\rm SO}}

\def\Sp{{\rm Sp}}
\def\Mp{{\rm Mp}}

\usepackage[all]{xy}

\def\A{{\mathbb A}}

\def\R{{\mathbb R}}

\def\Z{{\mathbb Z}}
\def\T{{\mathbb T}}

\def\C{{\bf C}}

\def\C{{\mathbb C}}

\def\G{{\mathbb G}}

\title[The Langlands-Weissman  Program]{The Langlands-Weissman  Program \\for Brylinski-Deligne extensions} 
 \author{Wee Teck Gan and Fan Gao}

\address{  Department of Mathematics, National University of Singapore, 10 Lower Kent Ridge Road
Singapore 119076} \email{matgwt@nus.edu.sg}
\email{gaofan@nus.edu.sg}

\usepackage{picture}
\begin{document}
\begin{abstract}
We describe an evolving and conjectural extension of the Langlands program 
for a class of nonlinear covering groups of algebraic origin studied by Brylinski-Deligne.  In particular, we describe the construction of an L-group extension of such a covering group (over a split reductive group) due to Weissman, study some of its properties and discuss a variant of it. Using this L-group extension, we describe a local Langlands correspondence for covering (split) tori and unramified genuine representations, using work of Savin, McNamara, Weissman and W.W. Li. Finally, we define the notion of automorphic (partial) L-functions attached to genuine automorphic representations of  the BD  covering groups. 
\end{abstract}

 
\maketitle

\section{\bf Introduction}

One of the goals of the local  Langlands program is to provide an arithmetic classification of the set of isomorphism classes of irreducible  representations of a locally compact group $G = \mathbb{G}(F)$, where $\mathbb{G}$ is a connected reductive group over a local field $F$.   Analogously, if $k$ is a number field with ring of adeles $\A$, the global Langlands program postulates a classification of automorphic representations of $\mathbb{G}(\A)$ in terms of Galois representations.  In this proposed arithmetic classification, which has been realised in several important instances, a key role is played by the L-group ${^L}\mathbb{G}$ of $\mathbb{G}$. This key notion was introduced by Langlands in his re-interpretation of the Satake isomorphism in the theory of spherical functions and used by him to introduce the notion of {\em automorphic L-functions}.
\vskip 5pt

\subsection{\bf Covering groups}
The theory of the L-group is confined to the case when $\G$ is a connected reductive linear algebraic group. 
On the other hand, since Steinberg's beautiful paper \cite{S}, the structure theory of nonlinear covering groups of $G$ (i.e. topological central extensions of $G$ by finite groups) have been investigated by many mathematicians, notably Moore \cite{Mo}, Matsumoto \cite{Ma}, Deodhar \cite{De}, Deligne \cite{D}, Prasad-Raghunathan \cite{PR1,PR2,PR3}, and its relation to the reciprocity laws of abelian class field theory has been noted.  In addition,  nonlinear covering groups of $G$  have repeatedly made their appearance in representation theory and the theory of automorphic forms. This goes way back to  Jacobi's construction of his theta function, a holomorphic modular form of weight $1/2$, and a more recent instance is the Shimura correspondence between integral and half integral weight modular forms and the work of Kubota \cite{Ku}. Both these examples concern automorphic forms and representations of the metaplectic group $\Mp_2(F)$, which is a nonlinear double cover of $\SL_2(F) = \Sp_2(F)$. As another example, the well-known Weil representation of $\Mp_{2n}(F)$ gives a representation theoretic  incarnation of theta functions and has been a very useful tool in the construction of automorphic forms.  Finally, much of Harish-Chandra's theory of local harmonic analysis and Langlands' theory of Eisenstein series  continue to hold for such nonlinear covering groups (see \cite{MW} and \cite{L4}).

\vskip 5pt

It is thus natural to wonder if the framework of the Langlands program can be extended to encompass the representation theory and the theory of automorphic forms of covering groups. There have been many attempts towards this end, such as Flicker \cite{F}, Kazhdan-Patterson \cite{KP1,KP2}, Flicker-Kazhdan \cite{FK}, Adams \cite{A1, A2}, Savin \cite{Sa} among others.  
However, these attempts have tended to focus on the treatment of specific families of examples rather than a general theory.  This is understandable, for what is lacking is a structure theory which is sufficiently functorial.
For example, the classification of nonlinear covering groups given in \cite{Mo, De, PR1, PR2} is given only when $\G$ is simply-connected and isotropic, in which case  a universal cover exists.
\vskip 5pt

\subsection{\bf Brylinski-Deligne theory}
A functorial structure theory was finally developed by  Brylinski and Deligne   in their 2002 paper \cite{BD}. 
More precisely, Brylinski-Deligne considered the category of multiplicative $K_2$-torsors on a connected reductive group $\G$ over $F$; these are extensions of $\G$ by the sheaf $\mathbb{K}_2$ of Quillen's $K_2$ group  in the category of sheaves of groups on the big Zariski site of $Spec(F)$:
\[  \begin{CD}
\mathbb{K}_2 @>>> \overline{\mathbb{G}} @>>> \mathbb{G}.  \end{CD} \]

\noindent In other words, Brylinski-Deligne started with  an extension problem in the world of algebraic geometry.  Some highlights of [BD]  include:
\vskip 5pt

\begin{itemize}
\item an elegant and functorial classification of this category in terms of enhanced root theoretic data, much like the  classification of split connected  reductive groups by their root data. 
\vskip 5pt

\item the description of a functor from the category of multiplicative $K_2$-torsors $\overline{\mathbb{G}}$ on $\G$ (together with an integer $n$ such that $\# \mu_n(F) = n$, which determines the degree of the covering)  to the category of topological central extensions $\overline{G}$ of $G$:
\[ \begin{CD}
 \mu_n @>>> \overline{G}  @>>> G. \end{CD} \]
These topological central extensions may be considered of ``algebraic origin" and can be constructed using cocycles which are essentially algebraic in nature. 
\vskip 5pt

\item though this construction does not exhaust all topological central extensions, it  captures a sufficiently large class of such extensions, and essentially all interesting examples which have been investigated so far; for example, it captures all such coverings of $G$ when  $\G$ is split and simply-connected.
\end{itemize}
\vskip 5pt

\noindent 
We shall give a more detailed discussion of the salient features of the Brylinski-Deligne theory in \S 2 and \S 3. 
Hence, the paper \cite{BD} provides a structure theory which is essentially algebraic and categorical, and  may be perceived as a natural extension of Steinberg's original treatment \cite{S} from the split simply connected case to general reductive groups. 
\vskip 5pt

\subsection{\bf Dual and L-groups.}
One should expect that  such a natural structure theory would elucidate the study of representations and automorphic forms of the BD covering groups $\overline{G}$. The first person to fully appreciate this (that is, besides Brylinski and Deligne themselves) is probably our colleague  M. Weissman.
In a series of papers \cite{W1, W2, HW}, Weissman systematically exploited the Brylinski-Deligne theory to study the representation theory of covering tori, the unramified representations  and the depth zero representations. 
This was followed by the work of several people who discovered  a ``Langlands dual group" $\overline{G}^{\vee}$ for a BD covering group $\overline{G}$ (with $\G$ split) from different considerations. These include the work of Finkelberg-Lysenko \cite{FL} and  Reich \cite{Re}  in the framework of the geometric Langland program and the work of McNamara \cite{Mc2, Mc3} who established a Satake isomorphism and interpreted it in terms of  the dual group $\overline{G}^{\vee}$. The dual group $\overline{G}^{\vee}$ was constructed by making a metaplectic modification of the root datum of $\G$. 
\vskip 5pt

In another paper \cite{W3}, Weissman built upon \cite{Mc2}  and gave a construction of the ``L-group"  ${^L}\overline{G}$ of a {\em split}  BD covering group $\overline{G}$. The construction in \cite{W3}  is quite involved, and couched in the language of Hopf algebras. Moreover,  with hindsight, it gives  the correct notion only for a subclass of BD covering groups. 
 In a foundational paper \cite{W7} under preparation, Weissman give a simpler and completely general revised construction of the L-group for an arbitrary BD covering group (not necessarily split), using the language of \'etale gerbes, thus laying the groundwork for an extension of the Langlands program to the setting of BD covering groups. 
\vskip 5pt

As Weissman is the first person to make use of the full power of the Brylinski-Deligne structure theory for the purpose of representation theory,
we shall call this evolving area  the Langlands-Weissman program for BD extensions.

 \vskip 10pt

\subsection{\bf Define ``is"}
 We  shall describe in  \S 4 Weissman's construction of the L-group of $\overline{G}$  for split $\G$ (given in the letter  \cite{W4}, where one could be more down-to-earth and avoid the notion of gerbes).
At this point, let us note that since $\G$ is split, one is inclined  to simply take ${^L}\overline{G}$ as the direct product $\overline{G}^{\vee} \times W_F$, where $W_F$ denotes the Weil group of $F$. At least, this is what one is conditioned to do by the theory of L-groups for linear reductive groups. However, Weissman realised that such an approach would be overly naive.

\vskip 5pt

Indeed, the key insight  of \cite{W3} is that the construction of the L-group of a BD covering group   should be the functorial construction of an extension
\begin{equation}  \label{E:LG}   \begin{CD}
 \overline{G}^{\vee} @>>>  {^L}\overline{G}  @>>>  W_F, \end{CD} \end{equation}
 and an L-parameter for $\overline{G}$ should be a splitting of this short exact sequence. 
 The point is that, {\em even if ${^L}\overline{G}$ is isomorphic to the direct product $\overline{G}^{\vee} \times W_F$},  it is not supposed to be equipped with a canonical isomorphism to  $\overline{G}^{\vee} \times W_F$. 
This reflects the fact that there is no  canonical irreducible genuine representation of $\overline{G}$, and hence there should not be any canonical L-parameter. Hence it would not be right to say that the L-group of $\overline{G}$ ``is" the direct product $\overline{G}^{\vee}$ with $W_F$.
 
\vskip 5pt

\subsection{\bf Goals of this paper}
Against this backdrop, the purpose of this paper is to supplement the viewpoint of   \cite{W3,W4,W7}  concerning the L-group ${^L}\overline{G}$ in several ways. We summarise our results here:
\vskip 5pt

\begin{itemize}
\item[(i)] Firstly, we point out that, though the L-group extension of a split BD covering group constructed in \cite{W3, W4,W7}  is a split extension, the L-group is {\em not} isomorphic to the direct product $\overline{G}^{\vee} \times W_F$ in general. This phenomenon can be seen already in the following simple family of BD covering groups
\[  \overline{G}_{\eta}  =  (\GL_2(F) \times \mu_2) / i_{\eta} (F^{\times}) \]
where $\eta \in F^{\times}/ F^{\times 2}$ and 
\[  i_{\eta}(t)  = ( t, (\eta, t)_2), \]
with $(-,-)_2$ denoting the quadratic Hilbert symbol. Observe that the first projection defines a topological central extension
\[  \begin{CD} 
\mu_2 @>>>  \overline{G}_{\eta} @>>>  \PGL_2(F).  \end{CD} \]
Then it turns out that $\overline{G}_{\eta}^{\vee}   \cong \SL_2(\C)$ and  
\[  {^L}\overline{G}_{\eta}  \cong  \SL_2(\C)  \rtimes_{\eta} W_F \]
where the action of $W_F$  on $\SL_2(\C)$ is by the conjugation action  via the map
$W_F  \longrightarrow \GL_2(\C)$ given by
\[  w \mapsto \left(
\begin{array}{cc}
\chi_{\eta}(w)  & \\
&  1  \end{array}  \right),  \]
with $\chi_{\eta}$ the quadratic character of $W_F$ determined by $\eta$. 
\vskip 5pt

This family of BD covering groups is quite instructive, as it illustrates several interesting phenomena. For example, one can show that the covering splits over the hyperspecial maximal compact subgroup $\PGL_2(\mathcal{O}_F)$ if and only if $\eta \in \mathcal{O}_F^{\times}$.  This seems to contradict \cite{Mc1, Mc3}  where it was claimed that a BD covering for a split $\G$ is always split over a hyperspecial maximal compact subgroup.  

\vskip 5pt

\item[(ii)]
Secondly, we would like to argue   that the L-group of a split BD covering group should be (the isomorphic class of)  the extension (\ref{E:LG}), {\em together with a finite set of distinguished splittings} (which we will define and which give rise to isomorphisms  ${^L}\overline{G} \cong \overline{G}^{\vee}  \times W_F$ if they exist).  If the degree of the covering is $n= 1$, so that $\overline{G} = G$, then the set of distinguished 
splittings is a singleton, so that ${^L}\overline{G}$ ``is"  $\overline{G}^{\vee} \times W_F$ in this case. 
 For a subclass of BD covering groups, we show that  the L-group short exact sequence (\ref{E:LG}) has such a family of distinguished splittings; this completes the results of \cite{W3} where the case for odd $n$ and  $n=2$ was treated.   
\vskip 5pt

\item[(iii)]   Thirdly,  we show that the distinguished splittings of ${^L}\overline{G}$ are in natural bijection with a family of distinguished genuine characters of the covering torus $\overline{T}$ (where $T$ is a maximal $F$-split torus of $G$). 
These distinguished genuine characters of $\overline{T}$ are `` as close to being trivial characters as possible", and are invariant under the natural Weyl group action for {\em certain} $\overline{G}$. Thus,  they serve as natural base-points for the definition of principal series representations, extending the results of Savin \cite{Sa} to general $n$ and these $\overline{G}$. We also give an explicit construction of such distinguished genuine characters of $\overline{T}$, using the Weil index and the $n$-th Hilbert symbol.  
\vskip 5pt

We point out, however, that for general $\overline{G}$, there may not exist distinguished splittings of ${^L}\overline{G}$ or Weyl-invariant genuine characters of $\overline{T}$. For example, for the covering groups 
$\overline{G}_{\eta}$ discussed in (i), Weyl-invariant genuine characters exist if and only if $(\eta, -1)_2  = 1$. 
\vskip 5pt

When $\G = \T$ is a split torus, the above passage between distinguished splittings of ${^L}\overline{T}$  and distinguished genuine characters of $\overline{T}$ extends to give a more explicit treatment of the local Langlands correspondence for covering tori shown in \cite{W1}. 

\vskip 5pt

\item[(iv)] Fourthly, we would like to suggest a slightly different treatment of the L-group given in \cite{W3,W4,W7}, by treating several closely related BD covering groups together. For example, the representation theory of all the groups $\overline{G}_{\eta}$  in (i) can clearly be treated together in terms of the representation theory of $\GL_2(F)$.  Extending this instructive example, one can slightly enhance Weissman's construction in \cite{W4, W7} to give an enlarged  L-group ${^L}\overline{G}^{\#}$:
\begin{equation}  \label{E:LG2}   \begin{CD}
 \overline{G}^{\vee} @>>>  {^L}\overline{G}^{\#}  @>>>  W_F \times (T_{Q,n}^{sc})^{\vee}[n] 
  \end{CD} \end{equation}
  where $(T_{Q,n}^{sc})^{\vee}[n]$ is a finite group: it is the group of $n$-torsion points in the maximal split torus of the adjoint quotient of $\overline{G}^{\vee}$.
  One can write the above exact sequence as:
  \begin{equation} \label{E:LG3}
  \begin{CD}
  \overline{G}^{\#} = \overline{G}^{\vee}  \rtimes (T_{Q,n}^{sc})^{\vee}[n] @>>> {^L}\overline{G}^{\#} @>>> W_F  \end{CD} \end{equation}
  where the action of $(T_{Q,n}^{sc})^{\vee}[n]$ on $\overline{G}^{\vee}$ is via the canonical adjoint action of $(\overline{G}^{\vee})_{ad}$ on $\overline{G}^{\vee}$.
  If one pulls back (\ref{E:LG2})  using the map 
  \[  {\rm id}  \times \chi_{\eta} : W_F \longrightarrow W_F \times \mu_2, \]
 one recovers the L-group extension ${^L}\overline{G}_{\eta}$. Thus, ${^L}\overline{G}^{\#}$ is an amalgam of all ${^L}\overline{G}_{\eta}$ and it always has a distinguished splitting.  Using this enlarged L-group, we can reduce the local Langlands correspondence for BD covering groups to a distinguished subclass of such groups, corresponding to $\eta =1$.  
 \vskip 5pt

   As an example, for  the groups $\overline{G}_{\eta}$ discussed in (i),    (\ref{E:LG2}) is:
  \[  \begin{CD} 
  \SL_2(\C)  @>>>  W_F \times \SL_2(\C)^{\pm}  @>\det>>    W_F \times \mu_2.
  \end{CD} 
  \]
whereas (\ref{E:LG3})  is:
\[  \begin{CD}
\SL_2(\C)^{\pm} @>>> W_F \times \SL_2(\C)^{\pm}  @>>> W_F  \end{CD} \]
\vskip 5pt

  \vskip 5pt

\item[(iv)]  Finally, we note that a distinguished splitting $s_0:  W_F \longrightarrow {^L}\overline{G}^{\#}$ gives rise to the notion  of {\em L-parameters relative to $s_0$}: for any splitting $s: W_F \longrightarrow {^L}\overline{G}^{\#}$, one sets
\[  \phi_s (w)  =  s(w)  / s_0(w)^{-1}  \in \overline{G}^{\#}  = \overline{G}^{\vee} \rtimes (T^{sc}_{Q,n})^{\vee}[n] \]
so that $\phi_s:  W_F \longrightarrow \overline{G}^{\#}$.  In particular, the Satake isomorphism furnishes the notion of {\em Satake parameters relative to $s_0$}. 
\vskip 5pt
 A  distinguished splitting $s_0$  thus allows one to define the notion of the $L$-function (with respect to $s_0$) associated to a representation $R$ of the (enlarged) dual group $\overline{G}^{\#}$.  Together with the Satake isomorphism, we thus have the notion of the ``automorphic L-function" associated to $s_0$ and $R$. 
In fact, if $R$ factors through the adjoint group of $\overline{G}^{\#}$, the resulting L-function is independent of the choice of $s_0$; in particular the adjoint L-function is well-defined. 
Another instance is the Langlands-Shahidi type L-functions, which will appear in the constant terms of Eisenstein series; this is the subject matter of the PhD thesis \cite{Ga} of the second author, which extends the results of Langlands' famous monograph \cite{La} to BD covering groups.
One would like to show, as in the linear case, that these automorphic L-functions are nice, i.e. have meromoprhic continuation and satisfy a functional equation of the usual type. This question is wide open at this moment.  
\vskip 5pt

\end{itemize}
\vskip 5pt  
\noindent We conclude the paper with a number of examples, illustrating the above constructions and results, as well as highlighting a few basic questions which we feel are crucial for carrying the theory forward. 
 \vskip 10pt
 
 \noindent{\bf Acknowledgments:} We thank Marty Weissman for many enlightening and helpful discussions, and for sharing with us his letter \cite{W4} to Deligne as well as his working manuscript \cite{W7}. We also like to acknowledge many conversations with Gordan Savin in which he shared his many insights on the topic of covering groups.
 \vskip 5pt
 
 This project was largely initiated when both authors were participating in the workshop ``Harmonic analysis and automorphic forms of covering groups" at the American Institute of Mathematics (AIM) in June 2013.  
 The authors would like to thank AIM for supporting the authors' participation.  The paper is completed when the first author is visiting the MSRI at Berkeley for the program ``New geometric methods in number theory and automorprhic forms"; the first author would like to thank MSRI for supporting his visit and providing excellent working conditions. 
 \vskip 5pt
 
 The research of both authors are partially supported by a Singapore government MOE Tier 2 grant   R-146-000-175-112.

 \vskip 10pt
 
\section{\bf Brylinski-Deligne Extensions}
Let $F$ be a field and $\G$ a connected reductive linear algebraic group over $F$. We shall assume that $\G$ is split over $F$ in this paper. Fix a maximal split torus $\T$ contained in a Borel subgroup $\mathbb{B}$ of $\G$.  Let  Let $Y = \Hom(\G_m, \T)$ be the cocharacter lattice and $X = \Hom(\T, \G_m)$ 
the character lattice. Then one has the set of roots  $\Phi \subset X$ and the set of coroots $\Phi^{\vee} \subset Y$
 of $(\G, \T)$ respectively.  The Borel subgroup $\B$ determines the set of simple roots $\Delta \subset \Phi$ and the set of simple coroots $\Delta^{\vee} \in \Phi^{\vee}$.  For each $\alpha \in \Phi$, one has the associated root subgroup $\mathbb{U}_{\alpha}  \subset \G$ which is normalised by $\T$.  
We shall  fix an \'epinglage or pinning for $(\G, \T)$, so that for each $\alpha \in \Phi$, one has an isomorphism $x_{\alpha}:  \G_a \longrightarrow \mathbb{U}_{\alpha}$.  

\vskip 5pt

Hence our initial data for this paper is a pinned connected split reductive group $(\G, \T, \B, x_{\alpha})$ over $F$. 
 \vskip 10pt

\subsection{\bf Multiplicative $K_2$-torsors.}
 The algebraic group  $\G$ defines a sheaf of groups on the big Zariski site on $Spec(F)$. Let $\mathbb{K}_2$ denote the sheaf of groups on ${\rm Spec}(F)$ associated to the 
$K_2$-group in Quillen's K-theory.  Then a multiplicative $K_2$-torsor is an extension
\[  \begin{CD}
1 @>>>  \mathbb{K}_2 @>>>    \overline{\G} @>>>  \G @>>>1 \end{CD}
\]
of sheaves of groups on $Spec(F)$.   We consider the category {\bf CExt}$(\G, \mathbb{K}_2)$ of such extensions where the morphisms between objects are given by morphisms of extensions. Given two such central extensions, one can form their Baer sum: this equips  {\bf CExt}$(\G, \mathbb{K}_2)$ with the structure of a commutative Picard category. 
\vskip 5pt

In \cite{BD}, Brylinski and Deligne made a deep study of   {\bf CExt}$(\G, \mathbb{K}_2)$  and obtained an elegant  classification of this category when  $\G$ is a connected reductive group.  We recall their results briefly in the case when $\G$ is split.
\vskip 5pt

\subsection{\bf Split torus.}  
Suppose $\T$ is a split torus, with cocharacter lattice $Y = \Hom(\G_m, \T)$ and character lattice $X = \Hom(\T, \G_m)$ 
the character lattice. Then we have:
\vskip 5pt

\begin{prop}  \label{P:tori}
Let $\T$ be a split torus over $F$. The category {\bf Cext}$(\T, \mathbb{K}_2)$ is equivalent as a commutative Picard category (by an explicit functor)   to the category whose objects are pairs $(Q, \mathcal{E})$, where
\vskip 5pt

\begin{itemize}
\item $Q$ is a $\Z$-valued quadratic form on $Y$, with associated symmetric bilinear form $B_Q(y_1, y_2)  = Q(y_1+y_2) - Q(y_1)  -Q(y_2)$;
\vskip 5pt

\item $\mathcal{E}$ is an extension of groups
\[  \begin{CD}
1 @>>>  F^{\times}  @>>>  \mathcal{E} @>>>  Y @>>> 1 \end{CD} \]
whose associated commutator map $[-,-]: Y \times Y \rightarrow F^{\times}$ is given by
\[  [y_1, y_2]  = (-1)^{B_Q(y_1,y_2)}. \]
\end{itemize}
The set of morphisms between $(Q,\mathcal{E})$ and $(Q', \mathcal{E}')$ is empty unless $Q  = Q'$, in which case it is given by the set of isomorphisms of extensions from $\mathcal{E}$ to $\mathcal{E}'$.  
The Picard structure is defined by
\[  (Q, \mathcal{E})  +  (Q' , \mathcal{E}')  = (Q + Q', \text{Baer sum of $\mathcal{E}$ and $\mathcal{E}'$}). \]
\end{prop}
\vskip 5pt

Observe that the isomorphism class of the extension $\mathcal{E}$ is completely determined by the commutator map and hence by the quadratic form $Q$.  The extension $\mathcal{E}$ is obtained from $\overline{\T}$ as follows.  Let $F(\!( \tau)\!)$ denote the field of Laurent series in the variable $\tau$ over $F$. Then one has
\[  \begin{CD}
K_2(F(\!(\tau)\!)) @>>> \overline{\T}(F(\!(\tau)\!))  @>>>  \T(F(\!(\tau)\!))= Y \otimes_{\Z} F(\!(\tau)\!)^{\times}.  \end{CD} \]
The map $y \mapsto y(\tau)$ defines a group homomorphism $Y \longrightarrow T(F(\!(\tau)\!))$. Pulling back by this morphism and pushing out by the residue map
\[  {\rm Res}: K_2(F(\!(\tau)\!))  \longrightarrow K_1(F) = F^{\times} \]
defined by
\[  {\rm Res}(f,g) = (-1)^{v(f) \cdot v(g)}  \cdot \left( \frac{f^{v(g)}}{g^{v(f)}}(0) \right), \]
one obtains the desired extension $\mathcal{E}$. 
\vskip 5pt

\subsection{\bf Simply-connected groups.}
Suppose now that $\G$ is a split simply-connected semisimple group over $F$ and recall that we have fixed the \'epinglage $(\T, \B, x_{\alpha})$.
 Let $W = N(\T)/\T$ be the corresponding Weyl group. 
Since $\G$ is simply-connected, the coroot lattice is equal to $Y$, so that the set of simple coroots  $\Delta^{\vee}$ is a basis for $Y$. 
  \vskip 5pt

Now we have:

\vskip 5pt

\begin{prop}  \label{P:sc}
The category {\bf CExt}$(\G, \mathbb{K}_2)$ is equivalent (as commutative Picard categories) to the category whose objects are $W$-invariant $\Z$-valued quadratic form $Q$ on $Y$, and whose only morphisms are the identity morphisms on each object. 
\end{prop}

\vskip 5pt

As a result of this proposition, whenever we are given a quadratic form $Q$ on $Y$, $Q$ gives rise to a multiplicative $K_2$-torsor $\overline{\G}_Q$ on $\G$, unique up to unique isomorphism, which may be 
pulled back to a multiplicative $K_2$-torsor $\overline{\T}_Q$ on $\T$ and hence gives rise to an extension $\mathcal{E}_Q$ of $Y$ by $F^{\times}$.   
The automorphism group of the extension $\mathcal{E}_Q$ is $\Hom(Y, F^{\times})$. Following [BD, \S 11], one can rigidify $\mathcal{E}_Q$ by giving it an extra structure, as we now explain.

\vskip 5pt

\vskip 5pt

 \subsection{\bf Rigidifying $\mathcal{E}_Q$.}
We continue to assume that $\G$ is simply-connected. 
 We have already fixed the \'epinglage 
 $\{ x_{\alpha}:  \alpha \in \Phi \}$ for $\G$, so that 
 \[  x_{\alpha}:  \G_a \longrightarrow \mathbb{U}_{\alpha} \subset \G. \]
 Indeed, one has an embedding
 \[  i_{\alpha}: \SL_2 \hookrightarrow  \G \]
 which restricts to $x_{\pm \alpha}$ on the upper and lower triangular subgroup of  unipotent matrices.
 By \cite{BD}, one has a canonical lifting
 \[  \tilde{x}_{\alpha} :  \G_a \longrightarrow \overline{\mathbb{U}}_{\alpha} \subset \overline{\G}. \]
 For $t \in \G_m$, we set
 \[   n_{\alpha}(t)   = x_{\alpha}(t)  \cdot x_{-\alpha}(-t^{-1})  \cdot x_{\alpha}(t) = i_{\alpha} \left(  \begin{array}{cc}
 0 & t \\
 -t^{-1} & 0  \end{array}  \right) \in N(\T_Q),  \]
 and
 \[  \tilde{n}_{\alpha}(t)  =  \tilde{x}_{\alpha}(t)  \cdot \tilde{x}_{-\alpha}(-t^{-1})  \cdot \tilde{x}_{\alpha}(t) \in \overline{\G}_Q    \]
 Then one has a map
 \[  s_{\alpha}:   T_{\alpha} := \alpha^{\vee}(\G_m)   \longrightarrow     \overline{\T}_{Q,\alpha} \]
 given by
 \[  \alpha^{\vee}(t)  \mapsto   \tilde{n}_{\alpha}(t) \cdot \tilde{n}_{\alpha}(-1). \]
 This is a section of $\overline{\G}_Q$ over $\T_{\alpha}$,  which is trivial at the identity element. 
  The section $s_{\alpha}$ is  useful in describing the natural conjugation action of $N(\T_Q)$ on $\overline{\T}_Q$.  By [BD, Prop. 11.3], one has the nice formula:
  \begin{equation} \label{E:nice-weyl}
   \tilde{n}_{\alpha}(1) \cdot \tilde{t}  \cdot \tilde{n}_{\alpha}(1)^{-1}  = \tilde{t} \cdot s_{\alpha}(\alpha^{\vee}(\alpha(t)). \end{equation}
  Moreover, 
  the collection of sections $\{ s_{\alpha}:  \alpha \in \Delta \}$ provides a collection of elements 
 $s_{\alpha}(\alpha^{\vee}(a)) \in \overline{\T}_Q$ with $a \in \G_m$,  and $\overline{\T}_Q$ is generated by $K_2$ and the collection of $s_{\alpha}(\alpha^{\vee}(a))$.

 \vskip 5pt

   Taking points in $F(\!(\tau)\!)$, we have the element 
 \[    s_{\alpha}( \alpha^{\vee}(\tau))  \in \overline{\T}_Q(F(\!(\tau)\!)), \]
 which gives rise (via the construction of $\mathcal{E}_Q$) to an element
 \[  s_Q(\alpha^{\vee})  \in   \mathcal{E}_Q.  \]
 Then we rigidify $\mathcal{E}_Q$ by equipping it with the set  $\{s_Q(\alpha^{\vee}):  \alpha^{\vee} \in \Delta^{\vee}\}$:  there is a unique automorphism of $\mathcal{E}_Q$ which fixes all these elements. 
 \vskip 5pt

   \vskip 5pt
{\em In the following, we shall fix a choice of the data  $(\overline{G}_Q, \overline{T}_Q,  \mathcal{E}_Q)$ for each $W$-invariant quadratic form $Q$ on $Y$ when $\G$ is split and simply-connected.} 
This is not a real choice, because any two choices are isomorphic by a unique isomorphism. 
The section $s_Q$ constructed above provides $\mathcal{E}_Q$ with a system of generators:
$\mathcal{E}_{Q}$ is generated by $s_Q(\alpha^{\vee}) \in \Delta^{\vee}$ and $a \in F^{\times}$ subject to the relations:
\vskip 5pt
\begin{itemize}
\item $a \in F^{\times}$ is central;
\item $[s_Q(\alpha^{\vee}),  s_Q(\beta^{\vee})]  =(-1)^{B_Q(\alpha^{\vee}, \beta^{\vee})}$ for $\alpha^{\vee}, \beta^{\vee} \in \Delta^{\vee}$.
\end{itemize}

\vskip 5pt
 
\vskip 5pt

\subsection{\bf General reductive groups.}
Now let $\G$ be a split connected reductive group over $F$, with fixed \'epinglage $(\T, \B, x_{\alpha})$. 
 Let 
\[   i_{sc}:   Y^{sc}  := \Z [ \Delta^{\vee}]  \subset Y \]
be the inclusion of the coroot lattice $Y^{sc}$ into $Y$, and let $X^{sc} \subset X \otimes_{\Z} \Q$ be the dual lattice of $Y^{sc}$. Then the quadruple $(X^{sc}, \Delta, Y^{sc}, \Delta^{\vee})$ is the root datum of the simply-connected cover $\G^{sc}$ of the derived group of $\G$, and one has a natural map 
\[  q: \G^{sc} \rightarrow \G.  \] 
Let $\T^{sc}$ be the preimage of $\T$ in $\G^{sc}$, so that one has a commutative diagram
\[  \begin{CD}
\T^{sc}  @>>> \G^{sc} \\
@VVV  @VVqV \\
\T  @>>> \G.  
\end{CD} \]

\vskip 5pt

The classification of the multiplicative $K_2$-torsors on $\G$ is an amalgam of the two Propositions above.
Given a multiplicative $K_2$-torsor $\overline{\G}$ of $\G$,  the above commutative diagram induces by pullbacks a commutative digram of multiplicative $K_2$-torsors:
\[  \begin{CD}
\overline{\T}^{sc}  @>>> \overline{\G}^{sc} \\
@VVV  @VVqV \\
\overline{\T}  @>>> \overline{\G}.  
\end{CD} \]
Then 

\begin{itemize}
\item  the multiplicative $K_2$-torsor $\overline{\T}$  gives  a pair $(Q, \mathcal{E})$ by Proposition \ref{P:tori};
\vskip 5pt

\item the multiplicative $K_2$-torsor $\overline{\G}^{sc}$  corresponds to a quadratic form $Q^{sc}$, and it was shown in \cite{BD} that $Q^{sc}$ is simply the restriction of $Q$ to $Y^{sc}$.

\vskip 5pt

\item we have fixed the data $(\overline{G}_{Q^{sc}}, \overline{T}_{Q^{sc}}, \mathcal{E}_{Q^{sc}})$  associated to $Q^{sc}$. Thus we have 
a canonical isomorphism 
\[ f: \overline{\G}_{Q^{sc}}  \longrightarrow \overline{\G}^{sc}, \]
restricting to  an isomorphism
 \[  f:  \overline{\T}_{Q^{sc}}  \longrightarrow  \overline{\T}^{sc} \]
 which then induces an isomorphism
 \[  f:   \mathcal{E}_{Q^{sc}}  \longrightarrow \mathcal{E}^{sc} =  q^*(\mathcal{E}). \]
 This isomorphism is characterised as the unique one which sends the elements $s_{Q^{sc}}(\alpha^{\vee}) \in \mathcal{E}_{Q^{sc}}$ (for $\alpha^{\vee} \in \Delta^{\vee}$) to the corresponding elements $s(\alpha^{\vee}) \in q^*(\mathcal{E})$. 
 In particular, we have a commutative diagram
 \[  \begin{CD}
 F^{\times} @>>>  \mathcal{E}_{Q^{sc}}  @>>>  Y^{sc} \\
 @|   @VVfV   @|  \\
 F^{\times} @>>> \mathcal{E}  @>>> Y^{sc}  
 \end{CD} \] 
 \end{itemize}
   
We have thus attached to a multiplicative $K_2$-torsor $\overline{\G}$ a  triple $(Q, \mathcal{E}, f)$.
Now we have:

 \vskip 5pt

\begin{thm}   \label{T:BDG}
The category {\bf CExt}$(\G, \mathbb{K}_2)$ is equivalent (via the above construction)  to the category ${\bf BD}_{\G}$ whose objects are triples $(Q, \mathcal{E}, f)$, where
\vskip 5pt

\begin{itemize}
\item $Q: Y \rightarrow \Z$ is a $W$-invariant quadratic form;

\item $\mathcal{E}$ is an extension of $Y$ by $F^{\times}$ with commutator map $[y_1, y_2] = (-1)^{B_Q(y_1, y_2)}$;

\item $f :  \mathcal{E}_{Q^{sc}}  \cong q^*(\mathcal{E})$ is an isomorphism of extensions of $Y^{sc}$ by $F^{\times}$. . 
\end{itemize}
The set of morphisms from $(Q, \mathcal{E}, f)$  to $(Q', \mathcal{E}', f')$ is empty unless $Q=Q'$, in which case it consists of isomorphisms of extensions $\phi:  \mathcal{E}  \longrightarrow \mathcal{E}'$ such that $f = f' \circ q^*(\phi)$.  In particular, the automorphism group of an object is $\Hom(Y/Y^{sc}, F^{\times})$.
\end{thm}
\vskip 5pt

\subsection{\bf Bisectors and Incarnation.}
While the above results give a nice classification of multiplicative $K_2$-torsors over split reductive groups $\G$ over $F$, it is sometimes useful and even necessary to work with explicit cocycles for computation. The paper \cite{BD} does provide a  category of nice algebraic  cocycles, as explicated in \cite{W3}, and one may replace the category of triples $(Q,\mathcal{E}, f)$ by a slightly simpler category with more direct connections to cocycles. 
\vskip 5pt

Extending the treatment in \cite{W3}, we consider (not necessarily symmetric) $\Z$-valued bilinear forms $D$ on $Y$ satisfying:
 \[  \text{$D(y,y)  = Q(y)$ for all $y \in Y$,} \]
 so that
\[  B_Q(y_1, y_2)  = D(y_1, y_2)  + D(y_2, y_1). \]
Such a $D$ is called a bisector of $Q$.   As shown in [W4], for any $Q$, there exists an associated bisector. 
\vskip 5pt

We consider a category
\[  {\bf Bis}_{\G}  = \bigcup_Q  {\bf Bis}_{\G,Q}  \]
where the full subcategory {\bf Bis}$_{\G,Q}$ consists of pairs $(D, \eta)$ where $D$ is a 
 bisector of $Q$ and 
\[  \eta:  Y^{sc}  \longrightarrow F^{\times} \]
is a group homomorphism. Given two pairs $(D_1, \eta_1)$ and $(D_2, \eta_2)$, the set of morphisms is the set of functions $\xi:  Y  \longrightarrow F^{\times}$ such that 
\vskip 5pt

\begin{itemize}
\item[(a)] $\xi(y_1 + y_2) \cdot \xi(y_1)^{-1} \cdot \xi(y_2)^{-1}  = (-1)^{D_1(y_1, y_2) - D_2(y_1, y_2)}$;
\item[(b)] $\xi(\alpha^{\vee}) = \eta_2(\alpha^{\vee}) / \eta_1(\alpha^{\vee})$ for all $\alpha^{\vee} \in \Delta^{\vee}$. 
 \end{itemize}
and where the composition of morphisms is given by multiplication: $\xi_1 \circ \xi_2(y)  = \xi_1(y) \cdot \xi_2(y)$.  
Note that in (b), $\xi$ may not be a group homomorphism when restricted to $Y^{sc}$, but we are only requiring the identity in (b) to hold as functions of sets when restricted to $\Delta^{\vee}$. 
Moreover, as shown in \cite{W3}, given two bisectors $D_1$ and $D_2$ of $Q$, one can always find $\xi$ such that (a) holds.  Thus, up to isomorphism, there is no loss of generality in fixing $D$ (for a fixed $Q$). 

\vskip 5pt
Then it was shown in  \cite{W3} that there is an incarnation functor
\[  {\bf Inc}_{\G}  :  {\bf Bis}_{\G} \longrightarrow {\bf BD}_{\G} \longrightarrow {\bf CExt}(\G, \mathbb{K}_2). \]
The second functor is a quasi-inverse to the functor in Theorem \ref{T:BDG}. 
On the level of objects, the first functor sends the pair $(D, \eta)$ in {\bf Bis}$_{\G,Q}$ to the triple $(Q, \mathcal{E}, f)$  defined as follows:
\vskip 5pt

\begin{itemize}
\item $\mathcal{E} = Y \times F^{\times}$ with group law
\[  (y_1, a_1) \cdot (y_2, a_2)  = (y_1 + y_2,  a_1a_2  (-1)^{D(y_1, y_2)}). \]

\item $f:  \mathcal{E}_{Q^{sc}}  \rightarrow \mathcal{E}$ is given by 
\[  f(s_{Q^{sc}}(\alpha^{\vee}))  =  ( \alpha^{\vee} ,  \eta(\alpha^{\vee})) \quad \text{for any $\alpha \in \Delta$}. \]
\end{itemize}
It was shown in [W4] that this  functor is an essentially surjective functor. In fact, our definition of morphisms in {\bf Bis}$_{\G}$ differs slightly from that of [W4]. With our version here, this functor is fully faithful as well so that {\bf Inc} is  in fact an equivalence of categories.
\vskip 5pt

Moreover, one can choose  {\bf Inc} so that if $\overline{\G}$ is the BD extension corresponding to $(D,\eta)$, then   $\overline{\T}$ can be described explicitly using $D$.  Namely, if
\[  D  = \sum_i  x^i_1  \otimes x^i_2  \in X \otimes X,  \]
then one may regard $\overline{\T} =  \T \times K_2$ with group law:
\[  (t_1, 1) \cdot  (t_2,1)  = (t_1t_2 ,  \prod_i  \{  x^i_1(t_1),  x^i_2(t_2) \}). \]
The associated extension $\mathcal{E}$ is then described as above in terms of the bisector $D$.
Further, we have the following explicit description for  the section $q \circ s_{\alpha}$:

\begin{prop}
In terms of the above realisation of $\overline{T}$ by $D$, the section 
\[  q \circ s_{\alpha}: \alpha^{\vee}(\G_m) \longrightarrow  \overline{\T} \quad \text{for $\alpha \in \Delta$} \]
is given by
\[  q\circ s (\alpha^{\vee}(a))  = (  \alpha^{\vee}(a),  \{ \eta(\alpha^{\vee}), a \}). \]
\end{prop}
\vskip 5pt

 \begin{proof}
 Suppose that $\overline{\G}$ is incarnated by $(D,\eta)$ and  corresponds to the triple $(Q, \mathcal{E}, f)$ where $\mathcal{E}$ and $f$ are described by $(D,\eta)$ as above.
 For fixed $\alpha\in \Delta$, one may write  
$$  q \circ s_{\alpha}(a) =  (\alpha^{\vee}(a),  \aleph_\alpha(a)),$$
where $\aleph_\alpha \in \Hom_\text{Zar}(\G_m, \mathbb{K}_2)$ is a homomorphism of sheaves of abelian groups for the big Zariski site. 
 By \cite[\S 3.7-3.8]{BD}, or in more details \cite[Thm. 1.1]{Bl},  we have 
\[
F^{\times} = \mathbb{K}_1(F) \cong  \Hom_\text{Zar}(\G_m, \mathbb{K}_2) \]
 where the isomorphism is given by
\[  b \mapsto  (a \mapsto \{ b,a \}). \]
 Hence, there exists some $\lambda_{\alpha} \in F^{\times}$ such that
$$\aleph_\alpha(a)=\{\lambda_\alpha, a\} \quad \text{for $a \in \G_m$}.$$
From this, it follows from the definition of $f$ that 
\[  f(s_{Q^{sc}}(\alpha^{\vee})) =  (\alpha^{\vee},  \lambda_{\alpha} ) \in \mathcal{E}. \]
By hypothesis, however, we have: 
\[  f(s_{Q^{sc}}(\alpha^{\vee})) = (\alpha^{\vee},  \eta(\alpha^{\vee})).\]
Hence we deduce that
\[  \lambda_{\alpha}  = \eta(\alpha^{\vee}) \]
and thus
\[  
 q \circ s_{\alpha}(a) =  (\alpha^{\vee}(a),  \aleph_\alpha(a)) =  (\alpha^{\vee}(a),  \{ \eta(\alpha^{\vee}), a\}), \]
 as desired.
\end{proof}

 Thus, the category {\bf Bis}$_{\G}$ provides a particularly nice and explicit family of  cocycles for BD extensions, and the essential surjectivity of {\bf Inc} says that every BD extension possesses such a cocycle, at least on the maximal torus $\T$.  This will be useful for computation.
\vskip 10pt

 \subsection{\bf Fair bisectors.}
In \cite{W3}, Weissman singled out a property of bisectors which he called {\em fair}. By definition, a bisector $D$ is fair if it satisfies the following:
\vskip 5pt
\begin{itemize}
\item for any $\alpha \in \Delta$ such that $Q(\alpha^\vee) \equiv 0 \mod 2$, $D(\alpha^\vee, y) \equiv D(y, \alpha^\vee)  \equiv 0 \mod 2$ for all $y \in Y$.
\end{itemize}
Weissman showed that  for any $Q$, ${\bf Bis}_{\G, Q}$ contains a fair bisector. We shall henceforth fix a fair bisector for a given $Q$.  The value of fairness will be apparent later on.  
 
\vskip 5pt

At the moment, we simply note that the
  objects  $(D, 1) \in {\bf Bis}_{\G,Q}$ with $D$ fair and  $\eta$ the trivial homomorphism) are quite special (as we shall see). Thus, we have a distinguished class of multiplicative $K_2$-torsors with invariant $Q$.    \vskip 5pt
  
  When $Q = 0$, for example, the  bisector $D = 0$ is fair, and $(D, 1)$ gives the 
split extension $\G \times K_2$.  In some sense,  the $K_2$-torsor with invariants $(D, 1)$ should be regarded as ``closest to being a split extension  among those with invariants $(D, \eta)$". 
As we shall illustrate in the rest of this section, the general BD extensions can often be described in terms of these distinguished BD extensions. 

\vskip 5pt
\subsection{\bf The case $Q = 0$.}
Let us consider the example when $Q = 0$. Then we may take the bisector $D = 0$ and regard the objects of ${\bf Bis}_Q$ as the set of homomorphisms $\eta:  Y^{sc} \rightarrow F^{\times}$. Let $\overline{\G}_{\eta}$ be the corresponding multiplicative $K_2$-torsor on $\G$. Then $\overline{\G}_{\eta_1}$ and $\overline{\G}_{\eta_2}$
are isomorphic precisely when $\eta_1/\eta_2$ can be extended to a homomorphism of $Y$ to $F^{\times}$. 
\vskip 5pt

How can we characterize the distinguished BD-extension in ${\bf Bis}_Q$ using the BD data $(Q,\mathcal{E}, f)$ when $D=0$? Since $D=0$, $\mathcal{E}$ is an abelian group and hence is a split extension:  $\mathcal{E}  = Y \times F^{\times}$.  
Each $\eta \in \Hom(Y^{sc}, F^{\times})$ gives a map 
\[  f_{\eta} :  \mathcal{E}_{Q^{sc}} \longrightarrow \mathcal{E} \]
defined by
\[  f_{\eta}(s_{Q^{sc}}(\alpha^{\vee}))  = (\alpha^{\vee},  \eta(\alpha^{\vee})) \quad \text{for $\alpha \in \Delta$.} \]
Applying the functor $\Hom_{\Z}(-, F^{\times})$ to the short exact sequence $Y^{sc} \rightarrow Y \rightarrow Y/Y^{sc}$, we obtain as part of the long exact sequence:
 \[ \begin{CD}
  \Hom_{\Z}(Y, F^{\times})  @>>> \Hom_{\Z}(Y^{sc}, F^{\times}) @>\delta>>  {\rm Ext}^1_{\Z}(Y/Y^{sc}, F^{\times}) @>>> 0. \end{CD} \]
  For any $\eta \in \Hom(Y^{sc}, F^{\times})$, its image under $\delta$ is the extension 
  \[  \begin{CD}
  F^{\times} @>>>  (Y \times F^{\times})/ f_{\eta}(Y^{sc}) @>>>  Y/Y^{sc}  \end{CD} \]
  The long exact sequence thus implies that this extension is split precisely when $\eta$ is equivalent to $1$.
  This gives a way of characterising the distinguished isomorphism class in ${\bf Bis}_Q$ when $D=0$. 
 \vskip 5pt

How can we construct the other BD extensions in ${\bf Bis}_Q$ when $D = 0$? 
We assume further that $\G$ is semisimple and let $q:  \G^{sc}  \longrightarrow \G$ be the natural isogeny with kernel  $Z = {\rm Tor}_{\Z}(Y/Y^{sc}, \G_m) \hookrightarrow \T^{sc} = Y^{sc} \otimes_{\Z}  \G_m$. 
Then $q^*(\overline{\G}_{\eta})$ is isomorphic to the split extension $\G^{sc}  \times K_2$ by a unique isomorphism. We may thus construct $\overline{\G}_{\eta}$ by starting with $\G^{sc}  \times K_2$ and then considering a quotient of this by a suitable embedding $Z \hookrightarrow \G^{sc} \times K_2$.  
 For this, we note that $\eta$ induces a map
 \[  i_{\eta} :  Z = {\rm Tor}_{\Z}(Y/Y^{sc}, \G_m) \hookrightarrow Y^{sc} \otimes_{\Z} \G_m  \longrightarrow F^{\times}\otimes_{\Z} \G_m \longrightarrow K_2. \]
 Then we have
 \[  \overline{\G}_{\eta}  =( \overline{\G} \times K_2 ) /  \{  (z, i_{\eta}(z) : z \in Z \}. \]
 \vskip 5pt

 \subsection{\bf $z$-extensions.}
 We consider another example which will play a crucial role later on, namely when $Y/Y^{sc}$ is a free abelian group. In this case, for any $(D, \eta) \in {\bf Bis}_Q$, $\eta$ can be extended to a homomorphism of $Y$, and so any two ($D, \eta_1)$ and $(D, \eta_2)$ are isomorphic. This means that there is a unique isomorphism class of objects in ${\bf Bis}_Q$, just like the case when $\G$ is simply-connected (where $Y = Y^{sc}$). However, the automorphism group $\Hom(Y/Y^{sc}, F^{\times})$ of an object is not trivial (unless $Y = Y^{sc}$).
 \vskip 5pt
 
 As we shall see later, some questions about a BD extension can be reduced to the case when $Y/Y^{sc}$ is free. This is achieved via the consideration of $z-$extensions.  More precisely, it is not hard to see that given any $\G$, one can find a central extension of connected reductive groups: 
 \[  \begin{CD}
 Z @>>>    \mathbb{H} @>\pi>> \G  \end{CD}  \]
 where 
 \vskip 5pt
 
 \begin{itemize}
 \item  $\mathbb{H}$ is such that $Y_{\mathbb{H}}/ Y_{\mathbb{H}}^{sc}$ is free;
 \item $\pi_*:  Y_{\mathbb{H}}^{sc} \longrightarrow Y_{\mathbb{H}}$ is an isomorphism;
 \item $Z$ is a split torus which is central in $\mathbb{H}$.
 \end{itemize}
 Such an extension is called a $z$-extension. 
 \vskip 5pt
 
 Given such a $z$-extension, and a BD extension $\overline{G}_{\eta}$ with invariant $(D, \eta)$, we obtain a BD extension 
 $\mathbb{H}_{\eta} := \pi^*(\overline{G}_{\eta})$ on $\mathbb{H}$ with  BD invariant $(D \circ \pi, \eta \circ \pi)$, so that one has
 \[  \begin{CD}
 Z_{\eta}  =Z  @>>>  \overline{\mathbb{H}}_{\eta}  @>>> \overline{\G}_{\eta} \end{CD} \]
 By our discussion above, we may choose an isomorphism $\xi:  \overline{\mathbb{H}}_1  \cong \overline{\mathbb{H}}_{\eta}$. Thus, via $\xi$, we have
 \[  \begin{CD}
  \xi^{-1}(Z_{\eta} ) @>>> \overline{\mathbb{H}}_1 @>>>   \overline{\G}_{\eta}  \end{CD} \]
  This shows that for any given bisector $D$,  all the BD extensions $\overline{\G}_{\eta}$ can be described as the quotient of a fixed BD extension $\overline{\mathbb{H}}_1$ (with $\eta =1$) by a suitable splitting of the split torus $Z$. 
  \vskip 5pt
  
 \subsection{\bf Running example.}
  Let us illustrate the above discussion using a simple example, where 
  \[  \G = \PGL_2,    \quad  Y = \Z \supset Y^{sc}  =2\Z  \quad \text{and} \quad D= 0. \]
   Then we take the  $z$-extension to be 
  \[  \begin{CD} \G_m  @>>> \GL_2 @>\pi>> \PGL_2 \end{CD} \]
 The distinguished BD-extension with $\eta=1$ is the split extension  $\overline{\G}_1 = \PGL_2 \times K_2$ and its pullback to $\GL_2$ is the split extension
$\overline{\mathbb{H}} = \pi^*(\overline{\G}_1) = \GL_2 \times K_2$. 
  For any $\eta \in \Hom(Y^{sc}, F^{\times})  \cong  F^{\times}$, the BD extension $\overline{\G}_{\eta}$ can then be described as
  \[  \overline{\G}_{\eta} =  (\GL_2 \times K_2) /  \{ (z, i_{\eta}(z)): z \in \G_m \}, \]
  where
  \[  i_{\eta}(z)  = \{ \eta, z \} \in K_2. \]
 This is a rather trivial family of BD extensions since their pullback to $\G^{sc} = \SL_2$ is split.
 Nonetheless, we shall use them as our running examples, as they already exhibit 
 various properties we want to highlight in this paper.  
 \vskip 10pt

\section{\bf Topological Covering Groups}  \label{S:tcg}

In this section, we will pass from the algebro-geometric world of multiplicative $K_2$-torsors to the world of topological central extensions. 
Let $F$ be a local field. If $F$ is nonarchimedean, let $\mathcal{O}$ denote its ring of integers with residue field $\kappa$. 
\vskip 5pt

\subsection{\bf BD covering groups.}
Start with a multiplicative $K_2$-torsor 
 $K_2 \rightarrow \overline{\G}  \rightarrow \G$, with associated BD data $(Q,\mathcal{E}, f)$ or bisector data $(D,\eta)$. 
By taking $F$-points, we obtain (since $H^1(F, K_2)  = 0$) a short exact sequence of abstract groups
\[  \begin{CD}
K_2(F)  @>>>  \overline{\G}(F)  @>>> G = \G(F)  \end{CD} \]
Now let $\mu(F)$ denote the set of roots of unity contained in the local field $F\ne \C$; when $F = \C$, we let $\mu(F)$ be the trivial group.  
Then the Hilbert symbol gives a map
\[  (-,-)_F:  K_2(F)  \longrightarrow \mu(F). \]
For any $n$ dividing $\# \mu(F)$, one has the $n$-th Hilbert symbol
\[  (-,-)_n   = (-,-)_F^{\#\mu(F) / n} : K_2(F)  \longrightarrow \mu_n(F). \]
Pushing the above exact sequence out by the $n$-th Hilbert symbol, one obtains a short exact sequence of locally compact topological groups
\[  \begin{CD} 
\mu_n(F)  @>>> \overline{G}  @>>> G  \end{CD} \]
We shall call this the BD covering group associated to the BD data $(Q, \mathcal{E}, \phi, n)$, or to the bisector data $(D, \eta, n)$. 

 \vskip 5pt
 
 Since we are considering degree $n$ covers, it will be useful to refine certain notions taking into account  the extra data $n$:
 \vskip 5pt
 
 \begin{itemize}
 \item for a bisector data $(D,\eta)$, we write $\eta_n$ for the composite
\[  \eta_n :  Y^{sc} \longrightarrow F^{\times} \longrightarrow F^{\times}/F^{\times n}. \]
\vskip 5pt

\item with $D$ fixed, we say that $\eta_n$ and $\eta'_n$ are equivalent if $\eta_n/\eta'_n$ extends to a homomorphism $Y \longrightarrow F^{\times}/ F^{\times n}$. 
 \end{itemize}
 \vskip 5pt

 \subsection{\bf Canonical unipotent section.}
Because a BD extension is uniquely split over any unipotent subgroup, one has unique splittings:
\[  \tilde{x}_{\alpha}:   F \longrightarrow \overline{U}_{\alpha}   \quad \text{for each $\alpha \in \Phi$.} \]
Indeed, as shown in \cite{L2}, there is a unique section
\[  i:  \{  \text{all unipotent elements of $G$} \} \longrightarrow \overline{G}  \]
satisfying:
\vskip 5pt

\begin{itemize}
\item for each unipotent subgroup  $\mathbb{U} \subset \G$, the restriction of $i$ to $U = \mathbb{U}(F)$ is a group homomorphism;
\item the map $i$ is $G$-equivariant.
\end{itemize}
For example, for each $\alpha \in \Phi$, we have seen that there is a homomorphism 
\[  \tilde{x}_{\alpha}  :  \G_a \longrightarrow \overline{\G}  \]
which induces a homomorphism
\[  \tilde{x}_{\alpha}:  F \longrightarrow \overline{G} \]
lifting the inclusion $x_{\alpha}  : F \hookrightarrow G$. 
Then one has  $\tilde{x}_{\alpha}  =  i \circ  x_{\alpha}$.  
\vskip 10pt

\subsection{\bf Covering torus $\overline{T}$.}  \label{SS:coveringT}
We may consider the pullback of $\overline{G}$ to the maximal split torus $T$:
\[  \begin{CD}  
\mu_n(F) @>>> \overline{T}  @>>> T   \end{CD} \]
As we observe in the last section, the bisector $D$ furnishes a cocycle for the multiplicative $K_2$-torsor on $T$. 
Thus the covering torus $\overline{T}$ also inherits a nice cocycle,  giving us a rather concrete description of $\overline{T}$.

\vskip 5pt

More precisely,  $\overline{T}  = T \times_D \mu_n(F)$  is generated by elements $\zeta \in \mu_n(F)$ and $y(a)$, for $y \in Y$ and $a \in F^{\times}$, subject to the relations:

\begin{itemize}
\item the elements $\zeta$ are central;
\item $[ y_1(a), y_2(b)]  = (a,b)_n^{B_Q(y_1, y_2)}$ for all $y_1, y_2 \in Y$ and $a,b \in F^{\times}$;
 \item $y_1(a) \cdot y_2(a)  = (y_1 + y_2)(a) \cdot (a,a)_n^{D(y_1, y_2)}$; 
 \item  $y(a) \cdot y(b) =  y(ab) \cdot (a,b)_n^{Q(y)}$. 
 \end{itemize}
As the second relation shows, $\overline{T}$ is not necessarily an abelian group. 
Moreover, the sections $q \circ s_{\alpha}$ for $\alpha \in \Delta$ takes the form
\[  q \circ s_{\alpha} (\alpha^{\vee}(a)) = \alpha^{\vee}(a) \cdot (\eta(\alpha^{\vee}), a)_n. \]

\vskip 5pt

\subsection{\bf The torus $T_{Q,n}$.}  \label{SS:coveringTQn}
The center $Z(\overline{T})$ of $\overline{T}$ is generated by $\mu_n(F)$ and those elements $y(a)$ such that
\[  B_Q(y, z)  \in n\Z  \quad \text{for all $z \in Y$.} \]
Thus,  we define:
\[  Y_{Q,n} :=   Y \cap n Y^*, \]
where $Y^* \subset Y \otimes_{\Z}  \Q$ is the dual lattice of $Y$ relative to $B_Q$. 
Then the center of $\overline{T}$ is generated by $\mu_n(F)$ and the elements $y(a)$ for all $y \in Y_{Q,n}$ and $a \in F^{\times}$.  It is clear that $nY \subset Y_{Q,n}$.
\vskip 5pt

Let $T_{Q,n}$ be the split torus with cocharacter group $Y_{Q,n}$. The inclusion $Y_{Q,n} \hookrightarrow Y$ gives an isogeny of tori
\[  f:  T_{Q,n}  \longrightarrow  T. \]
We may pullback the covering $\overline{T}$ using $f$, thus obtaining a covering torus $\overline{T}_{Q,n}$:
\[  \begin{CD}
\mu_n(F) @>>> \overline{T}_{Q,n}  @>>> T_{Q,n}  \end{CD} \]
Then $\overline{T}_{Q,n}$ is generated by $\zeta \in \mu_n(F)$ and elements $y(a)$ with $y \in Y_{Q,n}$ with the same relations as those for $\overline{T}$. However, the second relation now becomes
\[  [y_1(a), y_2(b) ]  =1,  \quad \text{for all $y_1, y_2 \in Y_{Q,n}$ and $a,b\in F^{\times}$}. \]
Thus $\overline{T}_{Q,n}$ is an abelian group. Moreover, it follows from the definition of pullbacks that there is a canonical homomorphism ${\rm Ker}(f) \rightarrow  \overline{T}_{Q,n}$, so that one has a short exact sequence of topological groups
\[  \begin{CD}
{\rm Ker}(f)  @>>> \overline{T}_{Q,n}  @>>> Z(\overline{T})  \end{CD} \]
In particular, to give a character of $Z(\overline{T})$ is to give a character of $\overline{T}_{Q,n}$ trivial on the subgroup ${\rm Ker}(f)$.   

\vskip 10pt

\subsection{\bf The kernel of $f$}
We need to  have a better handle of ${\rm Ker}(f)$.  
  The inclusions $nY \longrightarrow Y_{Q,n}  \longrightarrow Y$ give rise to  isogenies 
\[ \begin{CD}
T @>g>>  T_{Q,n} @>f>> T \end{CD} \]
so that $f \circ g$ is the $n$-power map on $T$.
 We have:
\vskip 5pt

\begin{lemma}  \label{L:TQn}
The kernel of $f$ is contained in the image of $g$. Indeed, ${\rm Ker}(f)  = g(T[n])$.  
\end{lemma}
\vskip 5pt

\begin{proof}
By the elementary divisor theorem, we may pick a basis $\{  e_i\}$of $Y$ such that a basis of $Y_{Q,n}$ is given by $\{  k_i e_i \}$ for some positive integers $k_i$.  Such bases allow us to identify the maps
\[ \begin{CD}
T = (F^{\times})^r  @>g>>  T_{Q,n} = (F^{\times})^r @>f>> T  = (F^{\times})^r \end{CD} \]
explicitly as
\[  g(t_i)  =  (t_i^{n/k_i})  \quad \text{and}  \quad f(t_i)  =(t_i^{k_i}). \]
Thus, 
\[  {\rm Ker}(f)  = \{  ( \zeta_i):  \zeta_i^{k_i} =1 \}  = g(T[n])  = \{  (\zeta_i^{n/k_i}):  \zeta_i^n =1 \}. \] 
\end{proof}
\vskip 5pt
Using the generators and relations for $\overline{T}_{Q,n}$, it is easy to see that 
the map
\[    y(a)    \mapsto  (ny)(a)    \in \overline{T}_{Q,n}  \]
gives a group homomorphism $\tilde{g}:  T   \longrightarrow \overline{T}_{Q,n}$. 
 The lemma implies that a character of $\overline{T}_{Q,n}$ trivial on the image of $\tilde{g}$  necessarily factors through to a character of $\overline{T}$. 
 \vskip 5pt

\section{\bf Tame Case} \label{S:tame}

Suppose now that $F$ is a $p$-adic field with residue field $\kappa$.
Assume that  $p$ does not divide $n$: we shall call this the tame case. 
The main question we want to consider in this section is whether a tame BD cover $\overline{G}$ necessarily splits over a hyperspecial maximal compact subgroup $K$ of $G$.  
 \vskip 5pt

Since we have fixed a Chevalley system of \'epinglage for $\G$,  we 
have its associated maximal compact subgroup 
 $K $ generated by $x_{\alpha}(\mathcal{O}_F)$ for all $\alpha \in \Phi$ and the maximal compact subgroup  $T(\mathcal{O})  = Y \otimes_{\Z} \mathcal{O}^{\times}$ of $T$. In particular $K = \underline{\G}(\mathcal{O})$ for a smooth reductive group $\underline{\G}$ over $\mathcal{O}$.  To ease notation, we shall simply write $\G$ for $\overline{\G}$ in what follows.
 Let $\G_{\kappa}$ denote the special fiber $\G \times_{\mathcal{O}} \kappa$ of $\G$. 
One has a natural reduction map 
\[   \G(\mathcal{O}) \longrightarrow G_{\kappa} := \G_{\kappa}(\kappa), \]
whose kernel is a pro-$p$ group. Restricting the BD cover to $K$, one has a topological central extension
\[  \begin{CD}
\mu_n  @>>>  \overline{\G}(\mathcal{O})  @>>>  \G(\mathcal{O}).  \end{CD} \]
Here, observe that we have abused notation and write $ \overline{\G}(\mathcal{O})$ for the inverse image of $\G(\mathcal{O})$ in $\overline{G}$.  
We would like to determine if this extension splits. 
\vskip 5pt

\vskip 5pt

\subsection{\bf The tame extension.}
All extensions in the tame case arise in the following way from the multiplicative $K_2$-torsor $\overline{\G}$. 
The prime-to-$p$ part of $\mu(F)$ is naturally isomorphic to $\kappa^{\times}$, and 
 there is a tame symbol
$K_2(F)  \longrightarrow \kappa^{\times}$ defined by
\[    \{a,b \}  \mapsto \text{ the image of $(-1)^{{\rm ord}(a) \cdot {\rm ord}(b)} \cdot \frac{a^{{\rm ord}(b)}}{b^{{\rm ord}(a)}}$ in $\kappa^{\times}$.} \]
Pushing out by this tame symbol gives the tame extension
\[  \begin{CD}
 \kappa^{\times}  @>>>  \overline{G}^{{\rm tame}} @>>>  G  \end{CD} \]
Hence any degree $n$ BD extension $\overline{G}$  with $(n,p)  =1$ is obtained as a pushout of $\overline{G}^{{\rm tame}}$. 

\vskip 5pt

\subsection{\bf Residual extension.}
We shall consider first the case of the tame extension $\overline{G}^{{\rm tame}}$ so that 
$n  = \#  \kappa^{\times}$.
It was shown in \cite[\S 12]{BD} and \cite{W2} that there is an extension of reductive algebraic group over $\kappa$:
\[  \begin{CD}
\G_m @>>>  \tilde{\G}_{\kappa}  @>>>  \G_{\kappa}  \end{CD} \]
with the following property: 
\vskip 5pt

\begin{itemize}
\item for any unramified extension $F'$ of $F$ with ring of integers $\mathcal{O}'$ and residue field $\kappa'$,  the tame extension 
\[  \begin{CD}
{\kappa'}^{\times}  @>>>  \overline{\G}(\mathcal{O'}) @>>>    \underline{\G}(\mathcal{O}')  \end{CD} \]
is the pullback of the extension
\[ \begin{CD}
\G_m(\kappa') = {\kappa'}^{\times} @>>>    \tilde{G}_{\kappa'}  = \tilde{\G}_{\kappa}(\kappa') @>>>  G_{\kappa'}  = \G_{\kappa}(\kappa') 
\end{CD} \]
with respect to the reduction map $\G(\mathcal{O'})  \longrightarrow G_{\kappa'}$. 
 \end{itemize}
We call this  extension of algebraic groups over $\kappa$ the residual extension.  
\vskip 5pt

\subsection{\bf Classification.}  \label{SS:class}
 
In \cite{W2}, algebraic extensions of $\G_{\kappa}$ by $\G_m$  were classified in terms of enhanced root theoretic data similar in spirit to (but simpler than)  the BD data. We give a sketch in the case when $\G_{\kappa}$ is split. Then such extensions are classified by the category of pairs $(\mathcal{E}_{\kappa}, f_{\kappa})$ with
\vskip 5pt

\begin{itemize}
\item  $\Z \rightarrow  \mathcal{E}_{\kappa} \rightarrow Y$ is an extension of free $\Z$-modules;
\item $f_{\kappa}:  Y^{sc}  \longrightarrow   \mathcal{E}_{\kappa}$ is a splitting of  $\mathcal{E}_{\kappa}$ over $Y^{sc}$. 
\end{itemize}
Moreover, the extension $\tilde{\G}_{\kappa}$ in question is split if and only if the map $f_{\kappa}$ can be extended to a splitting $Y \longrightarrow   \mathcal{E}_{\kappa}$, or equivalently,  if and only if the extension
 \[  \begin{CD}
 \Z @>>>  \mathcal{E}_{\kappa} /  f_{\kappa}(Y^{sc}) @>>>  Y/Y^{sc}  \end{CD} \]
 is split. This holds for example if $Y/Y^{sc}$ is free. 
 In particular, if $\G = \G^{sc}$ is simply connected, then $\tilde{\G}_{\kappa}$ is split and the splitting is unique.

\vskip 5pt
Given an extension $\tilde{\G}_{\kappa}$, one obtains the above two data as follows.  If $\T_{\kappa}$ is a maximal split torus of $\G_{\kappa}$, with preimage $\tilde{\T}_{\kappa}$, then the cocharcter lattice of $\tilde{\T}_{\kappa}$ gives the extension $\mathcal{E}_{\kappa}$.  The pullback of $\tilde{\G}_{\kappa}$ to $\G^{sc}_{\kappa}$  is canonically split. On restricting  this canonical splitting to the maximal split torus ${\T}_{\kappa}^{sc}$
(which is the pullback of $\T_{\kappa}$), one obtains the splitting $f_{\kappa}$ on the level of cocharacter lattices. 
\vskip 5pt

\subsection{\bf Splitting of $\overline{K}$.}  \label{SS:split}
Since the kernel of the reduction map is a pro-$p$ group, the set of splittings of the topological extension $\overline{K} = \overline{\G}(\mathcal{O})$ 
is in bijection with those  of the abstract extension $\tilde{G}_{\kappa}$. Further, a splitting of the residual extension $\tilde{\G}_{\kappa}$ gives rise to a splitting of $\tilde{G}_{\kappa}$ and thus of $\overline{\G}(\mathcal{O})$.  We will investigate the existence of splittings for $\tilde{\G}_{\kappa}$:  they give rise to splittings of $\overline{\G}(\mathcal{O})$ of ``algebraic origin". 
\vskip 5pt

For example, when $\G = \G^{sc}$ is simply-connected, the unique splitting of the residual extension
$\tilde{\G}_{\kappa}$  gives rise to a unique compatible system of splittings of $\G(\mathcal{O}')$ for all unramified extensions $\mathcal{O}'$ of $\mathcal{O}$.  Indeed, one has a natural bijection
\[  \{  \text{splittings of residual extension $\tilde{\G}_{\kappa}$} \} \longleftrightarrow  \{  \text{compatible system of splittings of $\G(\mathcal{O}')$} \}. \]
  
\vskip 5pt

 \subsection{\bf Determining $\tilde{\G}_{\kappa}$.}
One can now figure out the residual extension $\tilde{\G}_{\kappa}$ obtained from a BD extension $\overline{\G}$ with associated BD data $(Q, \mathcal{E},  f)$. 
\vskip 5pt

\begin{prop}
If $\overline{\G}$ has BD data $(Q, \mathcal{E}, f)$, then the associated data $(\mathcal{E}_{\kappa}, f_{\kappa})$ for the residual extension $\tilde{\G}_{\kappa}$ is obtained as follows:
\vskip 5pt
\begin{itemize}
\item $\mathcal{E}_{\kappa}$ is the  pushout of $\mathcal{E}$ by the valuation map ${\rm ord}:  F^{\times} \longrightarrow \Z$;
\item $f_{\kappa}$ is deduced from the associated map ${\rm ord}_* \circ f : \mathcal{E}_{Q^{sc}} \rightarrow \mathcal{E} \rightarrow \mathcal{E}_{\kappa}$:
\[  f_{\kappa}(\alpha^{\vee})  =   {\rm ord}_*( f( s_{Q^{sc}}(\alpha^{\vee}))) \quad \text{for $\alpha \in \Delta$.} \]
\end{itemize}
\end{prop}
\vskip 5pt
 
\begin{proof}
This question is systematically and more elegantly addressed in \cite{W6}, and we shall give a more  ad hoc argument here.  
We assume that $\overline{\G}$ is incarnated by $(D,\eta)$ for concreteness. As we discussed above, for any unramified extension $\mathcal{O}'$ of $\mathcal{O}$,
 there is a commutative diagram of extensions:
\[  \begin{CD}
{\kappa'}^{\times} @>>> \overline{\G}(\mathcal{O'})  @>>>  \G(\mathcal{O'})  \\
@VVV  @VVV  @VVV  \\
\G_m(\kappa')  @>>>  \tilde{\G}_{\kappa}(\kappa')  @>>>  \G_{\kappa}(\kappa'), \end{CD}  \]
and our goal is to determine the invariants $(\mathcal{E}_{\kappa}, f_{\kappa})$ for the extension $\tilde{\G}_{\kappa}$.
\vskip 5pt

Now  we note:
\vskip 5pt

\begin{itemize}
\item[(a)]  From the construction in \cite[\S 12.11]{BD},  one has a commutative diagram of extensions and splittings ($s_{\kappa}$ of $\tilde{\G}_{\kappa}^{sc}$) over $\kappa$:
\[  \begin{CD}
\T^{sc}_{\kappa}  @>s_{\kappa}>> \tilde{\T}^{sc}_{\kappa}  @>>> \tilde{\T}_{\kappa}   \\
@VVV @VVV   @VVV  \\
\G_{\kappa}  @>s_{\kappa}>>\tilde{\G}^{sc}_{\kappa}   @>>>  \tilde{\G}_{\kappa}  \end{CD}  \]
pulling back to a compatible system
\[  
  \begin{CD}
\T^{sc}(\mathcal{O}')  @>s_{\mathcal{O'}}>> \overline{\T}^{sc}(\mathcal{O}') @>>>  \overline{\T}(\mathcal{O}') \\
@VVV  @VVV @VVV \\
\G^{sc}(\mathcal{O}') @>s_{\mathcal{O'}}>>  \overline{\G}^{sc}(\mathcal{O}')   @>>> \overline{\G}(\mathcal{O}').   \end{CD}  \]
 \vskip 5pt

\item[(b)]  Using the description of $\overline{\T}$ in terms of $(D,\eta)$, one sees immediately that
$\overline{\T}(\mathcal{O}')  = \T(\mathcal{O}') \times {\kappa'}^{\times}$ (a direct product of groups), since the tame symbol is trivial on $\mathcal{O'}^{\times} \times \mathcal{O'}^{\times}$. Thus $\tilde{\T}{\kappa}  =  \T_{\kappa}  \times \G_m$.  
From this, one deduces that 
\[  \mathcal{E}_{\kappa}  =  Y \times \Z \quad \text{(as groups). }  \]
\vskip 5pt

\item[(c)]  The invariant $f_{\kappa}$ is deduced from the composite map 
\[ \begin{CD}  \tilde{f}_{\kappa}:  \T^{sc}_{\kappa} @>s_{\kappa}>>  \tilde{\T}^{sc}_{\kappa} @>>>  \tilde{\T}_{\kappa} =  \T_{\kappa}  \times \G_m
\end{CD} \] 
from (a). Suppose that 
\[  f_{\kappa} \circ \alpha^{\vee} : t  \mapsto   ( \alpha^{\vee}(t),  t^{n_{\alpha}}) \quad \text{with $n_{\alpha} \in \Z$ and for $t \in \G_m$.} \] 
Then we need to show that 
\[ n_{\alpha}  =  {\rm ord}(\eta(\alpha^{\vee}))  \quad \text{for all $\alpha \in \Delta$.}.  \]
\vskip 5pt

\item[(d)]  Now the splitting $s_{\kappa}:  \G^{sc}_{\kappa} \longrightarrow \tilde{\G}^{sc}_{\kappa}$ from (a) is uniquely determined by its restriction to the root subgroups $\U_{\alpha, \kappa}$ for $\alpha \in \Delta$. Hence $s$ is determined by $s_{\kappa} \circ x_{\alpha,\kappa}$, where $x_{\alpha}  : \G_a \rightarrow \U_{\alpha}$ is part of the fixed \'epinglage.    
Since 
\[  \alpha^{\vee}(t)  = n_{\alpha}(t)  \cdot n_{\alpha}(-1) \in \T^{sc}_{\kappa}  \subset \T_{\kappa}  \quad  \text{with} \quad n_{\alpha}(t)  = x_{\alpha}(t) \cdot x_{-\alpha}(-t^{-1}) \cdot x_{\alpha}(t), \] 
this implies that 
\[  (\alpha^{\vee}(t),  t^{n_{\alpha}}) =  f_{\kappa} \circ  \alpha^{\vee}(t)  =  \text{image of $n_{\alpha}(t)  \cdot n_{\alpha}(-1)$  in $\tilde{\T}_{\kappa}$.} \]
\vskip 5pt

\item[(e)] Likewise,
the induced system of splittings $s_{\mathcal{O'}} : \G^{sc}(\mathcal{O'}) \rightarrow\overline{\G}^{sc}(\mathcal{O'})$ is   determined by the unique splitting
\[  \tilde{x}_{\alpha}:  F'  \longrightarrow \overline{G}_{F'}^{{\rm tame}} \quad \text{for all $\alpha \in \Delta$.}  \]
This implies that the composite
\[  \begin{CD}
  f_{\mathcal{O'}}  :  \T^{sc}(\mathcal{O'}) @>s_{\mathcal{O'}}>>  \overline{\T}^{sc}(\mathcal{O'}) @>>>  \overline{\T}(\mathcal{O'})  \end{CD} \]
  from (a) is given by
\[  f_{\mathcal{O'}} \circ \alpha^{\vee}(\tilde{t}) = \text{image of $\tilde{n}_{\alpha}(\tilde{t}) \cdot \tilde{n}_{\alpha}(-1)$ in $\overline{T}(\mathcal{O'})$}. \]
  The RHS is nothing but the section
 \[  s_{\alpha}(\tilde{t}) = (\alpha^{\vee}(\tilde{t}), (\eta(\alpha^{\vee}), \tilde{t})_n) \in \T(\mathcal{O'})  \]
for $\tilde{t} \in \mathcal{O'}^{\times}$, whose image under the reduction map is
\[ (\alpha^{\vee}(t) ,  t^{{\rm ord}(\eta(\alpha^{\vee}))}) \in  \tilde{T}_{\kappa}(\kappa') = \T_{\kappa}(\kappa')  \times {\kappa'}^{\times}. \]
By (a),  
\[      f_{\kappa} \circ  \alpha^{\vee}(t)  = \text{the image of $ f_{\mathcal{O'}} \circ \alpha^{\vee}(\tilde{t})$ under reduction map.} \]
Hence,  it follows   that
 \[ n_{\alpha}  =  {\rm ord}(\eta(\alpha^{\vee}))  \]
 for $\alpha \in \Delta$, as desired.
\end{itemize}

\end{proof}
\vskip 5pt

In terms of the bisector data $(D,\eta)$, one has $\mathcal{E}  = Y \times_D  F^{\times}$ and 
$f(s_{Q^{sc}}(\alpha^{\vee}))  =  (\alpha^{\vee}, \eta(\alpha^{\vee}))$ for $\alpha \in \Delta$.  On pushing out by the valuation map, one has:
\[  \mathcal{E}_{\kappa}  = Y\times   \Z  \quad \text{(direct product of groups)} \]
and
\[  f_{\kappa}(\alpha^{\vee})  = (\alpha^{\vee}, {\rm ord}(\eta(\alpha^{\vee}))) \quad \text{for $\alpha \in \Delta$.} \]
In particular,  one has:
\vskip 5pt

\begin{cor}
 If $Y/Y^{sc}$ is free or if $\eta$ takes value in $\mathcal{O}_F^{\times}$, then the algebraic extension $\overline{\G}_{\kappa}$ is split. Thus, the topological central extension $\overline{K}$  of $K$ is also split. 
\end{cor}

\vskip 5pt

We have assumed that $n  = \# \kappa^{\times}$ above. In general, when $p$ does not divide $n$, the 
$n$-th Hilbert symbol map $K_2(F)  \longrightarrow \mu_n$ factors through
$K_2(F)  \longrightarrow \kappa^{\times} \longrightarrow \mu_n$. 
So the degree $n$ BD covering $\overline{\G}$ is obtained from the one of degree $\# \kappa^{\times}$ as a pushout. In particular, when the conditions of the above corollary holds, the degree $n$ cover  $\overline{K}$ is split as well. Indeed, whenever $\eta$ takes value in $\mathcal{O}^{\times} \cdot F^{\times n}$,  the cover $\overline{K}$ splits.  

\vskip 5pt
Note that we have merely given some simple sufficient conditions for $\overline{K}$ to be split.  These conditions may not be necessary in a given case, but as we will see below, it is possible for $\overline{K}$ to be non-split when they fail.  
Moreover, note that the splitting of $\overline{K}$ is not necessarily unique (if it exists). 
\vskip 10pt

\subsection{\bf Running example.}  We illustrate the discussion in this section with our running example: $\G = \PGL_2$, $D=0$ and $n=2$.  Then we have the BD extensions 
\[  \overline{\G}_{\eta}  = (\GL_2 \times K_2) /  \{  (z,  \{ \eta, z \}): z \in \G_m  \}. \]  
The associated BD covering groups are:
\[  \overline{G}_{\eta}  = (\GL_2(F)  \times \mu_2) /  \{  (z,  (\eta, z)_2):  z \in F^{\times}  \}. \]
Let $\pi_{\eta} :  \GL_2(F) \times \mu_2 \longrightarrow \overline{G}_{\eta}$ be the natural projection map, and let 
\[  A =  \{  \left( \begin{array}{cc}
a & 0 \\
0 & 1  \end{array}  \right) :  a \in F^{\times} \} \subset \GL_2(F) . \]
The projection map identifies $A$ with a maximal split torus $T$ of $\PGL_2(F)$ and 
$\pi_{\eta}$ identifies $A \times \mu_2$ with  $\overline{T}$ of $\PGL_2(F)$. In this case, $\overline{T}$ is abelian (since $D= 0$) and so $\overline{T} =  \overline{T}_{Q,2}$ and $f: T_{Q,2} \longrightarrow T$ is the identity map.
\vskip 5pt

Now consider the issue of whether the covering splits over $K= \PGL_2(\mathcal{O})$.  We have already seen from general arguments that it does when $\eta \in \mathcal{O}^{\times} \cdot F^{\times 2}$. When $\eta = \varpi$ is a uniformizer, we shall show now that the covering $\overline{K}_{\eta}$ is non split.
\vskip 5pt

If a splitting $K \longrightarrow \overline{K}_{\eta}$ exists, we would have a group homomorphism
\[  \phi:  \GL_2(\mathcal{O})  \longrightarrow (\GL_2(F)  \times \mu_2) /  \{  (z,  (\eta, z)_2):  z \in F^{\times}  \} \]
which is trivial on the center $Z(\mathcal{O})$ of $\GL_2(\mathcal{O})$. 
For $k \in \GL_2(\mathcal{O})$, we may write 
\[  \phi(k)  =  \text{the class of $(k,  \mu(k))$}  \]
for some $\mu(k)   = \pm 1$. Now it is easy to check that $\mu:  \GL_2(\mathcal{O}) \longrightarrow \mathcal{O}^{\times}/ \mathcal{O}^{\times 2} = \{\pm 1\}$ is a group homomorphism and thus $\mu$ factors as
\[ \begin{CD}
 \mu:  \GL_2(\mathcal{O}) @>\det>> \mathcal{O}^{\times} @>>> \kappa^{\times} @>>> {\pm 1}.\end{CD} \]
 If now $k  = z \in Z(\mathcal{O})$ is a scalar matrix, then the fact that $\phi(z)$ is trivial means that
 \[  \mu(z)  = (\varpi, z)_2.  \]
Since $\mu$ factors through $\det$, we see  that $\mu(z)  =1$, but $(\varpi, z)_2$ is not $1$ for some $z \in \mathcal{O}^{\times}$. With this contradiction, we see that the covering $\overline{K}_{\eta}$ is not
split when $\eta = \varpi$ is a uniformizer. 

\vskip 10pt
\section{\bf Dual and $L$-Groups}
In this section, we shall  recall the definition of the L-group ${^L}\overline{G}$ of a BD covering $\overline{G}$ for a split $\G$ over a local field, following Weissman \cite{W3, W4}.
The construction in \cite{W3} is quite involved, using a double twisting of the Hopf algebra of a candidate dual group. In a letter to Deligne \cite{W4}, Weissman gave an retreatment of his construction in \cite{W3}, using only the BD data $(Q, \mathcal{E}, f, n)$. We shall follow this more streamlined treatment in \cite{W4}.  

\vskip 5pt
Since we will simply be presenting the construction of these objects in this section, 
 the definition of the dual group or L-group of $\overline{G}$ may seem rather unmotivated at the end of the section.  Whether they are the right objects or not will largely depend on whether they give the right framework to describe the representation theory of $\overline{G}$.  In the subsequent sections, we will address these concerns by discussing various reality checks.

\vskip 5pt

\subsection{\bf Dual group \`a la Finkelberg-Lysenko-McNamara-Reich.}
Let $\G$ be a split connected reductive group over $F$, with maximal split torus $\T$ and cocharacter lattice $Y$.
Let $\Phi^{\vee} \subset Y$ be the set of coroots of $\G$ and let $Y^{sc} \subset Y$ be the sublattice generated by $\Phi^{\vee}$. Similarly, let $\Phi \subset X$ be the set of roots  generating a sublattice $X^{sc}$ in the character lattice $X$ of $\T$.
\vskip 5pt

 Suppose that $\overline{\G}$ is a multiplicative $K_2$-torsor with associated 
a BD data $(Q, \mathcal{E}, \phi)$. The data $(Q, \mathcal{E}, \phi,n)$ (with $|\mu_n(F)|  = n$) then
 gives  a central extension of   locally compact groups
\[  \begin{CD}
 1@>>>  \mu_n(F)  @>>>  \overline{G}  @>>> G @>>> 1 \end{CD} \]
 \vskip 5pt
 
Using the data $(Y, \Phi^{\vee}, Q,n)$, we may define a modified root datum
as follows:
\vskip 5pt

\begin{itemize}
\item we have already set
\[  Y_{Q,n}  = Y \cap n Y^*  \]
where $Y^* \subset Y \otimes_{\Z}  \Q$ is the dual lattice of $Y$ relative to $B_Q$. 
Let $X_{Q,n} \subset X \otimes_{\Z}  \Q$ be the dual lattice to $Y_{Q,n}$.
\vskip 5pt

\item for each $\alpha^{\vee} \in \Phi^{\vee}$, set
\[   n_{\alpha}  =   \frac{n}{\text{gcd}(n, Q(\alpha^{\vee}))}, \]
and 
\[   \alpha_{Q,n}^{\vee}  = n_{\alpha}  \cdot \alpha. \]
Denote by $\Phi_{Q,n}^{\vee}$ the set of such $\alpha_{Q,n}^{\vee}$'s and observe that 
\[  \Phi_{Q,n}^{\vee}  \subset Y_{Q,n}. \]
We let $Y_{Q,n}^{sc}$ denote the sublattice of $Y_{Q,n}$ generated by  $\Phi_{Q,n}^{\vee}$.
\vskip 5pt

\item likewise, for $\alpha \in \Phi$,  set
\[   \alpha_{Q,n}  = n_{\alpha}^{-1} \cdot \alpha \]
 and denote by $\Phi_{Q,n}$ the set of such $\alpha_{Q,n}$'s, so that $\Phi_{Q,n}  \subset X_{Q,n}$.
 \end{itemize}
 \vskip 5pt
 
 \noindent  Then it was shown in \cite{Mc2} and \cite{W3} that the quadruple $(Y_{Q,n}, \Phi_{Q,n}^{\vee},  X_{Q,n},  \Phi_{Q,n})$
 is a root datum, and hence determine a split connected reductive group $\overline{G}^{\vee}$ over $\C$. 
 {\em The group $\overline{G}^{\vee}$  is by definition the dual group of the BD extension $\overline{G}$}. 
 Observe that it only depends on $(Q, n)$ and is independent of the third ingredient $f$ of a BD data $(Q, \mathcal{E}, f)$; equivalently, it only depends on $(D, n)$ but not on $\eta$.  
\vskip 5pt

Let $Z(\overline{G}^{\vee})$ be the center of $\overline{G}^{\vee}$. Then note that
\[ Z(\overline{G}^{\vee})  = \Hom_{\Z}(Y_{Q,n}/  Y_{Q,n}^{sc}, \mathbb{C}^{\times}). \]
 \vskip 5pt

\subsection{\bf L-group \`a la Weissman}
We can now describe Weissman's proposal for the  L-group of $\overline{G}$.  
This is done by defining an extension
\[ \begin{CD}   Z(\overline{G}^{\vee}) @>>> E  @>>>  F^{\times}/ F^{\times n} \end{CD} \]
followed by pushing  out by the natural inclusion $ Z(\overline{G}^{\vee})  \rightarrow \overline{G}^{\vee}$ and pulling back via the natural projection $W_F \rightarrow F^{\times}/F^{\times n}$. This results in an extension
\[  \begin{CD}
1 @>>>  \overline{G}^{\vee}  @>>>  {^L}\overline{G}  @>>>  W_F @>>> 1, \end{CD} \]
which we call Weissman's L-group extension. 
\vskip 5pt

The construction of $E$ is as a Baer sum $E_1 +  E_2$ of two extensions $E_1$ and $E_2$. These are defined as follows:
\vskip 5pt

\begin{itemize}
\item $E_1$ is defined explicitly using the cocycle
\[   c_1:  F^{\times}/ F^{\times n}  \times F^{\times}/ F^{\times n}  \longrightarrow Z(\overline{G}^{\vee})  = \Hom_{\Z}(Y_{Q,n}/Y_{Q,n}^{sc},  \C^{\times}) \]
given by
\[  c_1(a,b) (y)  =  (a,b)_n^{Q(y)}. \]
Since $2 \cdot Q(y) = B_Q(y,y) \in n\Z$ for $y \in Y_{Q,n}$, we see that this cocycle is trivial when $n$ is odd, and is valued in $\pm 1$ when $n$ is even.  Note that $E_1$ depends only on $(Q,n)$. 
\vskip 5pt

\item the construction of $E_2$ is slightly more involved and uses the full BD data $(Q, \mathcal{E},f)$, where we recall that $\mathcal{E}$ is an extension
\[ \begin{CD}   F^{\times}  @>>> \mathcal{E}  @>>> Y, \end{CD} \]
and
\[  f:  \mathcal{E}_{Q^{sc}}  \longrightarrow q^*(\mathcal{E}) \]
is an isomorphism, with $q:  G^{sc}  \rightarrow G$ the natural map. 
\vskip 5pt

Since we have the inclusion $Y_{Q,n}^{sc}  \rightarrow Y_{Q,n}  \rightarrow Y$, we may pullback 
the extensions $\mathcal{E}_{Q^{sc}}$ and   $\mathcal{E}$ and pushout via $F^{\times} \longrightarrow F^{\times}/F^{\times n}$ to obtain
\[   \begin{CD}  
F^{\times}/F^{\times n}  @>>> \mathcal{E}_{Q^{sc},n}  @>>>  Y^{sc}_{Q,n}   \\
@|  @VVfV  @VVV  \\
F^{\times}/F^{\times n}  @>>> \mathcal{E}_{Q,n}  @>>> Y_{Q,n}  
 \end{CD}  \]
 Note that both $\mathcal{E}_{Q,n}$ and $\mathcal{E}_{Q^{sc},n}$ are abelian groups.
 \vskip 5pt
 
 For each $\alpha^{\vee} \in \Phi^{\vee} \subset Y^{sc}$, we have defined before 
an element $s_{Q^{sc}}(\alpha^{\vee})  \in \mathcal{E}_{Q^{sc}}$ lying over $\alpha^{\vee}$. Indeed, 
 $s_{Q^{sc}}(\alpha^{\vee})$ is the image of the element $s_{\alpha}(\alpha^{\vee}(\tau)) \in \overline{\T}^{sc}(F(\!(\tau)\!))$ under pushout by the residue map ${\rm Res}: K_2(F(\!(\tau)\!)) \longrightarrow F^{\times}$. Analogously, we 
have  the element $s_{Q^{sc}}(n_{\alpha} \cdot \alpha^{\vee}) \in \mathcal{E}_{Q^{sc}}$ 
which is the image of the element 
\[ s_{\alpha}((n_{\alpha} \cdot \alpha^{\vee})(\tau)) = s_{\alpha}(\alpha^{\vee}(\tau^{n_{\alpha}})) \in \overline{\T}^{sc}(F(\!(\tau)\!)). \]
It  lies over $\alpha_{Q,n}  = n_{\alpha} \alpha^{\vee} \in Y_{Q, n}^{sc}$.
Weissman showed that this induces a group homomorphism
\[  s_{Q^{sc}}:  Y_{Q, n}^{sc}   \longrightarrow   \mathcal{E}_{Q^{sc},n}. \]
Composing this with $f$, one obtains
\[  s_f =  f \circ s_{Q^{sc}}:   Y_{Q, n}^{sc}  \longrightarrow \mathcal{E}_{Q,n}. \]

\vskip 5pt

Viewing $Y_{Q,n}^{sc}$ as a subgroup of $\mathcal{E}_{Q,n}$ by the splitting $s_f$, we inherit an extension 
\begin{equation}\label{E:E2} \begin{CD} 
 F^{\times}/F^{\times n}  @>>> \overline{\mathcal{E}}_{Q,n} =   \mathcal{E}_{Q,n}/ s_f(Y_{Q,n}^{sc}) @>>> Y_{Q,n}/ Y_{Q,n}^{sc} \\
  \end{CD} \end{equation}
 Taking $\Hom_{\Z}( -, \C^{\times})$ (which is exact, since $\C^{\times}$ is divisible and hence injective), we obtain the desired extension:
 \begin{equation} \label{E:E2'} 
 \begin{CD} 
 Z(\overline{G}^{\vee}) @>>>  E_2  @>>>   \Hom_{\Z}( F^{\times}/F^{\times n},\C^{\times}) \cong  F^{\times}/F^{\times n} \end{CD} \end{equation}
where the last isomorphism is via  the $n$-th Hilbert symbol:  $a \in  F^{\times}/F^{\times n}$ giving rise to the character $\chi_a: b \mapsto (b,a)_n$.
\end{itemize}
 \vskip 5pt

\subsection{\bf Description using bisectors} 
We may describe the construction of $E_2$ in terms of the bisector data $(D, \eta)$. 
The bisector $D$ allows us to realise the extension $\mathcal{E}$ as a set $Y \times F^{\times}$ with group law
\[  (y_1,a) \cdot (y_2,b)  = (y_1 + y_2,  ab\cdot (-1)^{D(y_1, y_2)}). \]
Pushing this out by $F^{\times}  \rightarrow F^{\times}/ F^{\times n}$ and pulling back to $Y_{Q,n}$
gives the extension $\mathcal{E}_{Q,n} =  Y_{Q,n}  \times F^{\times}/ F^{\times n}$ with 
  the same group law as above. In particular, if $n$ is odd, $-1 \in F^{\times n}$ so that the cocycle $(-1)^{D(y_1, y_2)}$ is trivial.
\vskip 5pt

 The map  $f:  \mathcal{E}_{Q^{sc}} \longrightarrow \mathcal{E}$ is defined by 
 \[  f(s_{Q^{sc}}(\alpha^{\vee}))  = (\alpha^{\vee} , \eta(\alpha^{\vee})) \in Y \times F^{\times},  \quad \text{for $\alpha \in \Delta$.} \]
 Then  the splitting $s_f$ is given by
 \[  s_f(\alpha_{Q,n}^{\vee})  =  (\alpha_{Q,n}^{\vee}, \eta(\alpha_{Q,n}^{\vee})) \in Y_{Q,n} \times F^{\times},  \quad \text{  for $\alpha \in \Delta$.}  \]
\vskip 5pt

It is instructive to note that the above constructions are functorial in nature.  Given any isomorphism $\xi :  (D, \eta)  \longrightarrow (D', \eta')$,  $\xi$ carries the map $s_f$ corresponding to $(D,\eta)$ to the map $s_{f'}$ corresponding to $(D', \eta')$.


\vskip 5pt

\subsection{\bf Running example.}  Again, we illustrate the discussion in this section using our running example: $G = \PGL_2(F)$, $Q=0$ and $n=2$. 
In this case,  $Y = \Z = Y_{Q,n}$ and $Y^{sc}  = 2\Z = Y_{Q,n}^{sc}$. So 
\[  Z(\overline{G}_{\eta}^{\vee})  = \mu_2 \subset \overline{G}_{\eta}^{\vee}  = \SL_2(\C) \quad \text{for any $\eta$.} \]
Moreover, $E_1^{\eta}$ is the split extension $Z(\overline{G}_{\eta}^{\vee}) \times F^{\times}/ F^{\times 2}$ 
and $\mathcal{E}_{\eta}  = Y  \times F^{\times}/F^{\times 2}$ is split.
 Hence
 \[  \overline{\mathcal{E}}_{\eta}  =( \Z \times F^{\times}/F^{\times 2} ) /  \{  (2y,   \eta^y ): y \in \Z \} \]
 and 
\[  E^{\eta}_1+ E^{\eta}_2  =  \Hom(  \overline{\mathcal{E}}_{\eta} ,  \C^{\times})
= \{  ( t, a) \in \C^{\times} \times F^{\times}/F^{\times 2}:  t^2 = (\eta, a)_2 \}. \]
This contains $Z(\overline{G}_{\eta}^{\vee})  = \mu_2$ as the subgroup of elements $(\pm 1, 1)$, and the associated quotient is via the second projection to   $F^{\times}/F^{\times 2}$.  
 \vskip 5pt
 
 Observe that when $F = \R$ and $\eta = -1 \in \R^{\times}$, then $E_1^{\eta} +  E^{\eta}_2$ is the cyclic group $\mu_4$ and so the above extension is non-split! However, when we push out via the natural map $Z(\overline{G}_{\eta}^{\vee})  = \mu_2 \hookrightarrow \overline{T}_{\eta}^{\vee}  = \C^{\times}$, then the pushout sequence is split. Indeed the sequence splits once we pushout by $\mu_2 \hookrightarrow \mu_4$. 
 
 \vskip 5pt
 
 For general local field $F$, one sees that when one pushes $ \overline{\mathcal{E}}_{\eta}$ out by $\mu_2 \hookrightarrow \SL_2(\C)$, one obtains:
 \[  \begin{CD}
 \SL_2(\C)  @>>> \{  (g, a) \in  \GL_2(\C) \times F^{\times}/F^{\times 2}:  \det(g)  = (\eta, a)_2 \} @>>> F^{\times}/F^{\times 2} \end{CD} \]
 Pulling back to $W_F$, one obtains:
 \[  {^L}\overline{G}_{\eta}=  \{  (g,w) \in \GL_2(\C) \times W_F:  \det(g)  = \chi_{\eta}(w)\} \cong
 \SL_2(\C)  \rtimes_{\eta} W_F \]
 where $w \in W_F$ acts on $\SL_2(\C)$ by the conjugation action of the diagonal matrix ${\rm diag}(\chi_{\eta}(w), 1)$. Thus, while the L-group extension is always split,  ${^L}\overline{G}_{\eta}$ is not isomorphic to the direct product $\SL_2(\C) \times W_F$ for general $\eta$. 
 \vskip 5pt

 \subsection{\bf Functoriality for Levis.}  \label{SS:levi}
Suppose that $\mathbb{M} \subset \G$ is a proper Levi subgroup, then a BD-covering $\overline{G}$ restricts to one on $M$. If the bisector data for $\overline{G}$ is $(D, \eta, n)$, then that for $\overline{M}$ is $(D,  \eta|_{Y_{\mathbb{M}}^{sc}}, n)$, where we have restricted $\eta$ to the sublattice $Y^{sc}_{\mathbb{M}}$  generated by the simple coroots of $M$ in $Y$. The 
above construction produces L-groups extensions ${^L}\overline{M}$ and ${^L}\overline{G}$.  
An examination of the construction shows:
\vskip 5pt

\begin{lemma}
 This is  a natural commutative diagram of short exact sequences:
\[  \begin{CD}
Z(\overline{G}^{\vee})  @>>>  E_G  @>>>  W_F  \\
@VVV   @VVV  @|  \\
Z(\overline{M}^{\vee})  @>>>  E_M  @>>>  W_F \\ 
@VVV  @VVV  @| \\
\overline{M}^{\vee}  @>>>  {^L}\overline{M}  @>>> W_F  \\
@VVV  @VVV  @|  \\
\overline{G}^{\vee}  @>>>  {^L}\overline{G} @>>>  W_F  \end{CD} \]
\end{lemma}
 \vskip 5pt
\begin{proof}
 It suffices to exhibit a natural map from the extension $E_G$ to $E_M$, which gives:
\[  E_M  \cong  \frac{Z(\overline{M}^{\vee} \times E_G}{\Delta Z(\overline{G}^{\vee})} \]
and thus  induces a map 
\[  {^L}\overline{M}  :=\frac{ \overline{M}^{\vee} \times E_M}{\Delta Z(\overline{M}^{\vee})}\cong \frac{\overline{M}^{\vee} \times E_G}{\Delta Z(\overline{G}^{\vee})}  \longrightarrow \frac{\overline{G}^{\vee}  \times E_G}{\Delta Z(\overline{G}^{\vee})} =:  {^L}\overline{G} \]
making the above diagram commute. 
\vskip 5pt

Write $E_G=E_1 + E_2$ and $E_M=E_{M,1} + E_{M,2}$ as  Baer sums. Let $Y_{M,Q,n}^{sc} \subseteq Y_{\mathbb{M}}^{sc} \cap Y_{Q,n}^{sc}$ be the sublattice of $Y_{\mathbb{M}}^{sc}$ generated by $\alpha_{Q,n}^\vee$ for $\alpha^\vee \in Y_{\mathbb{M}}^{sc}$. 
The cocycle defining $E_1$ takes value in $\text{Hom}_\Z(Y_{Q,n}/Y_{Q,n}^{sc}, \C^\times)$, and the cocycle for $E_{M,1}$ takes the same formula and is valued in $\text{Hom}_\Z(Y_{Q,n}/Y_{M,Q,n}^{sc}, \C^\times)$. Thus there is a natural map from $E_1$ to $E_{M,1}$, which in fact is canonically isomorphic to the push-out of $E_1$.

Consider $E_2$ and $E_{M,2}$. The pull-back of $\mathcal{E}_{Q^{sc},n}$ via $Y_{M,Q,n}^{sc} \subseteq Y_{Q,n}^{sc}$ gives an extension 
\[
\begin{CD}
F^\times/F^{\times n}  @>>>  \mathcal{E}_{M, Q^{sc},n}  @>>>   Y_{M,Q,n}^{sc}
\end{CD}
\]
which has the splitting $s_{M, Q^{sc}}$ from the restriction of $s_{Q^{sc}}$ to $Y_{M,Q,n}^{sc}$. Let $s_{M, f}:=f\circ s_{M, Q^{sc}}$ and $\overline{\mathcal{E}}_{M, Q, n}:=\mathcal{E}_{M, Q^{sc},n}/s_{M,f}(Y_{M,Q,n}^{sc})$. Then one obtains the following commutative diagram
\[  \begin{CD}
F^\times/F^{\times n}  @>>>  \overline{\mathcal{E}}_{Q, n}    @>>>   Y_{Q,n}/Y_{Q,n}^{sc}    \\
@AAA   @AAA  @AAA  \\
F^\times/F^{\times n}  @>>>  \overline{\mathcal{E}}_{M, Q, n}  @>>>   Y_{Q,n}/Y_{M,Q,n}^{sc}. \end{CD} \]

Applying $\text{Hom}_\Z(-, \C^\times)$ and followed by the identification $\text{Hom}_\Z(F^\times/F^{\times n}, \C^\times) \simeq F^\times/F^{\times n}$ via the $n$-th Hilbert symbol, we obtain a natural map from $E_2$ to $E_{M,2}$ as in
\[  \begin{CD}
Z(\overline{G}^\vee)  @>>>  E_2  @>>>   F^\times/F^{\times n}  \\
@VVV   @VVV  @|  \\
Z(\overline{M}^\vee)  @>>>  E_{M,2}  @>>>   F^\times/F^{\times n}   . \end{CD} \]
Combining the two extensions for the Baer sum, we see that there is a natural map from $E_G$ to $E_M$. This gives a natural isomorphism
\[  Z(\overline{M}^{\vee}) \times E_G /\Delta Z(\overline{G}^{\vee})  \longrightarrow E_M. \]
\end{proof}
\vskip 5pt

\subsection{\bf Functoriality for $z$-extensions}  \label{SS:z-ext}
If 
\[  \begin{CD}
Z @>>> \mathbb{H}  @>>> \G  \end{CD} \]
is a $z$-extension, and $\overline{G}$ is a BD covering with bisector data $(D,\eta)$, then we obtain a BD covering $\overline{H}$ with essentially the same bisector data, such that
\[  \begin{CD}
Z @>i>>  \overline{H}  @>>>  \overline{G}. \end{CD} \]
From the construction of the L-group extension, one deduces:
\vskip 5pt

\begin{lemma}  \label{L:z-ext}
There is a commutative diagram 
\[  \begin{CD}
\overline{G}^{\vee}  @>>> {^L}\overline{G}  @>>>  W_F  \\
@VVV  @VVV  @|  \\
\overline{H}^{\vee}  @>>> {^L}\overline{H} @>>>  W_F \\
@VVV   @VVV     \\
Z^{\vee}  @=  Z^{\vee}   \end{CD} \]
\end{lemma}

\vskip 10pt
\section{\bf Distinguished Splittings of L-Groups} \label{S:d-splitting}
In this section, we study the L-group extension proposed in the previous section. More precisely, we will investigate whether this extension actually splits. Since the L-group extension is defined from the extension $E_1 + E_2$, it is natural to consider whether $E_1 + E_2$ splits. We have seen in the last section that it does not in general, but we would like to understand where the obstruction lies. We shall assume that $\overline{\G}$ is defined by a pair $(D, \eta)$ where $D$ is a {\em fair} bisector. 
\vskip 5pt

\subsection{\bf Splittings of $E_1+E_2$}  \label{SS:splitting}
 What does it mean to give a splitting of $E_1 + E_2$?  
Since $E_1$ is given by an explicit cocycle $c_1$ (valued in $\pm 1$), the issue of equipping $E_1+E_2$ with a splitting is equivalent to  finding a set theoretic section of $E_2$ whose associated cocycle  is equal to $c_1$.
In addressing this question, we shall work explicitly using a bisector data $(D,\eta)$. 
\vskip 5pt

In the construction of $E_2$ above, the extension (\ref{E:E2}) is equipped with a cocycle inherited from the fair  bisector $D$. However, in taking $\Hom_{\Z}(-, \C^{\times})$ to obtain (\ref{E:E2'}), we have partially lost the data of a cocycle.  To obtain a cocycle, we need to choose a set theoretic section of (\ref{E:E2'}).  Hence, 
for each $a \in F^{\times}/ F^{\times n}$, we need to extend the  character
\[  \chi_a \in \Hom_{\Z}(F^{\times}/F^{\times n}, \C^{\times}) \]
  to a character $\tilde{\chi}_a$ of $\overline{\mathcal{E}}_{Q,n}$. Another way of expressing this is to consider the pushout of $\mathcal{E}_{Q,n}$  by $\chi_a$:
\[  \begin{CD}
F^{\times}/F^{\times n} @>>> \mathcal{E}_{Q,n} @>>>  Y_{Q,n}   \\
@VV\chi_aV   @VVV   @VVV   \\
\mu_n  @>>> \tilde{\mathcal{E}}_{Q,n}  @>>>  Y_{Q,n}.  \end{CD} \]
Then extending $\chi_a$ amounts to finding a genuine character of 
 $\tilde{\mathcal{E}}_{Q,n}$ which is trivial when precomposed with $f \circ s$. 
 
 \vskip 5pt
 
 Note that  $\tilde{\mathcal{E}}_{Q,n}  =  Y_{Q,n}  \times \mu_n$ has group law
 \[  (y_1, 1)  \cdot (y_2, 1)  
=  ( y_1+ y_2,  \chi_a(-1)^{D(y_1, y_2)})
= (y_1 +y_2,  (a,a)_n^{D(y_1, y_2)}). \]
Thus to give a genuine character $\tilde{\chi}_a$ of $\tilde{\mathcal{E}}_{Q,n}$ is to give a function
\[  \tilde{\chi}_a :   Y_{Q,n}   \longrightarrow \C^{\times}    \]
satisfying
\[  \tilde{\chi}_a(y_1)  \cdot \tilde{\chi}_a(y_2)  = \tilde{\chi}_a(y_1+y_2) \cdot (a,a)_n^{D(y_1, y_2)}. \]
For the composite $\tilde{\chi}_a \circ f \circ s$ to be trivial, one needs to  require that
\[  \tilde{\chi}_a(\alpha_{Q,n}^{\vee})  =(a, \eta(\alpha_{Q,n}^{\vee}))_n \quad \text{  for all $\alpha \in \Delta$.} \]
\vskip 5pt

Let us fix such a genuine character $\tilde{\chi}_a$ for each $a \in F^{\times}/ F^{\times n}$. Using this as a set theoretic section for (\ref{E:E2'}), we may write
\[  E_2  =     F^{\times}/ F^{\times n} \times  \Hom_{\Z}(Y_{Q,n}, \C^{\times})    \]
with cocycle
\[  c_2(a,b)(y)  =  \tilde{\chi}_a (y) \tilde{\chi}_b(y)  \tilde{\chi}_{ab}(y)^{-1} \quad \text{for all $y \in Y_{Q,n}$}.  \] 
To impose the condition that $c_2  =  c_1$ means:
\[   \tilde{\chi}_a (y) \tilde{\chi}_b(y)    =   \tilde{\chi}_{ab}(y) \cdot (a,b)_n^{Q(y)}  \]
for all $y \in Y_{Q,n}$. 
\vskip 10pt

To summarise, we see that  to obtain a splitting of $E_1 + E_2$, it is equivalent to give a function
\[  \chi:  F^{\times}   \times Y_{Q,n}   \longrightarrow \C^{\times} \]
such that
\begin{itemize}
\item[(a)]   
\[  \chi(ab^n, y)  = \chi(a, y) \quad \text{ for $a,b \in F^{\times}$ and $y \in Y_{Q,n}$.} \]
\item[(b)] 
\[  \chi(a,y_1) \cdot \chi(a, y_2) =  \chi(a, y_1+y_2) \cdot (a,a)_n^{D(y_1, y_2)}  \]
\item[(c)]
\[  \chi(a, y) \cdot \chi(b,y)    =  \chi(ab, y) \cdot (a,b)_n^{Q(y)}. \]
\item[(d)]  
\[  \chi(a, \alpha_{Q,n}^{\vee})  = (a, \eta(\alpha_{Q,n}^{\vee}))_n  \quad \text{ for $a \in F^{\times}$ and $\alpha \in \Delta$.} \]
\end{itemize}
\vskip 10pt

In \S \ref{SS:coveringT} and \S  \ref{SS:coveringTQn},  we have described the covering torus $\overline{T}_{Q,n}$ by generators and relations, using the bisector $D$.  It follows 
immediately from these that to give a splitting of $E_1 + E_2$ is equivalent to giving a genuine character $\chi$  of $\overline{T}_{Q,n}$ (by properties (b) and (c)) satisfying some properties (dictated by properties (a) and (d)). 
More precisely, properties (a) and (d) can be rephrased as:
\vskip 5pt

\begin{itemize}
\item[(a')]  the inclusion $nY_{Q,n}  \hookrightarrow Y_{Q,n}$ gives the $n$-power isogeny 
$n: T_{Q,n}  \rightarrow T_{Q,n}$ and this lifts to a group homomoprhism 
\[   i_n:  T_{Q,n}  \longrightarrow \overline{T}_{Q,n} \]
given by
\[  i_n( y(a))  =  (ny)(a)  = y(a^n) \quad y \in Y_{Q,n}. \] 
 Then property (a) says that $\chi \circ i_n$ is trivial.
 \vskip 5pt
 
 \item[(d')]  the inclusion $Y_{Q,n}^{sc} \hookrightarrow Y_{Q,n}$ induces  
an isogeny $T_{Q,n}^{sc}  \rightarrow T_{Q,n}$. and the map
\[  \alpha_{Q,n}^{\vee}(a)  \mapsto \alpha_{Q,n}^{\vee}(a) \cdot  (\eta(\alpha_{Q,n}^{\vee}), a)_n   \in \overline{T}_{Q,n},  \quad \text{for all $\alpha \in \Delta$,} \]
  defines  a splitting of this isogeny to give $s_{\eta}:  T_{Q,n}^{sc} \longrightarrow \overline{T}_{Q,n}$.    Then property (d) says that $\chi \circ s_{\eta}$ is trivial.  
  Because $D$ is fair, we have in fact, for any $y \in Y_{Q,n}^{sc}$,
  \[  s_{\eta}(y(a)) =  y(a) \cdot (\eta(y), a)_n  \in \overline{T}_{Q,n}. \]
  \end{itemize}
\vskip 5pt

\subsection{\bf Obstruction.}
The question is thus: does there exist genuine characters of $\overline{T}_{Q,n}$ such that (a') and (d') are satisfied?  We shall see that there will be some obstructions. More precisely, 
suppose that  $y(a) \in T_{Q,n}^{sc}$ ($y \in Y_{Q,n}^{sc}$) belongs to 
\[  {\rm Ker}(T_{Q,n}^{sc} \rightarrow T_{Q,n}) = {\rm Tor}_1(Y_{Q,n}/ Y_{Q,n}^{sc}, F^{\times}). \]
Note that this kernel is generated by such pure tensors $y(a)$.  
 Then (d') requires $\chi$ to be trivial on the element $(\eta(y), a)_n \in \mu_n(F)$. But $\chi$
is supposed to be a genuine character. So we have our first obstruction:
\vskip 5pt

\noindent{\bf Obstruction 1:}
A genuine character  $\chi$ satisfying (d')  exists if and only if 
\begin{equation} \label{E:Ob1}
 (\eta(y),a)_n =1 \quad \text{whenever $y \otimes a =0$ in $Y_{Q,n} \otimes F^{\times}$ with $y \in Y_{Q,n}^{sc}$.} \end{equation}
\vskip 5pt

\noindent This condition does not hold automatically. It does hold, however, if $\eta_n|_{Y_{Q,n}^{sc}}$ can be extended to a homomorphism $Y_{Q,n} \longrightarrow F^{\times}/F^{\times n}$. 
 \vskip 5pt

Another obstruction is the following. Suppose that $y \in nY_{Q,n} \cap Y_{Q,n}^{sc}$.  Then properties (a') and (d') require that 
\[ \chi (y(a)) =1  = \chi(y(a)) \cdot (\eta(y), a)_n \quad \text{for any $a \in F^{\times}$}. \]
 Thus we have our second obstruction:
 \vskip 5pt
 
 \noindent{\bf Obstruction 2:} 
 A genuine character $\chi$ satisfying (a') and (d') exists  if and only if (\ref{E:Ob1}) holds and
 \begin{equation}  \label{E:Ob2}
   (a, \eta(y))_n  =1 \quad \text{for any $y \in nY_{Q,n} \cap Y_{Q,n}^{sc}$ and any $a \in F^{\times}$.} \end{equation}
 \vskip 5pt
 
 \noindent  Again, this condition does not hold automatically, but it does hold if $\eta_n|_{Y_{Q,n}^{sc}}$ is extendable to $Y_{Q,n}$. 
 \vskip 5pt

 \vskip 5pt

\subsection{\bf Existence of splitting.}
 To summarise, we have shown:
 \vskip 5pt
 
 \begin{prop}
(i) To give a splitting of $E_1 + E_2$ is equivalent to giving a genuine character of $\overline{T}_{Q,n}$
satisfying conditions (a') and (d') above. 
 \vskip 5pt
 
 (ii) If $\eta_n|_{Y_{Q,n}^{sc}}$ is extendable to a homomorphism $Y_{Q,n}  \rightarrow F^{\times}/ F^{\times n}$, then  such genuine characters as in (i)  exist, so that the sequence $E_1 + E_2$ is split. 
\vskip 5pt

(iii)  Under the hypothesis in (ii), the set of such splittings is a torsor under the group 
\[  \Hom(F^{\times}/F^{\times n},  Z(\overline{G}^{\vee})) = \Hom(W_F,   Z(\overline{G}^{\vee})[n]). \]
\end{prop}

\vskip 5pt
\begin{proof}
We have already shown (i) . For (ii), the existence of such genuine characters follows by Pontrjagin duality, noting that $\overline{T}_{Q,n}$ is abelian. Finally (iii) is clear. 
\end{proof}

\vskip 5pt
We shall see later that these obstructions do occur in our running example. 
 
\vskip 10pt

\subsection{\bf Distinguished splitting,}

The definition of $E_1 + E_2$ uses essentially only the data $(n, Y_{Q,n}, \Phi_{Q,n}^{\vee}, Q, \mathcal{E}_{Q,n})$. In particular, it remembers nothing much of the original $(Y, Q)$.  
Recall that there is a short exact sequence
\[  \begin{CD}
{\rm Ker}(f)  @>>> \overline{T}_{Q,n}  @>>>  Z(\overline{T}).  \end{CD} \]
Since we are after all defining the L-group of $\overline{G}$, it is natural to consider those splittings of $E_1 + E_2$ which corresponds to those genuine characters of $\overline{T}_{Q,n}$ which are trivial on ${\rm Ker}(f)$ and thus factors to give characters of $Z(\overline{T})$.
For this, it follows by Lemma \ref{L:TQn} that it suffices to consider those genuine characters which are trivial on the image of $\tilde{g}:  T  \rightarrow \overline{T}_{Q,n}$, with
\[  \tilde{g}:  y(a)  \mapsto      (ny)(a)   \in \overline{T}_{Q,n}. \]
Then we have:
\vskip 5pt

\begin{itemize}
\item[(a'')] a genuine character $\chi$ of $\overline{T}_{Q,n}$ factors to $Z(\overline{T})$ if $\chi \circ \tilde{g} =1$. 
\end{itemize}
\noindent  Note that $\tilde{g}$ agrees with  $i_n$ when pulled back to $T_{Q,n}$, so that this requirement  subsumes condition (a').
It is now natural to make the following definition:
\vskip 10pt

\noindent{\bf Definition:} 
We call a genuine character of $\overline{T}_{Q,n}$ which satisfies (a'') and (d')  a {\em distinguished genuine character}. We call the corresponding splitting of $E_1 + E_2$ a {\em distinguished splitting}. 

\vskip 10pt

Assume that we have a genuine character satisfying (a') and (d') already. As before, we see that there is an obstruction to (a''). Namely, 
\vskip 5pt

\noindent{\bf Obstruction 3:}  
A genuine character satisfying (a'') and (d') exists if and only if (\ref{E:Ob1}) holds and
\[  (a, \eta(y))_n =1 \quad \text{for any $y \in nY \cap Y_{Q,n}^{sc}$ and any $a \in F^{\times}$.} \]
\vskip 5pt

This condition is not automatic, but is satisfied when $\eta_n : Y^{sc} \rightarrow F^{\times}/F^{\times n}$ can be extended to a homomorphism of $Y$. 

\vskip 5pt

We have thus shown:
\vskip 5pt

\begin{prop}  \label{P:d-split}
Assume that  $\eta_n : Y^{sc} \rightarrow F^{\times}/F^{\times n}$ can be extended to a homomorphism of $Y$, and let $J = nY +  Y_{Q,n}^{sc}$.
\vskip 5pt 
 
 (i) The set of distinguished genuine characters of $\overline{T}_{Q,n}$ is nonempty. 
\vskip 5pt

(ii)  Consider the subgroup 
\[  Z^{\heartsuit}(\overline{G}^{\vee}) :=  \Hom(Y_{Q,n}/J, \C^{\times}) \subset  Z(\overline{G}^{\vee}).\]
 Then the set of distinguished splittings of $E_1 + E_2$ is a torsor under 
\[  \Hom(F^{\times}/ F^{\times n},  Z^{\heartsuit}(\overline{G}^{\vee})) = \Hom(W_F, Z^{\heartsuit}(\overline{G}^{\vee})[n]). \]
\vskip 5pt

(iii)  Each distinguished splitting of $E_1 + E_2$ gives rise to a splitting 
$s:  W_F \longrightarrow {^L}\overline{G}$ which induces an isomorphism 
\[  \overline{G}^{\vee} \times W_F \cong  {^L}\overline{G}. \]
\end{prop}
\vskip 5pt

Henceforth, when $\eta_n : Y^{sc} \rightarrow F^{\times}/F^{\times n}$ can be extended to a homomorphism of $Y$,    we shall consider the L-group of a degree $n$ BD covering $\overline{G}$ as the extension
\[  \begin{CD}
1 @>>>  \overline{G}^{\vee}  @>>> {^L}\overline{G} @>>> W_F @>>> 1 \end{CD} \]
equipped with the set  of  distinguished splittings.   In particular, when $\G = \T$ is a split torus, then the set of distinguished splittings is nonempty since $Y^{sc}  = 0$. 

 \vskip 10pt
\subsection{\bf Weyl invariance.}
Under the hypothesis of Proposition \ref{P:d-split}, we have distinguished genuine characters of 
$\overline{T}_{Q,n}$. These distinguished characters possess another desirable property.
Namely, the action of $N(\T)$ on $\overline{\T}$ gives rise to an action of $N(\T)$ on $\overline{T}$ which preserves the center $Z(\overline{T})$. Moreover, the action on $Z(\overline{T})$ factors through
the Weyl group $W = N(\T)/ \T$.  Hence it makes sense to ask if a genuine character of $Z(\overline{T})$ is $W$-invariant. 
\vskip 5pt

For $\alpha \in \Phi$, let $w_{\alpha}$ denote the element of $W$ corresponding to the element $q( n_{\alpha}(1))$. Then we have seen in (\ref{E:nice-weyl}) that for $\tilde{t} \in Z(\overline{T})$, 
\[  {\rm Int}(w_{\alpha})( \tilde{t})   = \tilde{t}  \cdot   q(s_{\alpha}(\alpha^{\vee}(\alpha(t)))). \]  
Now if $t = y(a) \in T_{Q,n}$, with $y \in Y_{Q,n}$, then 
\[ q(s_{\alpha}(\alpha^{\vee}(\alpha(t)))) = s_{\eta}(\alpha^{\vee}(a)^{\langle \alpha, y \rangle}). \]
\vskip 5pt

\begin{prop}
For $y \in Y_{Q,n}$, $n_{\alpha}$ divides $\langle \alpha, y \rangle$ for each $\alpha \in \Phi$. Hence,
$\chi$ is $W$-invariant if $\chi$ satisfies (d').  In particular, distinguished genuine characters of $T_{Q,n}$ gives rise to $W$-invariant genuine characters of $Z(\overline{T})$.
\end{prop} 
\vskip 5pt

\begin{proof}
Since $Q$ is $W$-invariant, we have
\[  Q(y)  = Q(s_{\alpha} y)  = Q(y - \langle y, \alpha \rangle \alpha^{\vee}) = Q(y) - \langle y, \alpha \rangle \cdot B(y, \alpha^{\vee}) +  \langle y, \alpha \rangle^2 Q(\alpha^{\vee}). \]
This implies that
\[  \langle \alpha, y \rangle  = 0 \quad \text{or} \quad B(y, \alpha^{\vee}) = \langle \alpha, y \rangle \cdot Q(\alpha^{\vee}). \]
In the latter case, note that $B(y, \alpha^{\vee})$ is divisible by $n$ if $y \in Y_{Q,n}$, in which case $n_{\alpha}$ divides $\langle \alpha, y \rangle$ as desired.

\end{proof}
\vskip 5pt

 \subsection{\bf Splitting of ${^L}\overline{G}$.}
 We have so far considered splittings of the fundamental extension $E: = E_1 + E_2$ of $F^{\times}/ F^{\times n}$ by $Z(\overline{G}^{\vee})$ with good properties. Since the L-group of $\overline{G}$ is obtained from this fundamental extension by a combination of pushout and pullback, one may consider splittings of the extensions 
derived from these operations. Of course, these extensions will possess more splittings than the fundamental extension from which they are derived. 
\vskip 5pt

For example, one may pull back the fundamental extension to obtain
\[   \begin{CD} 
Z(\overline{G}^{\vee}) @>>>   \tilde{E} @>>> W_F \end{CD} \]
and one may ask for distinguished splittings of this central extension.  This amounts to dropping condition (a') above, so that one is finding genuine characters of $\overline{T}_{Q,n}$ which satisfy only (d'). To ensure that this character factors through to $Z(\overline{T})$, it would not be reasonable to require the condition (a'')  (since (a') is not assumed); one would simply require the character to be trivial on ${\rm Ker}(f)$.   \vskip 5pt

\vskip 5pt

  Now recall from \S \ref{SS:levi} that  the L-group construction is functorial with respect to inclusion of Levi subgroups, so that one has an embedding ${^L}\overline{T}  \hookrightarrow {^L}\overline{G}$.  
  Since we know that the extension ${^L}\overline{T}$ is split, we deduce that ${^L}\overline{G}$ is split as well.  In particular, ${^L}\overline{G}$ is abstractly a semi-direct product.  
As we have seen in our running example,  it is not a direct product in general. 
\vskip 5pt

Further, suppose we fix a fair bisector $D$ and consider all BD-extensions with BD invariants $(D, \eta)$.  
All these BD covering groups $\overline{G}_{\eta}$ have isomorphic covering tori  
$\overline{T}  = T \times_D  \mu_n$.  When $\eta_n =1$ is trivial,  we have seen that the set of distinguished splittings is nonempty. If we fixed a distinguished splitting of $E$ for $\eta_n=1$, we would have fixed a splitting of ${^L}\overline{T}$ and hence for ${^L}\overline{G}_{\eta}$ for all $\eta$.  If $\chi$ is the associated genuine character of $Z(\overline{T})$, note that $\chi$ is Weyl-invariant for the Weyl action associated to $\eta_n  =1$. 
But this $\chi$ need not be Weyl-invariant for general $\eta$, since the Weyl action on $Z(\overline{T})$ depends on $\eta$.

\vskip 10pt

\subsection{\bf Running example.}
We consider our running example to illustrate the discussion of this section. Recall that we have:
\[  \begin{CD}
\mu_2  @>>>  \overline{G}_{\eta}  =( \GL_2(F) \times \mu_2)/ \{  (z,  (\eta, z):  z \in F^{\times} \} @>>>  \PGL_2(F)  \end{CD} \]
and
\[  E_1^{\eta}  + E_2^{\eta} = \{ (t,a) \in \C^{\times} \times F^{\times}/ F^{\times 2}: t^2  =  (\eta, a) \}.\]
Moreover, $Y_{Q,2}  =  Y \supset Y_{Q,n}^{sc}  = Y^{sc} = 2 Y$.   In this case, Obstruction 1 says that 
\[  (\eta , -1) =1.  \]
Clearly, this may fail if $-1 \notin F^{\times 2}$. However, it does hold
  if $\eta \in F^{\times 2}$, or equivalently, if $\eta_2$ is trivial.  
 Obstruction 2 says that
 \[  (\eta, a)  =1 \quad \text{for all $a \in F^{\times}$.} \]
This is clearly a stronger condition than the one above, and it holds if and only if $\eta_2$ is trivial.  
 Thus, we see that if $\eta \notin F^{\times 2}$, then the sequence $E_1^{\eta}  + E_2^{\eta}$ does not split.  
 \vskip 5pt
 
 When $\eta_2$ is trivial, however, the above two obstructions, as well as Obstruction 3, are all absent
 and so a distinguished character of $\overline{T}_{\eta}  = Z(\overline{T}_{\eta})$ exists. 
 If we identify $\overline{T}$ with $A \times \mu_2$, so that genuine characters of $\overline{T}_{\eta}$ is in natural bijection with the characters of $A$, 
  then the distinguished characters correspond naturally to quadratic characters of $A$.  
 
 \vskip 5pt
 
 The Weyl group action on $\overline{T}_{\eta}$ is given by
 \[   \left( \begin{array}{cc}
 a & \\
  & 1  \end{array}  \right) \mapsto  \left( \left( \begin{array}{cc}
 a^{-1} & \\
  & 1  \end{array}  \right),  (\eta, a) \right), \]
 so that the Weyl action on genuine characters is given by
 \[  \chi  \mapsto \chi^{-1}  \cdot (\eta, -). \]
 In particular, we see that $\chi$ is $W$-invariant if and only if 
 \[   \chi^2  = (\eta, -), \]
 i.e. if $\chi$ is a square root of the quadratic character $(\eta,-)$. Such a square root exists if and only if 
 \[  (\eta, -1)  =1, \]
 i.e. if and only if Obstruction 1 is absent.
 \vskip 15pt

  \section{\bf Construction of Distinguished Genuine Characters}  \label{S:construction}
  It is useful in practice to have an explicit construction of the distinguished genuine characters (when they exist). Note that once one constructs one such distinguished character, the others are obtained by twisting it by a character of the finite group $T_{Q,n}/ T_J$. 
We shall see that a distinguished genuine character can be constructed using the Hilbert symbol $(-,-)_n$ and the Weil index $\chi_{\psi}$ associated to a nontrivial additive character $\psi$ of $F$. Hence, 
\[  \chi_{\psi}  :  F^{\times}  \longrightarrow \mu_4 \subset \C^{\times} \]
satisfies
\[  \chi_{\psi}(a)  \cdot \chi_{\psi}(b)  = \chi_{\psi}(ab)  \cdot (a,b)_2  \]
for all $a,b \in F^{\times}$. Moreover, $\chi_{\psi_{a^2}}  = \chi_{\psi}$. Such Weil indices play an important role in the classical theory of the metaplectic groups $\Mp_{2n}(F)$. 
\vskip 10pt

\subsection{\bf Reduction to rank 1 case.} 
By the elementary divisor theorem, we may pick a basis $\{ y_i \}$ of $Y_{Q,n}$ such that $\{ k_i y_i \}$ is a basis for the lattice $J = nY + Y_{Q,n}^{sc}$ for some $k_i \in \Z$. 
This gives a decomposition
\[  J=  k_1Y_1\oplus ...\oplus k_rY_r  \subset  Y_1 \oplus ...\oplus Y_r  = Y \]
and a map
\[  T_J  =  \prod_i (k_iY_i) \otimes_{\Z} F^{\times}  \longrightarrow \prod_i  Y_i \otimes_{\Z} F^{\times} \]
which is the $k_i$-power map on the $i$-th coordinate. 
Write $T_{Q,n,i}$ for the 1-dimensional torus corresponding to $Y_i$ and $T_{J,i}$ for that corresponding to $k_iY_i$. Now, because $\overline{T}_{Q,n}$ is abelian,  we see that
\[  \overline{T}_{Q,n}  \cong   \overline{T}_1 \times .....\times \overline{T}_r /  Z \]
where
\[  Z = \{  (\epsilon_i) \in \prod_{i=1}^r  \mu_n:  \prod_i \epsilon_i =1 \}. \]
Moreover, the group law on $\overline{T}_{Q,n,i}  = (Y_i \otimes F^{\times})  \times \mu_n$ is given by
\[  y_i(a) \cdot   y_i(b)  = y_i(ab) \cdot (a,b)_n^{Q(y_i)}. \]
It follows that the map $T_{J,i} \rightarrow T_{Q,n,i}$ splits naturally to give 
\[  T_{J,i} \longrightarrow \overline{T}_{Q,n,i}.\] 
\vskip 5pt

Thus, to construct a distinguished genuine character on $\overline{T}_{Q,n}$, it suffices to construct a genuine character of $\overline{T}_{Q,n,i}$ which is trivial on the image of $T_{J,i}$; the product of these characters will then be a distinguished genuine character of $\overline{T}_{Q,n}$. 
We are thus reduced to constructing genuine characters of 1-dimensional covering tori. 
\vskip 10pt

\subsection{\bf The definition.}
We set
\[  \chi((y_i(a)) =  \chi_{\psi}(a)^{f_i}  \]
for some $f_i \in \Z$ to be determined.  
Now we need to check various requirements:
\vskip 5pt

\begin{itemize}
\item we first need to check the relation for $\overline{T}_{Q,n,i}$:
\[  \chi( y_i(a))  \cdot   \chi(y_i(b))  =  \chi(y_i(ab)) \cdot (a,b)_n^{Q(y_i)}. \]
This amounts to the requirement that
\[  f_i  \equiv  A_i:= \frac{2}{n} \cdot Q(y_i)  \mod 2. \]
\vskip 5pt

\item next we need to ensure that $\chi$ is trivial on the image of $T_{J,i}$, i.e.
trivial on $y_i(a^{k_i})$. 
But a short computation gives:
\[  \chi(y_i(a^{k_i})) = \chi_{\psi}(a)^{k_i f_i +  k_i(k_i-1) \cdot A_i}. \]
Thus we need
\[  k_i \cdot (f_i + (k_i-1)A_i)  \equiv 0 \mod 4. \]
To ensure this, we shall simply take
\[  f_i := -(k_i-1) A_i. \]
Then when $k_i$ is even, it is automatic that  $f_i  \equiv A_i \mod 2$.
We need to ensure that this continues to hold when $k_i$ is odd.
\vskip 5pt

For this, we need to show that $A_i = 0 \mod 2$  when $k_i$ is odd; equivalently, we need to show that $Q(y_i) \equiv 0 \mod n$. Since we already know that 
$Q(y_i)$ is divisible by $n/2$, it remains to show that if $2^e$ exactly divides $n$, 
then $2^e$ divides $Q(y_i)$.  As  $k_iy_i \in J = Y_{Q,n}^{sc} + nY$, we know that $Q(k_i y_i)  \equiv 0 \mod n$. So $2^e$ divides $k_i^2 \cdot Q(y_i)$. Since $k_i$ is odd, we see that $2^e$ does divide $Q(y_i)$, as desired. 
\end{itemize}
 \vskip 5pt
 
 We have completed the construction of a distinguished genuine character $\chi_0$. A formula for 
$\chi_0$ can be given as follows.  For $y = \sum_i n_i y_i \in Y_{Q,n}$ and $a \in F^{\times}$,
\[  \chi_0( y(a))  =  \prod_i  \chi_{\psi}(a^{n_i})^{f_i}  \cdot (a,a)_n^{\sum_{i<j}  n_in_j D(y_i, y_j)}, \]
with
\[  f_i = -  (k_i-1) \cdot \frac{2}{n}  \cdot Q(y_i). \]
Though explicit, a slightly unsatisfactory aspect of this formula is that one first needs to find  compatible bases
for the lattice chain $J \subset Y_{Q,n}$.  
 
\vskip 15pt

 \section{\bf LLC for Covering Tori}
 We are going to specialise the investigation of \S \ref{S:d-splitting} to several examples. In this section, we assume that $\G  = \T$ is a spit torus.

\vskip 5pt

\subsection{\bf LLC for $\overline{T}_{Q,n}$}
 When $\G   = \T$, one has $Y^{sc} = 0$, so that $\eta  =1$,  and the extension $E_1 + E_2$ is
 \[  \begin{CD}
 \overline{T}_{Q,n}^{\vee}   @>>>  E_1+E_2  @>>>  F^{\times}/ F^{\times n}. \end{CD} \]
 In the previous section, we have seen that to give a splitting of this sequence is equivalent to giving a genuine character on $\overline{T}_{Q,n}$ satisfying conditions (a') and (d'). Since $Y^{sc} = 0$, (d') is vacuous, and since $\eta =1$, (a') holds. So we obtain a bijection between the set of splittings of $E_1 + E_2$ and the set of genuine characters of $\overline{T}_{Q,n}$ satisfying (a'). 
 \vskip 5pt
 
 On the other hand, to obtain the L-group of $\overline{T}$, we pull $E_1+ E_2$ back by $W_F \longrightarrow F^{\times}  \longrightarrow F^{\times}/ F^{\times n}$. The same considerations show that to give a splitting of ${^L}\overline{T}$  is equivalent to giving a genuine character of $\overline{T}_{Q,n}$, where we don't insist on condition (a') anymore.  Hence, we have obtained  a natural bijection 
 \[  \{  \text{Splittings of ${^L}\overline{T}$} \} 
 \longleftrightarrow  \{  \text{genuine characters of $\overline{T}_{Q,n}$} \}. \]
  This is a classification of the genuine characters  of the abelian group $\overline{T}_{Q,n}$.
   \vskip 5pt
 
\subsection{\bf Construction for $\overline{T}$}
 As shown in \cite{W1}, $\overline{T}$ is a Heisenberg type group, and one has a natural bijection
 \[  \{  \text{genuine characters of $Z(\overline{T})$} \}  \longleftrightarrow 
 \{  \text{irreducible genuine representations of $\overline{T}$} \}. \]
 This bijection is defined as follows.  Choose and fix a maximal abelian subgroup $H$ of $\overline{T}$ containing $Z(\overline{T})$.
 Given a genuine character $\chi$ of $Z(\overline{T})$,  extend $\chi$ arbitrarily to a character $\chi_H$ of $H$. Then the induced representation
 \[  i(\chi)  = {\rm ind}_H^{\overline{T}}  \chi_H \]
 is irreducible and independent of the choice of $(H, \chi_H)$. By the analog of the Stone-von-Neumann theorem, it is characterized as the unique irreducible genuine representation of $\overline{T}$ which has central character $\chi$.  
 \vskip 5pt

 \subsection{\bf LLC for $\overline{T}$}
  Combining these, we obtain the following result which is the LLC for covering (split) tori:
 \vskip 5pt
 
 \begin{thm}
 There is a natural injective map
 \[  \mathcal{L}_{\overline{T}}: 
  \{  \text{irreducible genuine representations of $\overline{T}$} \} \hookrightarrow 
  \{  \text{Splittings of ${^L}\overline{T}$} \}.   \]
 \end{thm}
 \vskip 5pt
 
 The image of this injection can be described as follows.  
 If we fix a distinguished splitting $s_0$, which corresponds to a genuine character $\chi_0$ of $Z(\overline{T})$, then all other splittings of ${^L}\overline{T}$ are of the form $s  = s_0 \cdot \rho$ 
 where $\rho:  W_F  \longrightarrow \overline{T}^{\vee}  = T_{Q,n}^{\vee}=  X_{Q,n}  \otimes_{\Z}  \C^{\times}$.  
 Now the isogeny $f: \T_{Q,n}  \longrightarrow \T$ induces an isogeny  
 \[  f^*:  T^{\vee}  =  X  \otimes \C^{\times}  \longrightarrow T_{Q,n}^{\vee}  = X_{Q,n} \otimes_{\Z} \C^{\times},\]
 and thus a map
 \[ f_*:  \Hom(W_F,  T^{\vee})  \longrightarrow  \Hom(W_F,  T_{Q,n}^{\vee}). \]
 A splitting $s$ of ${^L} \overline{T}$ is in the image of the map $\mathcal{L}_{\overline{T}}$ if and only if its corresponding $\rho \in \Hom(W_F,  T_{Q,n}^{\vee})$ lies in the image of $f_*$. Could this be made more explicit?
 \vskip 10pt
 
 \section{\bf  LLC for Unramified Representations} \label{S:LLC-unram}
 
 We consider the tame case in this section, so that $p$ does not divide $n$. We assume that there is a splitting
 \[  s:  K  = \G(\mathcal{O})  \longrightarrow \overline{G}. \]
 Then one may consider the  $s$-unramified genuine representations of $\overline{G}$, i.e, those with nonzero $s(K)$-fixed vectors. We would like to obtain an LLC for such $s$-unramified genuine representations.
  \vskip 5pt
 
 \subsection{\bf Torus case.}
 We first consider the case when $\G = \T$ is a a split  torus.  Suppose that $\overline{\T}$ has bisector data $D$, so that
 \[  \overline{T}  =  T \times_D   \mu_n  \]
 In the tame case, the trivial section $y(a) \mapsto (y(a), 1) \in T \times_D \mu_n$ is a splitting over $\T(\mathcal{O})$ and any  splitting $s:  \T(\mathcal{O}) \longrightarrow \overline{T}$ is given by
 \[  s_{\mu}(t)  = (t  , \mu(t))  \quad \text{for $t \in \T(\mathcal{O}^{\times})$}  \]
 where $\mu:   \T(\mathcal{O}) \longrightarrow \mu_n$ is a group homomorphism. Such an $s:= s_{\mu}$ will give rise to a splitting of $\T_{Q,n}(\mathcal{O})$, denoted by $s$ as well, via pulling back. 
 \vskip 5pt
 
 We call a genuine representation $i(\chi)$  of $\overline{T}$ $s_{\mu}$-unramified if $i(\chi)$ has a nonzero vector fixed by $s_{\mu}(T(\mathcal{O}))$. In this tame case, one can check that $H = Z(\overline{T}) \cdot s_{\mu}(T(\mathcal{O}))$ is a maximal abelian subgroup of $\overline{T}$. From this, one deduces:
 \vskip 5pt
 
 \begin{lemma}
 A representation $i(\chi)$ is $s$-unramified if and only if $\chi$ is trivial when restricted to $Z(\overline{T}) \cap s(T(\mathcal{O}))$.  In this case, the space of $ s(T(\mathcal{O}))$-fixed vectors is $1$-dimensional. 
  \end{lemma}
 \vskip 5pt
 
  We say that a genuine character of $Z(\overline{T})$ or $\overline{T}_{Q,n}$ is $s$-unramfied if it is trivial on the image of $s$.   Thus, the above lemma says that $i(\chi)$ is $s$-unramified if and only if $\chi$ is $s$-unramified. Such a $\chi$ will pullback to an $s$-unramified character of $\overline{T}_{Q,n}$. 
  Conversely, observe that an $s$-unramified genuine character of $\overline{T}_{Q,n}$ automatically factors through to a genuine  character of $Z(\overline{T})$ since ${\rm Ker}(f)  \subset s(T_{Q,n}(\mathcal{O}))$. 
  Under the LLC for $\T$ defined in the last section, the $s$-unramified genuine characters  corresponds to a subset of splittings of ${^L}\overline{T}$. We would like to explicate this subset. 
 \vskip 5pt
 
 Examining the construction of the LLC for covering tori, we see that
 all the $s$-unramified characters $\chi$ have the property that their L-parameters $\rho_{\chi} :  W_F \longrightarrow {^L}\overline{T}$ have the same restriction to the inertia group $I_F$.  In other words, as a map from $W_F^{ab}  = F^{\times}$ to ${^L}\overline{T}$, all these $\rho_{\chi}$'s have the same restriction $\rho_s$ on $\mathcal{O}^{\times}$.     This restriction $\rho_s$ can be described   explicitly as follows.
 Suppose an $s$-unramified $\chi$ is given. 
 For $a \in \mathcal{O}^{\times}$ , $\rho_{\chi}(a) \in {^L}\overline{T}$ is basically given by the 
 function $\tilde{\chi}_a : Y_{Q,n} \longrightarrow \C^{\times}$.  The condition (b) in \S \ref{SS:splitting}  says that $\tilde{\chi}_a \in \Hom(Y_{Q,n}, \C^{\times})$, and  the $s$-unramified condition says that 
 \[    \tilde{\chi}_a(y)  = \mu(y(a))^{-1}. \]
 We shall call a splitting $\rho$ of ${^L}\overline{T}$ {\em $s$-unramified} if $\rho|_{I_F} = \rho_s$.
  
 \vskip 5pt
 
To summarize, we have:
\begin{prop}
 Under the LLC for covering tori, one has a bijection
 \[  \{  \text{irreducible $s$-unramified genuine representations of $\overline{T}$}\}  \]
 \[  \updownarrow \]
 \[   \{  \text{$s$-unramified splittings $\rho$ of ${^L}\overline{T}$ } \}. \]
 \end{prop}

 \subsection{\bf Satake isomorphism.}
 We now consider the case of general $\G$. The key step in understanding the $s$-unramified representation of $\overline{G}$ is the Satake isomorphism. More precisely, 
 let $\mathcal{H}(\overline{G})$ be the $\C$-algebra of anti-genuine locally constant, compactly supported functions on $\overline{G}$ which are bi-invariant under $s(T(\mathcal{O})$.  Let $\mathcal{H}(\overline{T})$ denote the analogous $\C$-algebra for $\overline{T}$. One can check that for an element $f \in \mathcal{H}(\overline{T})$, the support of $f$ is contained in $Z(\overline{T}) \cdot s(T(\mathcal{O}))$, and thus $f$ is completely determined by its restriction to 
 $Z(\overline{T})$.  Moreover, the Weyl group $W$ acts naturally on $Z(\overline{T}) \cdot s(T(\mathcal{O}))$ and thus on $\mathcal{H}(\overline{T})$.  
 
One has an explicit $\C$-algebra morphism
$$ \mathcal{S}  :   \mathcal{H}(\overline{G}) \longrightarrow  \mathcal{H}(\overline{T})$$
given by
$$\mathcal{S}(f)(t)=\delta(t)^{1/2} \int_U f(tu)du \text{ for all } f\in \mathcal{H}(\overline{G}).$$

It is shown in \cite{Mc2, L2, W7} that $\mathcal{S}$ gives an isomorphism of $\C$-algebras:
 \[  \mathcal{S}  :  
 \mathcal{H}(\overline{G}) \longrightarrow  \mathcal{H}(\overline{T})^W. \] 

 As a consequence of this, one deduces a bijection
   \[  \{  \text{irreducible $s$-unramified genuine representations of $\overline{G}$}\} \]
   \[ \updownarrow \]
   \[ \{  \text{irreducible modules of $\mathcal{H}(\overline{G})$} \} \]
 \[ \updownarrow \]
\[  \{  \text{$W$-orbits of $s$-unramified genuine characters of $Z(\overline{T})$} \}. \]
On the other hand, by pulling back, one has a bijection of the latter set with
\[  \{  \text{$W$-orbits of $s$-unramified genuine characters of $\overline{T}_{Q,n}$} \} \]
Applying the LLC for $\overline{T}_{Q,n}$, one then obtains:
  \vskip 5pt
 
 \begin{thm}
 The Satake isomorphism gives natural bijections
 \[  \{  \text{irreducible $s$-unramified genuine representations of $\overline{G}$}\} \]
 \[  \updownarrow \]
\[  \{  \text{$W$-orbits of $s$-unramified splittings of ${^L}\overline{T}$} \} \]
\[  \updownarrow \]
\[  \{  \text{$\overline{G}^{\vee}$-orbits of $s$-unramified splittings of ${^L}\overline{G}$} \}.   \]
 \end{thm}

\vskip 10pt

\section{\bf L-Groups: Second Take}

After the discussion of the previous sections and a study of our running example, we may draw the following conclusions:
\vskip 5pt
\begin{itemize}
\item[(i)] For a fixed {\em fair}  bisector $D$, and among all BD covering groups (of degree $n$) with bisector data $(D,\eta)$, those with $\eta_n  =1$ are most nicely behaved. For example, their maximal covering tori $\overline{T}$ have certain  distinguished Weyl-invariant genuine representations and the hyperspecial maximal compact subgroup $G(\mathcal{O})$ splits in $\overline{G}$.  Moreover, their L-groups are isomorphic to a direct product $\overline{G}^{\vee} \times W_F$. 
\vskip 5pt

\item[(ii)]  The BD covering groups $\overline{G}_{\eta}$ for a fixed bisector data are closely related, and it may be useful to consider them together, both structurally as well as from the point of view of representation theory. For example, they all have the same dual group $\overline{G}_Q^{\vee}$.
\end{itemize}

In this section, we would like to suggest a slightly different take on
 the L-group extension, so as to treat  the closely related groups $\overline{G}_{\eta}$ together.  
\vskip 5pt

 \vskip 5pt

\subsection{\bf The case $Q=0$.}
To guide our efforts, we shall consider the genuine representation theory of  the covering groups in the case when $Q=0$. This expands upon our running example and will provide a clue about the modifications needed.
\vskip 5pt

When $Q = 0  =D$,  the objects in ${\bf Bis}_{\G,Q}$ are simply homomorphisms $\eta: Y^{sc} \rightarrow F^{\times}$. Choose a $z$-extension 
\[  \begin{CD} 
Z @>>> \mathbb{H} @>>> \G, \end{CD} \]
so that $Y_{\mathbb{H}}^{sc}  = Y^{sc}$ and $Y_{\mathbb{H}}/ Y^{sc}$ is free. 
For any $\eta$, we have a corresponding short exact sequence 
\[  \begin{CD}
Z = Z_{\eta} @>>> \overline{\mathbb{H}}_{\eta} @>>> \overline{\G}_{\eta} \end{CD} \]
Since all $\eta :  Y_{\mathbb{H}}^{sc} = Y^{sc} \rightarrow F^{\times}$ are equivalent to the trivial homomorphism $1$ as objects of ${\bf Bis}_{\mathbb{H}, Q}$, we may choose an isomorphism 
\[  \xi : \overline{\mathbb{H}}_1 = \mathbb{H} \times K_2  \longrightarrow \overline{\mathbb{H}}_{\eta}. \]
After taking $F$-points, and noting that $H^1_{Zar}(F, Z)  =0$, we then have
\[ \begin{CD}
  \overline{H}_1  = H  \times \mu_n @>\xi>> \overline{H}_{\eta}  @>>> \overline{G}_{\eta}, \end{CD} \]
  and the kernel of this map is the subgroup
  \[  \xi^{-1} (Z_{\eta}) = \{ (z,  \chi_{\eta, \xi}(z)^{-1} ): z\in  Z\}  \subset H \times \mu_n, \]
  where $\chi_{\eta, \xi}$ is the map
  \[ \begin{CD}
   Z = Y_Z \otimes F^{\times} @>\xi>> F^{\times} \otimes F^{\times} @>(-,-)_n>>  \mu_n. \end{CD} \]
   Hence, the set of  genuine representations of $\overline{G}_{\eta}$ can be identified (by pulling back) with a subset of the  genuine representations of the split extension $H \times \mu_n$ whose restriction to the central subgroup $Z \subset H$ is the character $\chi_{\eta,\xi}$. 
 
 \vskip 5pt
 
 Now the L-group of $\overline{H} = H \times \mu_n$  is a short exact sequence 
 \[  \begin{CD}  
 H^{\vee}  @>>> {^L}\overline{H} @>>>  W_F \end{CD} \]
 which is equipped with a finite set of distinguished splittings. For example one may take  the distinguished splitting $s_0$ which corresponds to the trivial character of the maximal torus $T_H$ of $H$. Then we may identify the set of all splittings (modulo conjugacy by $H^{\vee}$) with the set of L-parameters
 \[  W_F \longrightarrow H^{\vee}  \]
of $H$. Thus, if the LLC holds for the linear group $H$,  there is a finite-to-one map
\[  \bigcup_{\eta}  {\rm Irr}_{gen}(\overline{G}_{\eta})  \longrightarrow \{ \text{L-parameters for $H$} \} \]
which may be construed as a (weak) LLC for the family of covering groups $\overline{G}_{\eta}$ (as $\eta$ varies). 
Moreover, the image of  ${\rm Irr}_{gen}(\overline{G}_{\eta})$ for a particular $\eta$
can be described as follows.  By Lemma \ref{L:z-ext}, there is a natural short exact sequence:
\[  \begin{CD}
G^{\vee}  @>>>  H^{\vee}  @>\rho>> Z^{\vee} = \Hom(Y_Z, \C^{\times}). \end{CD}  \]
 Then  ${\rm Irr}_{gen}(\overline{G}_{\eta})$ is the set of L-parameters $\phi$ of $H$ such that $\rho \circ \phi : W_F \longrightarrow Z^{\vee}$ is the L-parameter of the character $\chi_{\eta,\xi}$. Observe that the genuine representations of $\overline{G}_1  = G \times \mu_n$ is then parametrized by L-parameters of $H$ which factors through $G^{\vee}$.
\vskip 5pt

  There is an obvious  notion of a $z$-extension $H$ dominating another $H'$, and one can easily check that the above classification of the genuine representations of $\overline{G}_{\eta}$ behave functorially with respect to dominance. This suggests that it is possible (and certainly desirable) to formulate the LLC for $\overline{G}_{\eta}$ without reference to the $z$-extension $H$.  
  \vskip 5pt

\subsection{\bf Modification of $L$-group.} Motivated by  the above discussion, we can revisit the L-group construction.   Fix the quadratic form $Q$ on $Y$.  The crucial $E_2$ construction starts with 
\[  \begin{CD}
F^{\times}/ F^{\times n}  @>>>  \mathcal{E}_{Q,n}  @>>>  Y_{Q,n} \end{CD} \]
and then use the section
\[  s_{\eta}:  Y_{Q,n}^{sc}  \longrightarrow  \mathcal{E}_{Q,n} \]
to form the quotient
\[  \begin{CD}
F^{\times}/ F^{\times n}  @>>>  \mathcal{E}_{Q,n}/ s_{\eta}(Y_{Q,n}^{sc})  @>>>  Y_{Q,n}/Y_{Q,n}^{sc}, \end{CD} \]
before applying $\Hom( -, \C^{\times})$.   
To incorporate all $\eta$'s together, we observe that  the section $s_{\eta}$ is independent of $\eta$ when restricted to the sublattice $nY^{sc}_{Q,n}$. 
Then  one has the commutative diagram with exact rows:  
\[ \begin{CD}
 F^{\times}/ F^{\times n}  @>>>  \mathcal{E}_{Q,n}/ s_{\eta}(nY_{Q,n}^{sc})  @>>>  Y_{Q,n}/nY_{Q,n}^{sc} \\
 @|      @VVV    @VVV  \\
 F^{\times}/ F^{\times n}  @>>>  \mathcal{E}_{Q,n}/ s_1(Y_{Q,n}^{sc})  @>>>  Y_{Q,n}/Y_{Q,n}^{sc}.
\end{CD} \]
Taking $\Hom(-, \C^{\times})$,  one obtains  the commutative diagram with exact rows and columns, which defines the modification $\tilde{E}_2$ of $E_2$:
\[ \begin{CD}
 Z(\overline{G}_1^{\vee})=   \Hom(Y_{Q,n}/Y_{Q,n}^{sc}, \C^{\times})   @>>>   E_2    @>>>  F^{\times}/ F^{\times n} \\
 @VVV      @VVV    @|  \\
  \Hom(Y_{Q,n}/nY_{Q,n}^{sc}, \C^{\times})  @>>> \tilde{E}_2   @>>>  F^{\times}/ F^{\times n} \\
  @VVV    @VVV  @VVV   \\
  (T_{Q,n}^{sc})^{\vee}[n] = \Hom(Y_{Q,n}^{sc}/nY_{Q,n}^{sc}, \C^{\times})  @=  (T_{Q,n}^{sc})^{\vee}[n] @>>> 1 
\end{CD} \]
The cocycle defining $E_1$ also defines an extension
\[  \begin{CD}
\Hom(Y_{Q,n}/nY_{Q,n}^{sc}, \C^{\times})   @>>> \tilde{E}_1 @>>> F^{\times}/F^{\times n}. \end{CD} \]
Then we can form the Baer sum and obtain
\[ \begin{CD}
 Z(\overline{G}_1^{\vee})    @>>>  E =  E_1 + E_2    @>>>  F^{\times}/ F^{\times n} \\
 @VVV      @VVV    @|  \\
 \Hom(Y_{Q,n}/ nY_{Q,n}^{sc}, \C^{\times})  @>>>\tilde{E} = \tilde{E}_1 +  \tilde{E}_2   @>>>  F^{\times}/ F^{\times n} \\
  @VVV    @VVV  @VVV   \\
   (T_{Q,n}^{sc})^{\vee}[n]  @=    (T_{Q,n}^{sc})^{\vee}[n]  @>>> 1 
\end{CD} \]
From this, we infer the short exact sequence:
\[  \begin{CD}
Z(\overline{G}_1^{\vee}) @>>>  \tilde{E}  @>>>  F^{\times}/F^{\times n}  \times (T_{Q,n}^{sc})^{\vee}[n]. 
\end{CD}  \]
This is our enlarged fundamental extension.  Pushing this out by $Z(\overline{G}_1^{\vee}) \hookrightarrow \overline{G}^{\vee}$ and pulling back to $W_F$,  one obtains:
\[  \begin{CD}
\overline{G}_Q^{\vee}  @>>>  {^L}\overline{G}_Q^{\#}  @>>>  W_F \times (T_{Q,n}^{sc})^{\vee}[n]  \end{CD} \]
which is our enlarged L-group extension for the family of BD covers with fixed BD-invariant $Q$.   
\vskip 5pt

Note that $(T_{Q,n}^{sc})^{\vee}$ is a maximal  torus in the adjoint quotient $(\overline{G}^{\vee})_{ad}$ of $\overline{G}^{\vee}$, so that
$ (T_{Q,n}^{sc})^{\vee}[n] $ is its $n$-torsion  subgroup.
\vskip 5pt

\subsection{\bf Relation with ${^L}\overline{G}_{\eta}$}
How can one recover the L-group of $\overline{G}_{\eta}$, as previously defined, from the enlarged L-group defined here?   Given an $\eta :  Y_{Q,n}^{sc}  \longrightarrow F^{\times}$, one obtains a natural map
\[  \varphi_{\eta}:  W_F  \longrightarrow W_F^{ab}  = F^{\times}  \longrightarrow (T_{Q,n}^{sc})^{\vee}[n]  = \Hom(Y_{Q,n}^{sc}, \mu_n) \]
given by
\[  \varphi_{\eta}(a) (y)  = (\eta(y),  a)_n  \quad \text{for all $a \in F^{\times}$ and $y \in Y_{Q,n}^{sc}$.} \]
Pulling back the enlarged L-group extension by the diagonal map 
\[  {\rm id}  \times \varphi_{\eta}:  W_F \longrightarrow W_F \times (T_{Q,n}^{sc})^{\vee}[n], \]
one obtains the L-group extension ${^L}\overline{G}_{\eta}$.   
 \vskip 10pt

 \subsection{\bf The modified dual group.}
 Consider the kernel  $\overline{G}_Q^\# \subseteq {^L}\overline{G}_Q^\#$   of the following composition of surjections:
$$\begin{CD}
 {^L}\overline{G}_Q^{\#}  @>>>  W_F \times (T_{Q,n}^{sc})^{\vee}[n]  @>>> W_F,\end{CD}$$
where the second map is the projection on the first component. By the definition of ${^L}\overline{G}_Q^\#$, the group $\overline{G}_Q^\#$ lies in the exact sequence
\begin{equation} \label{E:mod-dual}
\begin{CD}
\overline{G}_Q^\vee @>>>  \overline{G}_Q^{\#}  @>>>  (T_{Q,n}^{sc})^{\vee}[n] ,\end{CD}
\end{equation}
which is the push out of 
\begin{equation} \label{modi dual}
\begin{CD}
Z(\overline{G}_Q^\vee) @>>>  \text{Hom}(Y_{Q,n}/nY_{Q,n}^{sc}, \C^\times)  @>>>  (T_{Q,n}^{sc})^{\vee}[n]. 
\end{CD}
\end{equation}
Since the adjoint group $(\overline{G}^{\vee})_{ad}$ acts naturally on $\overline{G}^{\vee}$, preserving $Z(\overline{G}^{\vee})$, 
there is a canonical splitting of $\overline{G}_Q^\#$ over $(T_{Q,n}^{sc})^{\vee}[n]$, which gives a canonical isomorphism
$$\overline{G}_Q^\# \simeq \overline{G}_Q^\vee \rtimes (T_{Q,n}^{sc})^{\vee}[n].$$
\vskip 5pt

Moreover, the action of  $(T_{Q,n}^{sc})^{\vee}[n]$ on $\overline{G}^{\vee}$ is identity on its maximal torus $\overline{T}^{\vee} = X_{Q,n} \otimes \C^{\times}$ and preserves the maximal unipotent subgroup corresponding to the set of simple roots $\Delta_{Q,n}^{\vee}$.
 This shows that every irreducible representation is invariant under the action of $ (T_{Q,n}^{sc})^{\vee}[n]$, and thus extends (non-canonically) to $\overline{G}^{\#}$. In other words, the representation theory of the disconnected group $\overline{G}^{\#}$ is not more complicated than that of $\overline{G}^{\vee}$.


 \subsection{\bf Running example}
Consider the case $\G=\PGL_2$ and $D=0$, $n=2$. The dual group is $\overline{G}_Q^\vee=\SL_2(\C)$ and the exact sequence (\ref{modi dual}) is
$$\begin{CD}
\mu_2 @>>> \mu_4  @>>>  \mu_2. 
\end{CD}$$
We obtain $\overline{G}_Q^\#\simeq \SL_2(\C) \rtimes \mu_2$. The action of the nontrivial element $\epsilon \in \mu_2$ on $\SL_2(\C)$ is given by
$$\epsilon : g \mapsto 
\left( \begin{array}{cc}
 \xi & \\
  & \xi^{-1}  \end{array}  \right) 
g
\left( \begin{array}{cc}
 \xi & \\
  & \xi^{-1}  \end{array}  \right)^{-1} ,$$
where $\xi \in \mu_4$ is any square root of $\epsilon$. As there is no splitting of $\overline{G}_Q^\#$ over $\mu_2$ valued in the center $Z(\overline{G}_Q^\vee)$, the group $\overline{G}_Q^\#$ is not isomorphic to the direct product $\SL_2(\C) \times \mu_2$.

This example shows that in general $\overline{G}_Q^\# \simeq \overline{G}_Q^\vee \rtimes (T_{Q,n}^{sc})^\vee[n]$ is not a direct product of the two groups.

\vskip 10pt

\section{\bf The LLC}
 After the discussion in the previous sections, we can now formulate the LLC for BD covering groups. 
 
  \vskip 5pt
 
\subsection{\bf L-parameters.}
After introducing the L-group extension, one now has the following notions:
\vskip 5pt

\begin{itemize}
\item An L-parameter for the covering group $\overline{G}_{\eta}$ is a splitting $\phi: W_F \longrightarrow {^L} \overline{G}_{\eta}$ of the extension ${^L} \overline{G}_{\eta}$, taken up to conjugacy by $\overline{G}^{\vee}$.
Equivalently, it is a splitting $\phi:  W_F \longrightarrow {^L}\overline{G}_Q^{\#}$ of the enlarged L-group extension  such that 
\[  p \circ \phi   =  \varphi_{\eta}, \]
 where $p :  {^L}\overline{G}_Q^{\#}\longrightarrow (T_{Q,n}^{sc})^{\vee}[n]$ is the natural projection.
\vskip 5pt

\item we have demonstrated the existence of a finite set of distinguished splittings for ${^L}\overline{G}_1$ and thus for ${^L}\overline{G}_Q^{\#}$.  If we fix one such splitting $\phi_0$, then all  splittings of ${^L}\overline{G}_Q^{\#}$ are of the form $\phi_0 \cdot \phi$ where
\[  \phi:  W_F \longrightarrow   \overline{G}_Q^{\#}  = \overline{G}_Q^{\vee} \rtimes (T_{Q,n}^{sc})^{\vee}[n]. \]
We call such $\phi$'s the L-parameters relative to the distinguished splitting $\phi_0$. 
\end{itemize}

 \vskip 10pt

\subsection{\bf  Local $L$-factors}
Given a representation 
\[  R:  {^L}\overline{G}_{\eta}  \longrightarrow \GL(V)  \]
where $V$ is a complex finite-dimensional vector space over $\C$, and a splitting $\phi$ of ${^L}\overline{G}_{\eta}$, one obtains a complex representation $R \circ \phi$ of $W_F$ and hence an Artin L-factor   $L(s, \phi, R)$. Alternatively, if $\phi:  W_F  \longrightarrow \overline{G}^{\#}$ is an L-parameter relative to a distinguished splitting $\phi_0$ of ${^L}\overline{G}^{\#}$ over $W_F$, and $R:  \overline{G}^{\#} \longrightarrow \GL(V)$ is a representation, then one has an associated L-factor $L_{\phi_0}(s,\phi, R)$. As we noted before, the irreducible representations of $\overline{G}^{\#}$ are simply extensions of those of $\overline{G}^{\vee}$.  We shall give a more detailed treatment of this in the next section, where we introduce automorphic L-functions.   
\vskip 5pt

\subsection{\bf The LLC}
In view of the unramified LLC discussed in Section \ref{S:LLC-unram}, one is tempted to conjecture the existence of a finite-to-one map giving rise to a commutative diagram:
\[ \begin{CD}
 \mathcal{L}_{\eta}: {\rm Irr} \overline{G}_{\eta}  @>\mathcal{L}_{\eta}>>  \{  \text{splittings of ${^L}\overline{G}_{\eta}$}\} \\
 @VVV    @VVV  \\
 \bigcup_{\eta}  {\rm Irr} \overline{G}_{\eta}  @>\mathcal{L}>>\{  \text{splittings of ${^L}\overline{G}_Q^{\#}$} \}.  
 \end{CD} \]
 This is a weak local LLC. 
 As shown in the case when $\G  = \T$ is a split torus, one should not expect this map to be surjective. Thus, one would also like to have a conjectural parametrization of the fibers of this map.  This would be a strong LLC. 
\vskip 5pt

\subsection{\bf Reduction to $\eta = 1$.}
We shall show that this weak LLC for general $\overline{G}_{\eta}$ can be reduced to the case of trivial $\eta$.  This is similar to the discussion at the beginning of the section and relies on the consideration of $z$-extensions. 
\vskip 5pt

More precisely, we choose a $z$-extension $Z \longrightarrow \mathbb{H} \longrightarrow \G$ giving rise to
\[  \begin{CD}
Z @>>> \overline{\mathbb{H}}_{\eta}  @>p>> \overline{\G}_{\eta}.   \end{CD} \]
Choose an isomorphism $\xi:  \overline{\mathbb{H}}_1  \longrightarrow  \overline{\mathbb{H}}_{\eta}$, which realises
\[   \overline{\G}_{\eta}  \cong  \overline{\mathbb{H}}_1  /  \xi^{-1} (Z). \]
Thus one has an  injection
\[ \xi^* \circ p^* :  {\rm Irr} \overline{G}_{\eta}  \hookrightarrow {\rm Irr}  \overline{H}_1  \]
whose image consists of those irreducible genuine representations of $\overline{H}_1$ whose restriction to $Z$ is a prescribed character $\chi_{\xi}$.  If the LLC holds for the case of trivial $\eta$, then one would have a map
\[ \mathcal{L} \circ  \xi^* \circ p^*:  {\rm Irr} \overline{G}_{\eta}   \longrightarrow  \{  \text{splittings of ${^L}\overline{H}_1$} \}. \]
Now recall that by Lemma \ref{L:z-ext}, there is a natural map
\[ p:  {^L}\overline{H}_1  \longrightarrow Z^{\vee}. \]
If one assumes the (weak) LLC for $\overline{H}_1$ satisfies the natural property that the restriction of the central character of $\pi \in {\rm Irr}(\overline{H}_1)$ to $Z$ corresponds  to the parameter $p \circ \mathcal{L}(\pi)$ under the usual LLC for the (linear) torus $Z$, then one sees that 
\[
  {\rm Irr} \overline{G}_{\eta}   \longrightarrow  \{  \text{splittings $\phi$ of ${^L}\overline{H}_1$:  $p \circ \phi$ corresponds to $\chi_{\xi}$} \} = \{ \text{splittings of ${^L}G_{\eta}$} \}. \]

 \vskip 5pt

\subsection{\bf Reduction to discrete series.}
 The existence of the (weak) LLC map $\mathcal{L}$ can be reduced to the case of (quasi)-discrete series representation, much like the case of linear reductive groups.  More precisely, Ban-Jantzen \cite{BJ} have established the analog of the Langlands classification for general covering groups. This says that every irreducible representation is uniquely expressed as the unique irreducible quotient of a standard module.  As in the linear case, this reduces the definition of $\mathcal{L}$ to (quasi)-tempered representations. On the other hand, as shown in the work of W.W. Li \cite{L4}, any tempered representation is  contained in a representation parabolically induced from (quasi-)discrete series representations; moreover the decomposition of this induced representation is governed by an R-group.  This reduces the definition of $\mathcal{L}$ to the case of (quasi-)discrete series representations. 
\vskip 10pt

\section{\bf   Automorphic $L$-functions}
While we have considered only the case of local fields for most of this paper, we shall now consider the global setting, so that $k$ is a number field with ring of adeles $\A$. We shall briefly explain how the construction of the L-group  extension extends to the global situation, referring to \cite{W7}  for the details.
The goal of the section is to give  a definition of the notion of {\em automorphic L-functions}.
\vskip 5pt

\subsection{\bf Adelic BD covering}
Starting with a BD extension $\overline{\G}$ over ${\rm Spec} k$  and a positive integer $n$ such that $|\mu_n(k)| = n$, Brylinski and Deligne showed that  one inherits the following data:
\vskip 5pt

\begin{itemize}
\item for each place $v$ of $k$,   a local BD covering group $\overline{G}_v$ of degree $n$;
\vskip 5pt

\item for almost all $v$,  a splitting $s_v:  \G(\mathcal{O}_v) \longrightarrow \overline{G}_v$;
\vskip 5pt

\item a restricted direct product $\prod'_v \overline{G}_v$ with respect to the family of subgroups $s_v( \G(\mathcal{O}_v) )$, from which one can define:
\[  \overline{G}(\A): =  \prod'_v  \overline{G}_v /  \{  (\zeta_v) \in \oplus_v \mu_n(k_v):  \prod_v \zeta_v  =1 \}, \]
which gives a topological central extension
\[  \begin{CD}
\mu_n(k)  @>>>  \overline{G}(\A) @>>> G(\A). \end{CD} \]
We shall call this the adelic or global BD covering group;
\vskip 5pt
\item a natural inclusion
\[  \begin{CD}
\mu_n(k_v)  @>>>  \overline{G}_v @>>>  G(k_v) \\
@| @VVV  @VVV \\
\mu_n(k)  @>>>  \overline{G}(\A)  @>>>  G(\A)  \end{CD} \]
for each place $v$ of $k$;
\vskip 5pt

\item a natural splitting 
\[   i:  G(k)  \longrightarrow   \overline{G}(\A),\]
which allows one to consider the space of automorphic forms on $\overline{G}(\A)$.

\end{itemize}

\vskip 5pt

\vskip 10pt

\subsection{\bf  Global L-group extension}
One may define the L-group extension for the adelic BD cover $\overline{G}(\A)$ following the same procedure as in the local setting. We briefly summarise the process, highlighting the differences. 
Suppose that $\overline{\G}$ has BD invariant $(Q, \mathcal{E}, f)$ or bisector data  $(D, \eta)$. Then one has:
\vskip 5pt

\begin{itemize}
\item The dual group of $\overline{G}(\A)$ is defined in exactly the same way as in the local setting. 
Namely,  one may define the lattice $Y_{Q,n}$ and  the modified coroot lattice $Y_{Q,n}^{sc}$ in the same way.  This gives the dual group $\overline{G}^{\vee}$. Indeed, since $\G$ is split,  the definition of these objects work over any $k$ or $k_v$  and gives the same complex dual group $\overline{G}^{\vee}$.
\vskip 5pt

\item The role of $F^{\times}/ F^{\times n}$ in the local setting is replaced by $\A^{\times} / k^{\times} \A^{\times n}$. More precisely, 
with $(-,-)_n$ denoting the global n-th Hilbert symbol, the 2-cocycle
\[  c_1(a,b) (y)  = (a,b)_n^{Q(y)} \]
for $a,b \in  \A^{\times}/  \A^{\times n}$ defines an extension
 \[  \begin{CD}
  \Hom( Y_{Q,n}/  Y_{Q,n}^{sc}, \C^{\times})  @>>> E_1'  @>>>  \A^{\times}/  \A^{\times n}   \end{CD} \]
 Since $c_1$ is trivial on $k^{\times} \times k^{\times}$, this sequence splits canonically over the image of $k^{\times}$ in $  \A^{\times}/ \A^{\times n}$.
 Dividing out by the image of $k^{\times}$ under the splitting gives the extension $E_1$:
 \[  \begin{CD}
    \Hom( Y_{Q,n}/  Y_{Q,n}^{sc}, \C^{\times})  @>>> E_1  @>>>  \A^{\times}/ k^{\times} \A^{\times n}   \end{CD} \]
 
 \vskip 5pt

\item The extension $E_2$ is defined in the same way, applying $\Hom(-, \C^{\times})$ to the sequence
\[  \begin{CD}
k^{\times}/  k^{\times n}  @>>>  \mathcal{E}/  s_f( Y_{Q,n}^{sc})  @>>> Y_{Q,n}/ Y_{Q,n}^{sc},  \end{CD} \]
and pulling back by 
\[ \A^{\times}/  k^{\times} \A^{\times n}   \longrightarrow  \Hom(k^{\times}/ k^{\times n}, \C^{\times}). \]
Forming Baer sum with $E_1$ gives the global fundamental sequence:
\[  \begin{CD} 
Z(\overline{G}^{\vee})  @>>>  E  @>>>  \A^{\times}/ k^{\times} \A^{\times n}. \end{CD} \]
Pulling back to the global Weil group $W_k$ and pushing out to $\overline{G}^{\vee}$ gives the global L-group extension, which fits into a commutative diagram:
\[  \begin{CD}
\overline{G}^{\vee} @>>>  {^L}\overline{G}_v  @>>>  W_{k_v} \\
@| @VVV  @VVV \\
\overline{G}^{\vee}  @>>>  {^L}\overline{G}_{\A}  @>>>  W_k  \end{CD} \]
for each place $v$ of $k$.
\vskip 5pt

\item This global L-group extension satisfies the same functoriality with respect to Levi subgroups and $z$-extensions as in the local case. 
\end{itemize}

\vskip 5pt
\subsection{\bf Distinguished splittings.}
One may examine the question of splitting of the global fundamental sequence. 
As in the local case, this amounts to finding a genuine automorphic character 
\[  \chi:  T_{Q,n}(k) \cdot T_J(\A) \backslash \overline{T}_{Q,n}(\A)  \longrightarrow \C^{\times} \]
where $J  = Y_{Q,n}^{sc} +  nY \subset Y_{Q,n}$. 
Such a character exists  when $\eta_n  =1$ is trivial, and thus a distinguished splitting of the fundamental sequence $E$ exists in this case. If we fix such an automorphic character $\chi  = \prod_v \chi_v$, then 
each $\chi_v$ corresponds to a distinguished splitting of the local fundamental sequence. Moreover, 
$\chi$ is invariant under the Weyl group $W(k)  = N(\T)(k)/  \T(k)$. 
The explicit construction of a distinguished genuine character given in \S \ref{S:construction} produces an automorphic character.
\vskip 5pt

We deduce from the above discussion that the L-group extension is always split, and thus is abstractly a semi-direct product, but it may not be a direct product in general.

\vskip 5pt

\subsection{\bf  Automorphic L-functions}
We now have all the ingredients to define the notion of (partial) automorphic L-functions.
Let $\pi = \otimes_v  \pi_v$ be a cuspidal automorphic representation of $\overline{G}_{\A}$. 
For almost all $v$, $\pi_v$ is $s_v$-unramified. By the unramified LLC, $\pi_v$ gives rise to an $s_v$-unramified splitting 
\[  \rho_{\pi,v}:  W_{k_v}  \longrightarrow {^L}\overline{G}_v  \subset {^L}\overline{G}_{\A}. \]  
\vskip 5pt

 \vskip 5pt

Thus, if $R :  {^L}\overline{G}_{\A} \longrightarrow  \GL(V)$ is any finite dimensional representation, we may form the local Artin L-factor  for the representation $R \circ \rho_{\pi,v}:  W_{k_v} \longrightarrow \GL(V)$:
\[  L(s, \pi_v, R)  =  \frac{1}{\det(1-  \rho_{\pi,v}(Fr_v)  \cdot q_v^{-s} |V^{I_v}) }.  \]
Then we may form the  partial global L-function relative to $R$:
\[  L^S(s, \pi, R)  = \prod_{v \notin S}  L(s, \pi_v, R). \]

\vskip 5pt

 We want to highlight some instances where one can write down these automorphic L-functions. 
\vskip 5pt

\begin{itemize}
\item[(i)] (L-functions relative to a distinguished splitting) 
 If ${^L}\overline{G}_{\A}$ possesses a distingsuihed splitting $\rho_0$  (e.g. if $\eta_n  =1$), then 
$\rho_0$ is $s_v$-unramified for almost all $v$, and so we have an unramified homomorphism
\[  \rho_{\pi,v}/  \rho_{0,v}  :  W_{k_v}  \longrightarrow k_v^{\times}  \longrightarrow \Z  \longrightarrow \overline{G}^{\vee}.   \]
In other words, for almost all $v$, we have a Satake parameter $s_{\pi_v} \in  \overline{G}^{\vee}$, well-defined up to conjugacy, and depending on $\rho_{0,v}$.
\vskip 5pt

In this setting, if one has a representation $R:  \overline{G}^{\vee}  \longrightarrow \GL(V)$, one can form the partial L-function
\[  L^S(s, \pi,  R, \rho_0)
  :=  \prod_{v \notin S}  \frac{1}{\det(1-  s_{\pi_v}  \cdot q_v^{-s} |V) } . \]
  Because $\pi_v$ is unitary, it follows that this (partial) Euler product converges when $Re(s)$ is sufficiently large. We call this the $(R, \rho_0)$ L-function of $\pi$.  
  \vskip 5pt
  
  More generally, if one fixes a distinguish splitting $\rho_0$ of the enlarged L-group ${^L}\overline{G}^{\#}$ (which always exists), one has the notion of unramified L-parameters relative to $\rho_0$:
  \[  \rho_{\pi,v}: W_{k_v}  \rightarrow k_v^{\times} \rightarrow \Z \rightarrow \overline{G}^{\#},  \]  
  which gives rise to a Satake parameter $s_{\pi_v}  \in \overline{G}^{\#}$.  If one extends $R$ above to the disconnected group $\overline{G}^{\#}$, 
  one can define the partial L-function $L^S(s, \pi, R, \rho_0)$ as above. 
  \vskip 5pt

\item[(ii)]  (Adjoint type L-functions)   
If $\overline{G}^{\vee}_{ad}  := \overline{G}^{\vee}/  Z(\overline{G}^{\vee})$ denotes the adjoint quotient of $\overline{G}^{\vee}$, there is a natural commutative diagram of extensions
 \[  \begin{CD}
 \overline{G}^{\vee}  @>>>  {^L}\overline{G}_{\A} @>>>  W_k  \\
 @VVV  @VVpV  @VVV  \\
 \overline{G}^{\vee}_{ad}  @>>>   \overline{G}^{\vee}_{ad}  \times W_k  @>>>  W_k.  \end{CD}  \] 
  Thus, if  $R:   \overline{G}^{\vee}_{ad}\longrightarrow \GL(V)$ is any  representation, we may pull it back to ${^L}\overline{G}_{\A}$ and obtain a partial L-function $L^S(s, \pi, R)$.    
\vskip 5pt

\item[(iii)]  (Langlands-Shahidi  L-functions)
More generally, suppose that $\mathbb{P} = \mathbb{M} \cdot \mathbb{N}  \subset \G$ is a Levi subgroup, and $\pi$ is an automorphic representation of the BD covering $\overline{M}_{\A}$.  
By functoriality of the L-group construction, one has inclusions
 \[   E_{\overline{G}} \hookrightarrow  {^L}\overline{M}_{\A}  \hookrightarrow {^L}\overline{G}_{\A} \]
 where $E_{\overline{G}}$ is the fundamental sequence for $\overline{G}$. 
 As in (ii) above, one has a natural commutative diagram of extensions
  \[  \begin{CD}
 \overline{M}^{\vee} @>>>   {^L}\overline{M}_{\A}  @>>> W_k  \\
 @VVV  @VVV  @VVV  \\
 \overline{M}^{\vee}/  Z(\overline{G}^{\vee}) @>>> \overline{M}^{\vee}/  Z(\overline{G}^{\vee})  \times W_k  @>>> W_k, \end{CD}  \]
 i.e. a canonically split extension.  Then any representation of $\overline{M}^{\vee}$ which is trivial on $Z(\overline{G}^{\vee})$ pulls back to a representation of ${^L}\overline{M}_{\A}$. 
 \vskip 5pt
 
 A source of such representations is the adjoint action of ${^L}\overline{M}_{\A}$ on ${\rm Lie}( N^{\vee})$. Let  $R$ be  an irreducible summand, so that  $R$ is trivial on $Z(\overline{G}^{\vee})$. 
Then we obtain a partial automorphic L-function $L^S(s, \pi, R)$.  
As shown in the PhD thesis \cite{Ga} of the second author,
these Langlands-Shahidi type L-functions appear in the constant term of the Eisenstein series on $\overline{G}_{\A}$, as in the case of linear groups.
  \end{itemize}
 \vskip 15pt
 
 A basic open question is whether such automorphic L-functions associated to automorphic representations of BD covering groups have the usual nice properties such as meromorphic continuation and functional equations. A related question is whether such an automorphic L-function agrees with one for a linear reductive group.

\vskip 15pt

\section{\bf Examples}
We conclude this paper by giving a number of examples to illustrate the above discussion.  These examples are the ones which have been studied in the literature. As these groups arise as the cover $\overline{\G}$ of $\G$ which has simply-connected derived group, we may assume that $\overline{\G}$ is incarnated by a fair $(D, 1)$ without loss of generalities. For  fixed $n\in \N$, we have the associated degree $n$ cover $\overline{G}$.
\vskip 5pt

We have seen that there always exist distinguished splittings of ${^L}\overline{G}$ for such $\overline{G}$, with respect to which ${^L}\overline{G} \simeq \overline{G} \times W_F$. In this section, we use $\chi_\psi$ to denote a distinguished character constructed in section \ref{S:construction}. It will be shown explicitly that our construction in the simply-connected simply-laced case agrees with the one given by Savin \cite{Sa}. It is also compatible with the one for the classical double cover $\overline{\Sp}_{2r}$, as in \cite{K, Ra}.
\vskip 5pt

The computation of the bilinear form $B_Q$ below uses crucially the identity
$$B_Q(\alpha^\vee, y)=Q(\alpha^\vee)\cdot \langle \alpha, y \rangle,$$
where $\alpha \in \Phi^\vee$ is any coroot and $y\in Y$.

\vskip 5pt

\subsection{\bf Simply-connected case}
Consider a simply-connected simple group $\G$ of arbitrary type. There is up to unique isomorphism a $\mathbb{K}_2$-torsor $\overline{\G}$ associated to a Weyl-invariant quadratic form on $Y^{sc}=Y$. Consider $\overline{\G}$ incarnated by $(D, \eta)$. As indicated above, there is no loss of generality in assuming $D$ fair and $\eta=1$, and we will do so in the following.

For simplicity we assume $n=2$ except for the case of the exceptional $G_2$ where the computation is very simple for general $n$. We also assume that $Q$ is the unique Weyl-invariant quadratic form which takes value 1 on the short coroots of $\G$. The general case of $n$ and $Q$ follows from similar computations.

Note that whenever we have assumed $n=2$, we will write $Y_{Q,2}$ and $Y_{Q,2}^{sc}$ for the lattices $Y_{Q,n}$ and $Y_{Q,n}^{sc}$ which are of interest. We also have $J=2Y+Y_{Q,2}^{sc}=Y_{Q,2}^{sc}$ since $Y=Y^{sc}$.

\subsubsection{\bf The simply-laced case $A_r, D_r, E_6, E_7, E_8$ and compatibility}
Now let $\G$ be a simply-laced simply-connected group of type $A_r$ for $ r\ge 1$, $D_r$ for $r\ge 3$, and $E_6, E_7, E_8$. Let $\Delta=\{\alpha_1, ..., \alpha_r\}$ be a fixed set of simple roots of $\G$. Let $\overline{\G}$ be the extension of $\G$ determined by the quadratic form $Q$ with $Q(\alpha_i^\vee)=1$ for all coroots $\alpha_i^\vee$. 
We obtain the two-fold cover $\overline{G}$ of $G$.

Clearly we have $n_\alpha=2$ for all  $\alpha \in \Phi$ in this case. Let $\alpha_i^\vee \in \Delta^\vee$ for $i=1, ..., r$ be the simple coroots of $\G$. It is easy to compute the bilinear form $B_Q$ associated with $Q$: 
\begin{equation} \label{sl B}
    B_Q(\alpha_i^\vee, \alpha_j^\vee) =
    \begin{cases}
      -1 & \text{if } \alpha_i \text{ and } \alpha_j \text{ connected in the Dynkin diagram,} \\
      0 & \text{otherwise}.
    \end{cases}
\end{equation}

In order to show compatibility with Savin, we may further assume that $\overline{G}$ is incarnated by the following fair bisector $D$ associated with $B_Q$ as given in \cite{Sa},
\begin{equation} \label{S-fair D}
    D(\alpha_i^\vee, \alpha_j^\vee) =
    \begin{cases}
      0 & \text{ if } i< j, \\
Q(\alpha^\vee_i) &\text{ if } i=j, \\
      B_Q(\alpha_i^\vee, \alpha_j^\vee) & \text{ if } i>j.
    \end{cases}
\end{equation}

The following lemma is in \cite{Sa} and reproduced here for convenience. The stated result can also be checked by straightforward computation.
\begin{lemma} \label{geom lm}
Let $\Omega$ be a subset of the vertices in the Dynkin diagram of $\G$ satisfying:
\begin{enumerate}
\item[(i)] No two vertices in $\Omega$ are adjacent,

\item[(ii)] Every vertex not in $\Omega$ is adjacent to an even number of vertices in $\Omega$. 
\end{enumerate}

Then the map given by $\xymatrix{\Omega \ar@{|->}[r] & e_\Omega}$ with $e_\Omega:=\sum_{\alpha_i \in \Omega} \alpha_i^\vee$ gives a well-defined correspondence between such sets $\Omega$ and the cosets of $Y_{Q,2}/J$. In particular, the empty set  corresponds to the trivial coset $J$.
\end{lemma}
By properties of $B_Q$ and (i) of $\Omega$ above, it follows that
$$Q(e_\Omega)=|\Omega|.$$

We now give a brief case by case discussion.

\vskip 10pt

\textbf{The $A_r$ case.} There are two situations according to the parity of $r$.
\vskip 10pt

\underline{Case 1:} $r$ is even.  As an illustration, we first do the straightforward computation. Let $\sum_i k_i \alpha_i^\vee \in Y_{Q,2}$ for proper $k_i\in \Z$. Then $B_Q(\sum_i k_i \alpha_i^\vee, \alpha^\vee_j) \in 2\Z$ for all $1\le j\le r$ by the definition of $Y_{Q,2}$. In view of (\ref{sl B}), it is equivalent to 
\begin{equation} \label{comp}
\begin{cases}
2k_1 + (-1) k_2 &\in 2\Z, \\
(-1)k_1 + 2 k_2 + (-1) k_3 &\in 2\Z, \\
(-1)k_2 + 2 k_3 + (-1) k_4 &\in 2\Z, \\
\vdots & \ \\
(-1)k_{r-2} + 2 k_{r-1} + (-1) k_r &\in 2\Z, \\
(-1)k_{r-1} + 2k_r &\in 2\Z.
\end{cases}
\end{equation}
It follows that $k_2$ is even and so are the successive $k_4, ..., k_r$ (we have assumed $r$ to be even). Similarly, $k_{r-1}$ is even and therefore all $k_{r-3}, ..., k_1$ are also even. This gives $Y_{Q,2}=J$.

Note that we could simply apply the lemma to get $Y_{Q,2}=J$, which corresponds to the fact that only the empty set satisfies properties $(i)$ and $(ii)$. There is nothing to check in this case, and the character $\chi$ such that $\chi \circ s_\eta$ is trivial will be a distinguished character.
\vskip 10pt

\underline{Case 2:} $r=2k-1$ is odd. Recall the notation 
$$\alpha_{i, Q, 2}^\vee:=n_{\alpha_i} \cdot \alpha_i^\vee. $$

The consideration as in (\ref{comp}) works. However, for convenience we will apply the lemma to get $[Y_{Q,2}, J]=2$ with a basis of $Y_{Q,2}$ given by
$$\{\alpha_{r, Q, 2}^\vee, \alpha_{r-1, Q, 2}^\vee, ..., \alpha_{2, Q, 2}^\vee, e_\Omega=\sum_{m=1}^k \alpha_{2m-1}^\vee \}.$$
The nontrivial coset corresponds to the set $\Omega$ indicated by alternating bold circles below
$$
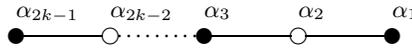
\begin{picture}(4.7,0.2)(0,0)
\put(1,0){\circle*{0.08}}
\put(1.5,0){\circle{0.08}}
\put(2,0){\circle*{0.08}}
\put(2.5,0){\circle{0.08}}
\put(3,0){\circle*{0.08}}
\put(1.04,0){\line(1,0){0.42}}
\multiput(1.55,0)(0.05,0){9}{\circle*{0.02}}
\put(2.04,0){\line(1,0){0.42}}
\put(2.54,0){\line(1,0){0.42}}
\put(1,0.1){\footnotesize $\alpha_{2k-1}$}
\put(1.5,0.1){\footnotesize $\alpha_{2k-2}$}
\put(2,0.1){\footnotesize $\alpha_3$}
\put(2.5,0.1){\footnotesize $\alpha_2$}
\put(3,0.1){\footnotesize $\alpha_1$}
\end{picture}
$$
A basis for $J=2Y^{sc}$ is given by
$$\{ \alpha_{r, Q, 2}^\vee, \alpha_{r-1, Q, 2}^\vee, ..., \alpha_{2, Q, 2}^\vee, 2 e_\Omega \}.$$

Now the construction of the distinguished $\chi_\psi$ in previous section is determined by
\begin{equation}
\begin{cases}
\chi_\psi(\alpha_{i, Q, 2}^\vee (a)) =1, \ 2\le i \le r, \\
\chi_\psi(e_\Omega(a))=\gamma_\psi(a)^{(2-1)Q(e_\Omega)}=\gamma_\psi(a)^{|\Omega|}.
\end{cases}
\end{equation}

This agrees with the formula in \cite{Sa}  when we substitute $a=\varpi$ in $\gamma_\psi(a)$ for $\psi$ of conductor $\mathcal{O}_F$. See \cite[pg 118]{Sa}.

\vskip 10pt

\textbf{The $D_r$ case.} We also have two cases.
\vskip 10pt

\underline{Case 1:} $r=2k-1$ is odd with $k\ge 2$. Then $[Y_{Q,2}, J]=2$ with the nontrivial $\Omega=\{ \alpha_1, \alpha_2 \}$:
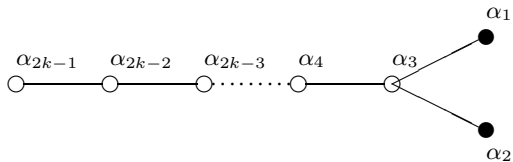
$$
\begin{picture}(4.7,0.4)(0,0)
\put(1,0){\circle{0.08}}
\put(1.5,0){\circle{0.08}}
\put(2,0){\circle{0.08}}
\put(2.5,0){\circle{0.08}}
\put(3,0){\circle{0.08}}
\put(3.5, 0.25){\circle*{0.08}}
\put(3.5, -0.25){\circle*{0.08}}
\put(1.04,0){\line(1,0){0.42}}
\put(1.54,0){\line(1,0){0.42}}
\multiput(2.05,0)(0.05,0){9}{\circle*{0.02}}
\put(2.54,0){\line(1,0){0.42}}
\put(3.00,0){\line(2,1){0.46}}
\put(3.00,0){\line(2,-1){0.46}}
\put(1,0.1){\footnotesize $\alpha_{2k-1}$}
\put(1.5,0.1){\footnotesize $\alpha_{2k-2}$}
\put(2,0.1){\footnotesize $\alpha_{2k-3}$}
\put(2.5,0.1){\footnotesize $\alpha_4$}
\put(3,0.1){\footnotesize $\alpha_3$}
\put(3.5,0.35){\footnotesize $\alpha_1$}
\put(3.5,-0.4){\footnotesize $\alpha_2$}
\end{picture}
$$
\vspace{0.5cm}

Consider the basis $\{ \alpha_{i,Q,2}^\vee: 2\le i \le r \} \cup \{ e_\Omega \}$ of $Y_{Q,2}$, then the construction of distinguished $\chi_\psi$ in previous section is determined by
\begin{align*}
\chi_\psi(e_\Omega(a))=\gamma_\psi(a)^{(2-1)Q(e_\Omega)}=\gamma_\psi(a)^{|\Omega|}.
\end{align*}
\vskip 10pt

\underline{Case 2:} $r=2k$ is even. Then $[Y_{Q,2}: J]=4$. There are three nontrivial sets $\Omega_i$ for $i=1, 2, 3$ as indicated by the bold circles below.
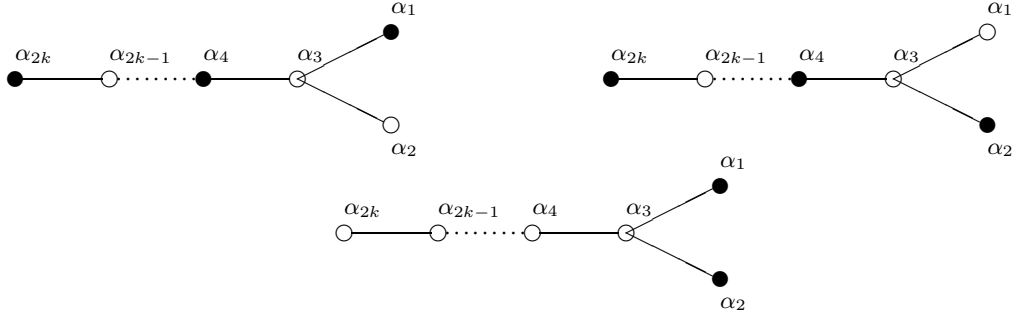
\begin{multicols}{2}
$$
\begin{picture}(4.7,0.1)(1,0)
\put(1.5,0){\circle*{0.08}}
\put(2,0){\circle{0.08}}
\put(2.5,0){\circle*{0.08}}
\put(3,0){\circle{0.08}}
\put(3.5, 0.25){\circle*{0.08}}
\put(3.5, -0.25){\circle{0.08}}
\put(1.54,0){\line(1,0){0.42}}
\multiput(2.05,0)(0.05,0){9}{\circle*{0.02}}
\put(2.54,0){\line(1,0){0.42}}
\put(3.00,0){\line(2,1){0.46}}
\put(3.00,0){\line(2,-1){0.46}}
\put(1.5,0.1){\footnotesize $\alpha_{2k}$}
\put(2,0.1){\footnotesize $\alpha_{2k-1}$}
\put(2.5,0.1){\footnotesize $\alpha_4$}
\put(3,0.1){\footnotesize $\alpha_3$}
\put(3.5,0.35){\footnotesize $\alpha_1$}
\put(3.5,-0.4){\footnotesize $\alpha_2$}
\end{picture}
$$
\goodbreak
$$
\begin{picture}(4.7,0.1)(1,0)
\put(1.5,0){\circle*{0.08}}
\put(2,0){\circle{0.08}}
\put(2.5,0){\circle*{0.08}}
\put(3,0){\circle{0.08}}
\put(3.5, 0.25){\circle{0.08}}
\put(3.5, -0.25){\circle*{0.08}}
\put(1.54,0){\line(1,0){0.42}}
\multiput(2.05,0)(0.05,0){9}{\circle*{0.02}}
\put(2.54,0){\line(1,0){0.42}}
\put(3.00,0){\line(2,1){0.46}}
\put(3.00,0){\line(2,-1){0.46}}
\put(1.5,0.1){\footnotesize $\alpha_{2k}$}
\put(2,0.1){\footnotesize $\alpha_{2k-1}$}
\put(2.5,0.1){\footnotesize $\alpha_4$}
\put(3,0.1){\footnotesize $\alpha_3$}
\put(3.5,0.35){\footnotesize $\alpha_1$}
\put(3.5,-0.4){\footnotesize $\alpha_2$}
\end{picture}
$$
\end{multicols}

$$
\begin{picture}(4.7,0.4)(0,0)
\put(1.5,0){\circle{0.08}}
\put(2,0){\circle{0.08}}
\put(2.5,0){\circle{0.08}}
\put(3,0){\circle{0.08}}
\put(3.5, 0.25){\circle*{0.08}}
\put(3.5, -0.25){\circle*{0.08}}
\put(1.54,0){\line(1,0){0.42}}
\multiput(2.05,0)(0.05,0){9}{\circle*{0.02}}
\put(2.54,0){\line(1,0){0.42}}
\put(3.00,0){\line(2,1){0.46}}
\put(3.00,0){\line(2,-1){0.46}}
\put(1.5,0.1){\footnotesize $\alpha_{2k}$}
\put(2,0.1){\footnotesize $\alpha_{2k-1}$}
\put(2.5,0.1){\footnotesize $\alpha_4$}
\put(3,0.1){\footnotesize $\alpha_3$}
\put(3.5,0.35){\footnotesize $\alpha_1$}
\put(3.5,-0.4){\footnotesize $\alpha_2$}
\end{picture}
$$

\vspace{1.3cm}

That is, $\Omega_1=\{\alpha_1\} \cup \{\alpha_{2m}: 2\le m\le k \}$, $\Omega_2=\{ \alpha_2 \}\cup \{\alpha_{2m}: 2\le m\le k \}$ and $\Omega_3=\{\alpha_1, \alpha_2 \}$. Note $|\Omega_1|=|\Omega_2|=k$ and $|\Omega_2|=2$.

A basis of $Y_{Q,2}$ is given by
$$\{\alpha_{i,Q,2}^\vee: 3\le i\le 2k-1 \}\cup \{e_{\Omega_1}, e_{\Omega_2}, e_{\Omega_3} \}.$$
However, the construction of distinguished characters utilizes the elementary divisor theorem. Thus we have to provide bases for $Y_{Q,2}$ and $J$ aligned in a proper way. To achieve this, consider the alternative basis of $Y_{Q,2}$ given by
$$\{\alpha_{i,Q,2}^\vee: 3\le i\le 2k-1\} \cup \{ e_{\Omega_1}+e_{\Omega_2} + e_{\Omega_3}, e_{\Omega_2} +e_{\Omega_3}, e_{\Omega_3} \}.$$
Then it is easy to check that the set
$$\{\alpha_{i,Q,2}^\vee: 3\le i\le 2k-1\} \cup \{e_{\Omega_1}+e_{\Omega_2} + e_{\Omega_3}, 2(e_{\Omega_2} +e_{\Omega_3}), 2e_{\Omega_3} \}$$
is a basis for $J$.
Note 
$$Q(e_{\Omega_2} +e_{\Omega_3})=|\Omega_2|+Q(2\alpha^\vee)=|\Omega_2|+4.$$

Thus a distinguished character could be determined by
\begin{equation*}
\begin{cases}
\chi_\psi(\alpha_{i,Q,2}^\vee(a))=1, \ 3\le i\le 2k-1, \\
\chi_\psi((e_{\Omega_1}+e_{\Omega_2} + e_{\Omega_3})(a))=1,\\
\chi_\psi((e_{\Omega_2} +e_{\Omega_3})(a))=\gamma_\psi(a)^{Q(e_{\Omega_2} +e_{\Omega_3})}= \gamma_\psi(a)^{|\Omega_2|}, \\
\chi_\psi(e_{\Omega_3}(a))=\gamma_\psi(a)^{|\Omega_3|}.
\end{cases}
\end{equation*}

However, since we have assumed that $D$ takes the special form given by (\ref{S-fair D}), we have
\begin{align*}
& D(e_{\Omega_1}, e_{\Omega_2} + e_{\Omega_3}) = |\Omega_1|\\
& D(e_{\Omega_2}, e_{\Omega_3}) =Q(\alpha_i^\vee) =1.
\end{align*}

Thus
\begin{align*}
\chi_\psi(e_{\Omega_1}(a)) \cdot \chi_\psi((e_{\Omega_2} +e_{\Omega_3})(a)) &= (a, a)_2^{|\Omega_1|} \cdot
\chi_\psi((e_{\Omega_1}+ e_{\Omega_2} +e_{\Omega_3})(a)) \\
\chi_\psi((e_{\Omega_2}(a)) \cdot \chi_\psi(e_{\Omega_3}(a)) &= (a, a)_2 \cdot
\chi_\psi((e_{\Omega_2} +e_{\Omega_3})(a)) 
\end{align*}

Recall $\gamma_\psi(a)^2=(a, a)_2$. This combined with the above results gives
\begin{equation*}
\begin{cases}
\chi_\psi(e_{\Omega_1}(a))=\gamma_\psi(a)^{|\Omega_1|}, \\
\chi_\psi(e_{\Omega_2}(a))=\gamma_\psi(a)^{|\Omega_2|}, \\
\chi_\psi(e_{\Omega_3}(a))=\gamma_\psi(a)^{|\Omega_3|}.
\end{cases}
\end{equation*}

It agrees with the genuine character given by Savin.

\vskip 10pt

\textbf{The $E_6, E_7, E_8$ case.} 

\vskip 10pt

For $E_6$ and $E_8$, $Y_{Q,2}=J$ and so the situation is trivial. Consider $E_7$, then $[Y_{Q,n}: J]=2$. The nontrivial $\Omega$ is given by $\Omega=\{\alpha_4, \alpha_6, \alpha_7 \}$.
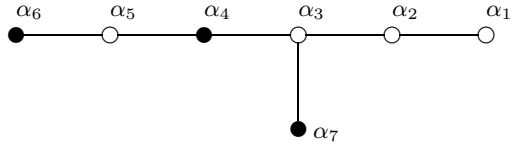
$$
\begin{picture}(4.7,0.2)(0,0)
\put(1,0){\circle*{0.08}}
\put(1.5,0){\circle{0.08}}
\put(2,0){\circle*{0.08}}
\put(2.5,0){\circle{0.08}}
\put(3,0){\circle{0.08}}
\put(3.5,0){\circle{0.08}}
\put(2.5,-0.5){\circle*{0.08}}
\put(1.04,0){\line(1,0){0.42}}
\put(1.54,0){\line(1,0){0.42}}
\put(2.04,0){\line(1,0){0.42}}
\put(2.5,-0.04){\line(0,-1){0.42}}
\put(2.54,0){\line(1,0){0.42}}
\put(3.04,0){\line(1,0){0.42}}
\put(1,0.1){\footnotesize $\alpha_6$}
\put(1.5,0.1){\footnotesize $\alpha_5$}
\put(2,0.1){\footnotesize $\alpha_4$}
\put(2.5,0.1){\footnotesize $\alpha_3$}
\put(3,0.1){\footnotesize $\alpha_2$}
\put(3.5,0.1){\footnotesize $\alpha_1$}
\put(2.58,-0.55){\footnotesize $\alpha_7$}
\end{picture}
$$
\vspace{1.2cm}

The set $\{\alpha_{i,Q,2}^\vee: 1\le i\le 6 \} \cup \{e_\Omega \}$ is a basis of $Y_{Q,2}$, while $\{\alpha_{i,Q,2}^\vee: 1\le i\le 6\} \cup \{2 e_\Omega \}$ a basis for $J$.

Our distinguished character is determined by
$$\chi_\psi(e_\Omega(a))=\gamma_\psi(a)^{|\Omega|}.$$
This agrees with Savin also.

\vskip 10pt

\subsubsection{\bf The case $C_r$.}
\vskip 10pt

Let $\Sp_{2r}$ be the simply-connected simple group with Dynkin diagram:

$$
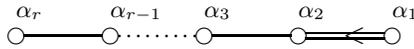
\begin{picture}(4.7,0.2)(0,0)
\put(1,0){\circle{0.08}}
\put(1.5,0){\circle{0.08}}
\put(2,0){\circle{0.08}}
\put(2.5,0){\circle{0.08}}
\put(3,0){\circle{0.08}}
\put(1.04,0){\line(1,0){0.42}}
\multiput(1.55,0)(0.05,0){9}{\circle*{0.02}}
\put(2.04,0){\line(1,0){0.42}}
\put(2.54,0.015){\line(1,0){0.42}}
\put(2.54,-0.015){\line(1,0){0.42}}
\put(2.74,-0.04){$<$}
\put(1,0.1){\footnotesize $\alpha_r$}
\put(1.5,0.1){\footnotesize $\alpha_{r-1}$}
\put(2,0.1){\footnotesize $\alpha_3$}
\put(2.5,0.1){\footnotesize $\alpha_2$}
\put(3,0.1){\footnotesize $\alpha_1$}
\end{picture}
$$
\

Let $\{\alpha_1^\vee, \alpha_2^\vee, ..., \alpha_r^\vee \}$ be the set of simple coroots with $\alpha_1^\vee$ the short one. Let $n=2$. Here $\overline{\Sp}_{2r}$ is determined by the unique Weyl-invariant quadratic form $Q$ on $Y$ with $Q(\alpha_1^\vee)=1$.

It follows $n_{\alpha_1}=2$. Also $Q(\alpha_i^\vee)=2$ and $n_{\alpha_i}=1$ for $2\le i\le  r$. Moreover, a basis of $Y_{Q,2}=Y^{sc}$ is given by
$$\{\alpha_1^\vee,  \alpha_2^\vee, ..., \alpha_r^\vee \}, $$

while a basis for $Y_{Q,2}^{sc}$ is
$$ \{ 2\alpha_1^\vee,  \alpha_2^\vee, ..., \alpha_r^\vee \}.$$

Since $J=Y_{Q,2}^{sc}$, by the construction of distinguished character $\chi_\psi$, it is determined by
\begin{equation}
\begin{cases}
\chi_\psi(1, \alpha^\vee_i(a)) =1, \text{ if } i=2, 3, ..., r; \\
\chi_\psi(1, \alpha^\vee_1(a)) =\gamma_\psi(a)^{(2-1)Q(\alpha^\vee_1)}=\gamma_\psi(a).
\end{cases}
\end{equation}

This uniquely determined a genuine character of $\overline{T}$ which is abelian. It can be checked  that this agrees with the classical one (cf. \cite{K, Ra} for example).

\subsubsection{\bf The $B_r$, $F_4$ and $G_2$ case.}

\vskip 10pt
 
For completeness, we also give the explicit form of the distinguished character constructed in previous section for the double cover  $\overline{G}$ of the simply connected group $G$ of type $B_r$, $F_4$ and $G_2$. Recall that when $n=2$ we have $J=Y_{Q,2}^{sc}$.

\vskip 5pt

\textbf{\bf The $B_r$ case.} Consider the Dynkin diagram of $B_r$:

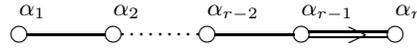
$$
\begin{picture}(4.7,0.2)(0,0)
\put(1,0){\circle{0.08}}
\put(1.5,0){\circle{0.08}}
\put(2,0){\circle{0.08}}
\put(2.5,0){\circle{0.08}}
\put(3,0){\circle{0.08}}
\put(1.04,0){\line(1,0){0.42}}
\multiput(1.55,0)(0.05,0){9}{\circle*{0.02}}
\put(2.04,0){\line(1,0){0.42}}
\put(2.54,0.015){\line(1,0){0.42}}
\put(2.54,-0.015){\line(1,0){0.42}}
\put(2.74,-0.04){$>$}
\put(1,0.1){\footnotesize $\alpha_1$}
\put(1.5,0.1){\footnotesize $\alpha_2$}
\put(2,0.1){\footnotesize $\alpha_{r-2}$}
\put(2.5,0.1){\footnotesize $\alpha_{r-1}$}
\put(3,0.1){\footnotesize $\alpha_r$}
\end{picture}
$$
\vskip 10pt

Let $Q$ be the unique Weyl-invariant quadratic form with $Q(\alpha^\vee_i)=1$ for $1\le i\le r-1$. It gives $Q(\alpha^\vee_r)=2$. We have also assumed that the double cover $\overline{G}$ is incarnated by a fair bisector $D$. The discussion now will be split into two cases according to the parity of $r$. 

\underline{Case 1:} $r$ is odd. Direct computation gives $Y_{Q,2}^{sc}=Y_{Q,n}$ and therefore this case is trivial.

\underline{Case 2:} $r$ is even. It is not difficult to compute the index $[Y_{Q,2}: Y_{Q,2}^{sc}]=2$. In fact, a basis of $Y_{Q,2}$ is given by
$$\{ \alpha^\vee_1+\alpha^\vee_3 +... +\alpha^\vee_{r-1} \} \cup \{ 2\alpha^\vee_i: 2\le i\le r-1 \} \cup \{\alpha_r^\vee \}.$$
This gives a basis of $J=Y_{Q,2}^{sc}$:
$$\{ 2(\alpha^\vee_1+\alpha^\vee_3 +... +\alpha^\vee_{r-1}) \} \cup \{2\alpha^\vee_i: 2\le i\le r-1 \} \cup \{\alpha_r^\vee \}.$$

We have 
$$Q(\alpha^\vee_1+\alpha^\vee_3 +... +\alpha^\vee_{r-1})=r/2.$$
By the construction of distinguished character $\chi_\psi$, it is determined by
\begin{equation}
\begin{cases}
 \chi_\psi\big((\alpha^\vee_1+\alpha^\vee_3 +... +\alpha^\vee_{r-1})(a)\big) =\gamma_\psi(a)^{r/2}; \\
 \chi_\psi\big((2\alpha_i^\vee)(a)\big) =1, \text{ for } 2\le i\le r-1; \\
 \chi_\psi\big(\alpha_r^\vee(a) \big) =1.
\end{cases}
\end{equation}

\vskip 5pt

\textbf{The $F_4$ case.} Consider the Dynkim diagram of $F_4$:

$$
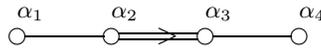
\begin{picture}(4.7,0.2)(0,0)
\put(2,0){\circle{0.08}}
\put(2.5,0){\circle{0.08}}
\put(3,0){\circle{0.08}}
\put(3.5,0){\circle{0.08}}
\put(2.04,0){\line(1,0){0.42}}
\put(2.54,0.015){\line(1,0){0.42}}
\put(2.54,-0.015){\line(1,0){0.42}}
\put(2.74,-0.04){$>$}
\put(3.04,0){\line(1,0){0.42}}
\put(2,0.1){\footnotesize $\alpha_1$}
\put(2.5,0.1){\footnotesize $\alpha_2$}
\put(3,0.1){\footnotesize $\alpha_3$}
\put(3.5,0.1){\footnotesize $\alpha_4$}
\end{picture}
$$

\vskip 10pt 

Let $Q$ be such that $Q(\alpha^\vee_i)=1$ for $i=1, 2$. It implies $Q(\alpha_i^\vee)=2$ for $i=3, 4$. Clearly $n_{\alpha_i}=2$ for $i=1, 2$ and $n_{\alpha_i}=1$ for $i=3, 4$. We can compute 
$$B_Q(\alpha^\vee_1, \alpha^\vee_2)=-1,\quad B_Q(\alpha^\vee_2, \alpha^\vee_3)=-2, \quad B_Q(\alpha^\vee_3, \alpha^\vee_4)=-2.$$ 
Also $B_Q(\alpha^\vee_i, \alpha^\vee_j)=0$ if $\alpha_i$ and $\alpha_j$ are not adjacent in the Dynkin diagram.

Moreover, any $\sum_i k_i\alpha^\vee_i \in Y^{sc}$ with certain $k_i \in \Z$ belongs to $Y_{Q,2}$ if and only if 
$$B_Q(\sum_i k_i\alpha^\vee_i, \alpha_j^\vee) \in 2\Z \text{ for all } 1\le j\le 4,$$
which explicitly is given by
\begin{equation*}
\begin{cases}
2k_1 + (-1) k_2 \in 2\Z, \\
(-1) k_1 + 2k_2 +(-2) k_3\in 2\Z, \\
(-2)k_2 +4k_3 +(-2)k_4 \in 2\Z, \\
(-2)k_3 +4k_4 \in 2\Z.
\end{cases}
\end{equation*}

Equivalently, $k_1, k_2\in 2\Z$. This shows $Y_{Q, 2}=Y_{Q,2}^{sc}$, and thus the situation is trivial. 
\vskip 5pt

\textbf{The $G_2$ case.} Consider the Dynkin diagram of $G_2$:
$$
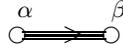
\begin{picture}(4.7,0.2)(0,0)
\put(2.5,0){\circle{0.08}}
\put(3,0){\circle{0.08}}
\put(2.53,0.018){\line(1,0){0.44}}
\put(2.54,0){\line(1,0){0.42}}
\put(2.53,-0.018){\line(1,0){0.44}}
\put(2.74,-0.035){$>$}
\put(2.5,0.1){\footnotesize $\alpha$}
\put(3,0.1){\footnotesize $\beta$}
\end{picture}
$$

\vskip 10pt 

Let $Q$ be such that $Q(\alpha^\vee)=1$. This determines $Q(\beta^\vee)=3$. Note $B_Q(\alpha^\vee, \beta^\vee)=-Q(\alpha^\vee)=-3$.

Since the computation is straightforward, we may assume $n \in \N_{\ge 1}$ is general instead of $2$. It follows $n_\alpha=n$ and $n_\beta=n/\text{gcd}(n, 3)$. Then $k_1\alpha^\vee + k_2\beta^\vee $ lies in $Y_{Q,n}$ if and only if
\begin{equation*}
\begin{cases}
2k_1 -3k_2 \in n\Z, \\
-3k_1 + 6k_2 \in n\Z.
\end{cases}
\end{equation*}

Equivalently, $k_1\in n\Z$ and $k_2$ divisible by $n/\text{gcd}(n,3)$. This exactly shows $Y_{Q,n}=Y_{Q,n}^{sc}$ for arbitrary $n$.  
Also in this case, it is trivial to define the distinguished character for the fair $D$.

\vskip 10pt
\subsection{\bf Kazhdan-Patterson coverings $\overline{\GL}_r$ \cite{KP1, KP2}} 
We consider the embedding of $\GL_r \hookrightarrow \SL_{r+1}$ given by $g\mapsto (g, \text{det}(g)^{-1})$. Using this embedding, we could describe the root data $(X, \Phi, Y, \Phi^\vee)$ of $\GL_r$ in terms of that of $\SL_{r+1}$. Let $\{\alpha_i^\vee: 1\le i\le r\}$ denote a set of simple roots of $\SL_{r+1}$. Then $\Delta^\vee=\{\alpha_i^\vee: 1\le i\le r-1\}$ will be a set of simple coroots of $\GL_r$. Moreover, a basis for $Y$ is given by $\{e_i\}, 1\le i\le r$ where
$$e_i:= \sum_{k\ge i}^r \alpha_k^\vee.$$

Let $\overline{\SL}_{r+1}$ be the degree $n$ covering associated with $(n, Q)$ where $Q(\alpha_i^\vee)=1$ for all $i$. Then the pull-back covering of $\GL_r$ will be the Kazhdan-Patterson $\overline{\GL}_r$, which we may assume is incarnated by a fair $(D, 1)$.

We first determine the inherited bilinear form $B_Q$ on $Y$:
\begin{equation*}
B(e_i, e_j)=
\begin{cases}
2 & \text{ if } i=j, \\
1 & \text{ otherwise}.
\end{cases}
\end{equation*}
Write $c_{n,r}:=n/\text{gcd}(r+1, n)$. It follows
$$Y_{Q,n}=\{\sum_{i=1}^r m_i e_i: m_i \equiv m_j \text{ mod } n, \text{ and }  c_{n,r}|m_i \text{ for all } i, j \}.$$
In particular, a basis for $Y_{Q,n}$ is given by
$$\{n e_i: 1\le i\le r-1\} \cup \{ c_{n,r} \cdot (\sum_{i=1}^r e_i)\}.$$
On the other hand, $Y_{Q,n}^{sc}$ is spanned by $\{n\cdot \alpha_i^\vee: 1\le i\le r-1\}$. It follows that $J=Y_{Q,n}^{sc} + nY$ has a basis given by
$$\{n e_i: 1\le i\le r-1\} \cup \{ n\cdot (\sum_{i=1}^r e_i)\}.$$

A distinguished character $\chi_\psi$ is thus determined by
\begin{equation*}
\begin{cases}
\chi_\psi\big( ne_i(a)\big) =1 \text{ for all } 1\le i\le r-1, \\
\chi_\psi\big((\sum_{i=1}^r c_{n,r} e_i)(a)\big)=\gamma_\psi(a)^{\frac{r(r+1)c_{n,r}}{n}(n-c_{n,r})}.
\end{cases}
\end{equation*}

An examination of the root datum shows that
\[  \overline{G}^{\vee}  \cong \{  (g, \lambda) \in \GL_r(\C) \times \GL_1(\C): \det(g)  = \lambda^{(r+1,n)} \} \subset \GL_r(\C)  \times \GL_1(\C). \]  
Thus, in general, the dual group $\overline{G}^\vee$ may not be $\GL_r$. For example, consider $n=2$. If $r$ is even, then
$$\overline{\GL}_r^\vee=\GL_r.$$
If $r$ is odd, then $\overline{\GL}_r^\vee \subseteq \GL_r \times \GL_1$ and is given by
$$\overline{\GL}_r^\vee=\{(g, a)\in \GL_r\times \GL_1: \text{det}(g)=a^2 \}.$$
When $r=1$ or $r=3$ there is an isomorphism $\overline{\GL}_r^\vee \simeq \GL_r$ given by $(g, a)\mapsto g a^{-1}$. However, for odd $r\ge 5$, there exists an isogeny $\overline{\GL}_r^\vee \to \GL_r$ of degree two given by $(g, a) \mapsto g$.

 \vskip 10pt

\subsection{\bf The cover $\overline{\GSp}_{2r}$}

Let $\G$ be the group $\GSp_{2r}$ of similitudes of symplectic type, and let $(X, \Delta, Y, \Delta^\vee)$ be its root data given as follows. The character group $X\simeq \Z^{r+1}$ has a standard basis $\{e_i^*: 1\le i\le r \} \cup \{e_0^* \}$, and the roots are given by 
$$\Delta=\{e_i^*-e_{i+1}^*: 1\le i \le r-1 \} \cup \{ 2e_r^*-e_0^* \}.$$

The cocharacter group $Y\simeq \Z^{r+1}$ is given with a basis $\{e_i: 1\le i\le r \} \cup \{e_0 \}$. The coroots are 
$$\Delta^\vee=\{ e_i-e_{i+1}: 1\le i \le r-1 \} \cup \{ e_r \}.$$

Write $\alpha_i=e_i^* - e_{i+1}^*$, $\alpha^\vee_i=e_i-e_{i+1}$ for $1\le i\le r-1$, and also $\alpha_r=2e_r^*-e_0^*$, $\alpha_r^\vee=e_r$. Consider a covering $\overline{\G}$ incarnated by $(D, 1)$. We are interested in those $\overline{\G}$ whose restricted to $\Sp_{2r}$ is the one with $Q(\alpha_r^\vee)=1$. That is, we assume
$$Q(\alpha_i^\vee)=2 \text{ for } 1\le i \le r-1, \quad Q(\alpha_r^\vee)=1.$$

Since $\Delta^\vee \cup \{e_0\}$ gives a basis for $Y$, to determine $Q$ it suffices to specify $Q(e_0)$. Let $n=2$, and we obtain a double cover $\overline{\GSp}_{2r}$ which restricts to the classical metaplectic double cover $\overline{\Sp}_{2r}$. Note also the number $Q(e_0) \in \Z/2\Z$ determines whether the similitude factor $F^\times$ corresponding to the cocharacter $e_0$ splits into $\overline{\GSp}_{2r}$ or not. To recover the classical double cover of $\GSp_{2r}$, we should take $Q(e_0)$ to be even.

Back to the case of $n=2$ and general $Q(e_0)$. We compute the root data for the complex dual group $\overline{\GSp}_{2r}^{\vee}$. We have
$$Y_{Q,2}=\{\sum_{i=1}^r k_i \alpha_i^\vee +k e_0 \in Y: k_i\in \Z \text{ for } 1\le i\le r-1, k_r, k \in 2\Z \}$$
and the sublattice $Y_{Q,2}^{sc}$ is spanned by $\{\alpha_{i, Q, 2}^\vee\}_{1\le i\le r}$, i.e.
$$\{\alpha_1^\vee, \alpha_2^\vee, ..., \alpha_{r-1}^\vee, 2\alpha_r^\vee \}.$$

The lattice $J=Y_{Q,2}^{sc}+ 2Y$ is thus equal to $Y_{Q,2}$. In this case, distinguished character $\chi_\psi$ will be determined by the condition (since we have assumed $(D,1)$ to be fair)
$$\chi_\psi\big(y(a)\big)=1 \text{ for all } y \in Y_{Q,2},\  a\in F^\times.$$
\vskip 5pt

An examination of the root datum gives:
\[  \overline{\GSp}_{2r}^{\vee} = \begin{cases}
\GSp_{2r}(\C), \text{  if $r$ is odd;} \\
{\rm PGSp}_{2r}(\C)  \times \GL_1(\C),  \text{  if $r$ is even.}  \end{cases} \]
This explains the difference in the representation theory of $\overline{\GSp}_{2r}$ for even and odd $r$  observed in the work of Szpruch \cite{Sz3}.

\vskip 5pt

\vskip 10pt

\section{\bf Questions}

In this final section, we highlight a few questions which we feel are important to carry the program forward:
\vskip 5pt

\begin{itemize}
\item[(a)] (Real Groups)  When $F = \R$, Harish-Chandra's classification of discrete series representation works equally well for BD covering groups. Might one be able to formulate this classification in terms of L-parameters in the spirit of this paper, analogous to what Langlands accomplished for linear real reductive groups?  
Moreover, can the results of the recent papers \cite{AH, ABPTV} be formulated in the framework of this paper?  
 \vskip 5pt

\item[(b)] (Hecke algebra isomorphisms) When$F$ is $p$-adic and  $\G$ is simply-connected, Savin has studied the Iwahori-Hecke algebra of covers of $G$ and established Hecke algebra isomorphisms with those of linear reductive groups.  One would expect to show such results in the generality of this paper. Namely, for any BD covering group $\overline{G}$  
such that the hyperspecial maximal compact subgroup $K$ splits in the covering,  one may consider its Iwahori Kecke algebra $\mathcal{H}(\overline{G}, I)$. 
The L-group of $\overline{G}$ is isomorphic to a semi-direct $\overline{G}^{\vee} \rtimes W_F$ and this determines a unique quasi-split linear reductive group group $G'$ over $F$  with this L-group. One might expect that there is a natural isomorphism (relative to a distinguished splitting) 
\[  \mathcal{H}(\overline{G},  I) \cong \mathcal{H}(G', I'). \]
\vskip 5pt

\item[(c)] (Supercuspidals) Depth zero supercuspidal representations of BD covering groups have been studied by Howard-Weissman \cite{HW}. 
One would expect that the construction of supercuspidal representations of J.K. Yu and S. Stevens for linear reductive groups can be extended to BD covering groups. This is an ongoing investigation of the first author with J.L.Kim. One expects that Kim's proof of the exhaustion of Yu's construction can be extended to give exhaustion in the covering case. For this, a better understanding of the Bruhat-Tits theory of BD covering groups is probably needed; a preliminary study has been conducted by Weissman.
\vskip 5pt

\item[(d)] (Harmonic Analysis) For covering groups, any conjugacy-invariant function or distribution is necessarily supported on the subset of "good" or "relevant" elements. These are elements $g \in \overline{G}$ such that $g$ is not conjugate to $g \cdot \epsilon$ for any $\epsilon \in \mu_n$ (for a degree $n$ cover). It is clear that 
to have a better understanding of invariant harmonic analysis, a better understanding of such "good" or "relevant" elements is necessary.  One might ask if the BD structure theory is robust enough to give one a classification of such elements.
\vskip 5pt

\item[(e)] (Automorphic L-functions) We have defined the global partial automoprhic L-function associated to an automorphic representation of a BD covering group. One would like to have a definition of the local L-factors at the remaining set of places so as to obtain a complete L-function.  One would also like to show that the partial  L-function has a meromorphic continuation, and the complete L-function satisfies a standard functional equation relative to $s \leftrightarrow 1-s$. For this, one may ask if Langlands-Shahidi theory can be extended to the covering case, but this is not clear since uniqueness of Whittaker models is false in general for covering groups. The only such success is the thesis work \cite{Sz1, Sz2} of D. Szpruch where the Langlands-Shahidi theory was extended to  the group $\Mp_{2n}$. One might also ask if the myriad of Rankin-Selberg integrals for various L-functions of linear groups have counterparts in the covering case.

\vskip 5pt

\item[(f)] (Functoriality) More generally, one would like to show that this class of automorphic L-functions from BD covering groups belong to the class of automorphic L-functions of linear reductive group.  In the context of (b), one might expect that there is a functorial  transfer of automorphic representations from $\overline{G}$ to $G'$ which respects (partial) automorphic L-functions. One might imagine comparing the trace formula for these two groups. For this, we note that  the work of Wenwei Li \cite{L2, L4, L5, L6} has carried the theory of the trace formula for covering groups to the point where one has the invariant trace formula. 
\vskip 5pt

\item[(g)] (Endoscopy)  The next step in the theory of the trace formula for covering groups is  undoubtedly the stable trace formula. For this, one needs to develop the theory of endoscopy for covering groups. This includes the definition of stable conjugation, the definition of endoscopic groups, the definition of correspondence of stable classes between a cove rig group and its endoscopic groups and the definition of the transfer factors. Since the theory of endoscopy for linear reductive groups is essentially of arithmetic origin and content, one might expect that a nice theory exists for the BD covering groups since these are of algebraic origin. The only covering group for which a theory of endoscopy exists is the group $\Mp_{2n}$, where the theory is due to Adams \cite{A1}, Renard \cite{R1, R2} and W.W. Li \cite{L1, L3}. The recent preprint \cite{L7} of W.W. Li has taken the first step towards the stabilisation of the trace formula for $\Mp_{2n}$.
 
\vskip 5pt

\item[(h)]  (Automorphic Discrete Spectrum) Naturally, one hopes to have an analog of the Arthur's conjecture for BD covering groups, including an analog of the Arthur multiplicity formula. This will very much depend on the shape of the  theory of endoscopy. The only BD covering group for which a precise conjecture exists is $\Mp_{2n}$.

\vskip 5pt

\item[(i)]  (General Covering Tori)  The various questions highlighted above are already highly non-trivial and interesting when $\G = \T$ is a (not necessarily split) torus.  The ongoing work of Weissman \cite{W5} and Hiraga-Ikeda aim to understand this case completely but many mysteries remain. 
\vskip 5pt

\item[(j)] (Applications) 
The impetus for a program naturally depends on its potential applications. The motivation for our investigations is simply in the naive hope of including the representation theory and automorphic forms of BD covering groups in the framework of the Langlands philosophy.  This is reasonable enough (if naive) for the point of view of representation theory.  But what about form the point of view of number theory?  Automorphic forms of covering groups have traditionally found applications in analytic number theory, such as in the work of Bump-Friedberg-Hoffstein \cite{BFH}. It is reasonable to demand concrete arithmetic applications of this potential theory. The only thought we have to offer is perhaps in various branching or period problems, such as the arithmetic information contained in Fourier coefficients or in the analogs of the Gross-Prasad conjecture.

\vskip 5pt

\item[(k)] (Geometric Counterpart) 
The definition of the dual group of a BD cover first appeared in the context of the Geometric Langlands Program, throughh the work of Finkelberg-Lysenko \cite{FL}. One might expect the geometric theory to offer more evidence for this program.   From the geometric side, iquantum groups seems to play an important role in the theory.  This is also reflected to some extent in the work of Brubaker-Bump-Freidberg-et-al  \cite{BBF} on the Whittaker-Fourier  coefficients of metaplectic Eisenstein series and the metaplectic Casselman-Shalika formula. However, quantum groups are conspicuously missing from the framework developed in this article, and one may wonder if and how they should be incorporated.
\vskip 5pt

\item[(l)] (Function Fields) 
On the other hand, one can consider classical function fields (of curves over finite fields) and ask whether V. Lafforgue's recent construction \cite{Laf} of the global Langlands correspondence for arbitrary linear reductive groups could be extended to the case of BD covering groups: this seems a very tantalizing problem whose resolution should shed much light on the Langlands-Weissman program. 
\end{itemize}


\vskip 15pt

\begin{thebibliography}{9999999}


 \bibitem[A1]{A1}  J. Adams, {\em Lifting of characters on orthogonal and metaplectic groups}, Duke Math. J. 92 (1998), no. 1, 129�178.

\bibitem[A2]{A2}  J. Adams, {\em Characters of covering groups of SL(n)}, J. Inst. Math. Jussieu 2 (2003), no. 1, 1�21. 

\bibitem[ABPTV]{ABPTV}    J. Adams, D. Barbasch, A. Paul, P. Trapa and David Vogan, {\em Shimura correspondences for split real groups}, Journal of the American Math Society 20 (2007), 701-751.

\bibitem[AH]{AH} J. Adams and R. Herb, {\em Lifting of characters for nonlinear simply laced groups}, 
Represent. Theory 14 (2010), 70�147 


\bibitem[BJ]{BJ} D. Ban and C. Jantzen,
\emph{The Langlands quotient theorem for finite central extensions of p-adic groups}, preprint,
2013.




\bibitem[BBF]{BBF}
B. Brubaker, D. Bump and S. Friedberg,
\emph{Weyl group multiple Dirichlet series, Eisenstein series and crystal bases},
Ann. of Math. 173 (2011), No. 2, 1081–1120.


\bibitem[BD]{BD}
J.-L. Brylinski and P. Deligne,
\emph{Central extensions of reductive groups by $K_2$},
Publ. Math. Inst. Hautes $\acute{\text{E}}$tudes Sci., \textbf{94} 5-85, 2001.

\bibitem[BFH]{BFH}  D. Bump, S. Friedberg and J. Hoffstein, {\em Eisenstein series on the metaplectic group and nonvanishing theorems for automorphic L-functions and their derivatives}, Ann. of Math. (2) 131 (1990), no. 1, 53-127. 


\bibitem[Bl]{Bl}
S. Bloch,
\emph{Some formulas pertaining to the $K$-theory of commutative groupschemes},
Journal of Algebra, (53) 304-326, 1978.



\bibitem[D]{D} P. Deligne, {\em Extensions centrales de groupes alg\'ebriques simplement connexes et cohomologie galoisienne}, 
Inst. Hautes \'Etudes Sci. Publ. Math. No. 84 (1996), 35-89.

\bibitem[De]{De} V. Deodhar, {\em On central extensions of rational points of algebraic groups}, 
 Amer. J. Math. 100 (1978), no. 2, 303-386. 

  \bibitem[FL]{FL}  M. Finkelberg and S. Lysenlo, {\em Twisted geometric Satake equivalence}, 
J. Inst. Math. Jussieu 9 (2010), no. 4, 719-739.

\bibitem[F]{F} Y. Flicker, {\em Automorphic forms on covering groups of ${\rm GL}(2)$}, Invent. Math. 57 (1980), no. 2, 119-182.


\bibitem[FK]{FK}
Y. Flicker and D. Kazhdan,
\emph{Metaplectic correspondence},
I.H.E.S., (64) 53-110, 1986.








\bibitem[Ga]{Ga} F. Gao, {\em The Gindikin-Karpelevich formula and constant terms of Eisenstein series for Brylinski-Deligne extensions}, PhD thesis, National Univ. of Singapore.

\bibitem[HW]{HW} T. K. Howard and M. H, Weissman, {\em 
Depth-zero representations of nonlinear covers of p-adic groups}, Int. Math. Res. Not. IMRN 2009, no. 21, 3979-3995.


\bibitem[KP1]{KP1}
D. A. Kazhdan and S. J. Patterson,
\emph{Metaplectic forms},
I.H.E.S. Publ. Math. (59) 35-142, 1984.

\bibitem[KP2]{KP2}
D. A. Kazhdan and S. J. Patterson,
\emph{Towards a generalized Shimura correspondence},
Adv. in Math. 60(2) 161-234, 1986.



\bibitem[K]{K}
S. Kudla,
\emph{Notes on the local theta correspondence},
available at www.math.toronto.edu/~skudla.

\bibitem[Ki]{Ki}  J. L. Kim, {\em Supercuspidal representations: an exhaustion theorem}, J. Amer. Math. Soc. 20 (2007), no. 2, 273-320 

\bibitem[Ku]{Ku} T. Kubota, {\em On automorphic functions and the reciprocity law in a number field}, 
Lectures in Mathematics, Department of Mathematics, Kyoto University, No. 2 Kinokuniya Book-Store Co., Ltd., Tokyo 1969 iii+65 pp. 


\bibitem[La]{La}
 R. Langlands, \emph{Euler products},
 Yale Univ. Press, 1971.

\bibitem[Laf]{Laf} V. Lafforgue, {\em Chtoucas pour les groupes r\'eductifs et param\'etrisation de Langlands globale}, available at  arXiv:1209.5352.

\bibitem[L1]{L1} W.W. Li, {\em Transfert d'int\'egrales orbitales pour le groupe m\'etaplectique}, Compositio Mathematica 147 (2011), 524-590.

\bibitem[L2]{L2} W.W. Li, {\em La formule des traces pour les rev\^etements de groupes r\'eductifs connexes. I. Le d\'eveloppement 
g\'eom\'etrique fin}, Journal f�r die reine und angewandte Mathematik. Issue 686 (2014), pp.37-109. 

\bibitem[L3]{L3} W.W. Li, {\em Le lemme fondamental pond\'er\'e pour le groupe m\'etaplectique}, 
Canadian Journal of Mathematics, Vol 64 (3), 2012, pp.497-543.

\bibitem[L4]{L4} W.W. Li, {\em La formule des traces pour les rev\^etements de groupes r\'eductifs connexes. II. 
Analyse harmonique locale}, Annales scientifiques de l'ENS 45, fascicule 5 (2012).

\bibitem[L5]{L5} W.W. Li, {\em La formule des traces pour les rev\^etements de groupes r\'eductifs connexes. III. 
Le d\'eveloppement spectral fin}, Mathematsiche Annalen, Vol 356 (3), 2013, pp.1029-1064.

\bibitem[L6]{L6} W.W. Li, {\em La formule des traces pour les rev\^etements de groupes r\'eductifs connexes. IV. 
Distributions invariantes,} to appear in Annales de l'Institut Fourier.

\bibitem[L7]{L7} W.W. Li, {\em La formule des traces stable pour le groupe m\'etaplectique: les termes elliptiques}, 
submitted (arXiv:1307.1032). 

 


\bibitem[Ma]{Ma}
H. Matsumoto,
\emph{Sur les sous-groupes arithm\'etiques des groupes semi-simples d\'eploy\'es},
Ann. Sci. \'Ecole Norm. Sup. (4) 1-62, 1969.


\bibitem[Mc1]{Mc1}
P. McNamara, 
{\em Metaplectic Whittaker functions and crystal bases},  Duke Math. J. 156 (2011), no. 1, 1-31.

\bibitem[Mc2]{Mc2}
P. McNamara,
\emph{Principal series representations of metaplectic groups over local fields}
in \emph{Multiple Dirichlet series, $L$-functions and automorphic forms},
Birkhauser, 2012, 299-328.

\bibitem[Mc3]{Mc3}
P. McNamara, {\em The Metaplectic Casselman-Shalika Formula}, preprint, arXiv:1103.4653.

\bibitem[Mo]{Mo}
C. C. Moore,
\emph{Group extensions of $p$-adic and adelic linear groups},
I.H.E.S. Publ. Math. (35) 157-222, 1968.

\bibitem[MW]{MW}
 C. Moeglin and J.-L. Waldspurger,
 \emph{Spectral decomposition and Eisenstein series},
 Cambridge University Press, 1995.

\bibitem[PR1]{PR1} G. Prasad and M.S. Raghunathan, {\em Topological central extensions of semisimple groups over local fields},
 Ann. of Math. (2) 119 (1984), no. 1, 143�201.

\bibitem[PR2]{PR2} G. Prasad and M.S. Raghunathan, {\em Topological central extensions of semisimple groups over local fields. II}, 
 Ann. of Math. (2) 119 (1984), no. 2, 203�268. 

\bibitem[PR3]{PR3} G. Prasad and M.S. Raghunathan,
 {\em Topological central extensions of $\SL_1(D)$},  Invent. Math. 92 (1988), no. 3, 645�689.


\bibitem[Ra]{Ra}
R. Rao,
\emph{On some explicit formulas in the theory of Weil representation},
Pacific Journal of Math. 157 (1993), No. 2, 335-371.

\bibitem[Re]{Re} R.C. Reich, \emph{Twisted geometric Satake equivalence via gerbes on the factorizable Grassmannian}, preprint.

\bibitem[R1]{R1} D. Renard, {\em Transfert d'int\'{e}grales orbitales entre $\Mp(2n,\R)$ et $\SO(n+1,n)$},  Duke Math. J. 95 (1998), no. 2, 425-450.

\bibitem[R2]{R2} D. Renard, {\em Endoscopy for $\Mp(2n,\R)$},  Amer. J. Math. 121 (1999), no. 6, 1215-1243.

 
\bibitem[Sa]{Sa}
G. Savin,
\emph{On unramified representations of covering groups},
J. Reine Angew. Math. 566 (2004), 111-134.






\bibitem[S]{S}
R. Steinberg,
\emph{G\'en\'erateurs, relations et rev\^etements de groupes alg\'ebriques}, in
\emph{Colloq. Th\'eorie des Groupes Alg\'ebriques} (Bruxelles, 1962), Louvain, 1962.

\bibitem[St]{St} S. Stevens, {\em The supercuspidal representations of p-adic classical groups}, Invent. Math. 172 (2008), no. 2, 289-352. 

\bibitem[Sz1]{Sz1}
D. Szpruch,
\emph{The Langlands-Shahidi method for the metaplectic group and applications},
thesis (Tel Aviv University), available at arXiv: 1004.3516v1.

\bibitem[Sz2]{Sz2} D. Szpruch, {\em  Some irreducibility theorems of parabolic induction on
 the metaplectic group via the Langlands-Shahidi method}, Israel J. Math. 195 (2013), no. 2, 897-971.

\bibitem[Sz3]{Sz3} D. Szpruch, {\em  Some results in the theory of genuine representations of the metaplectic double cover 
of $\GSp_{2n}(F)$ over p-adic fields},  J. Algebra 388 (2013), 160-193. 

\bibitem[W1]{W1}
M. Weissman,
\emph{Metaplectic tori over local fields},
Pacific J. Math., \textbf{241}(1) 169-200, 2009.

 
\bibitem[W2]{W2}
---------------,
\emph{Managing metaplectiphobia: covering $p$-adic groups} in
\emph{Harmonic analysis on reductive, $p$-adic groups}, 2011, pg. 237-277.

 
\bibitem[W3]{W3}
---------------,
\emph{Split metaplectic groups and their $L$-groups},
to appear in the Crelle's journal.


 \bibitem[W4]{W4}
---------------,
A letter to P. Deligne, 2012.


 
\bibitem[W5]{W5}
---------------,
\emph{Covers of tori over local and global fields},
available at http://arxiv.org/abs/1406.3738

\bibitem[W6]{W6}
---------------,
\emph{Covering groups and their integral models},
available at http://arxiv.org/abs/1405.4625


\bibitem[W7]{W7} M. H. Weissman, {\em Coverig the Langlands program}, in progress. 

\bibitem[Yu]{Yu} J. K. Yu,  {\em Construction of tame supercuspidal representations}, J. Amer. Math. Soc. 14 (2001), no. 3, 579-622


\end{thebibliography}
\end{document}